\documentclass[11pt]{article}

\usepackage{amssymb,amsmath}
\usepackage{epsfig}
\usepackage{graphicx}
\usepackage{psfrag}

\newtheorem{theorem}{Theorem}[section]

\newtheorem{proposition}[theorem]{Proposition}
\newtheorem{corollary}[theorem]{Corollary}

\newenvironment{proof}[1][Proof]{\begin{trivlist}
\item[\hskip \labelsep {\bfseries #1}]}{\end{trivlist}}

\newcommand{\qed}{\nobreak \ifvmode \relax \else
      \ifdim\lastskip<1.5em \hskip-\lastskip
      \hskip1.5em plus0em minus0.5em \fi \nobreak
      \vrule height0.75em width0.5em depth0.25em\fi}

\title{A slice genus lower bound from $sl(n)$ Khovanov-Rozansky homology}
\author{Andrew Lobb}
\date{}

\begin{document}

\maketitle

\tableofcontents

\section{Introduction}

\subsection{$sl(2)$ Khovanov homology}

\begin{figure}
\centerline{
{
\psfrag{a}{$a$}
\psfrag{-a}{$-a^{-1}$}
\psfrag{=b}{$=b$}
\includegraphics[height=1in,width=4in]{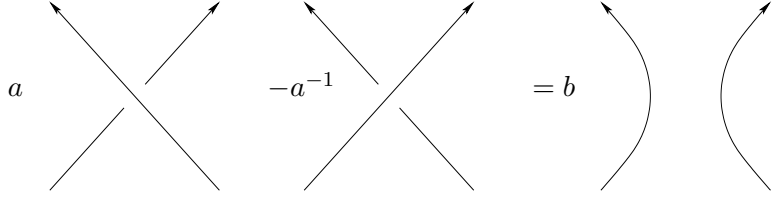}
}}
\caption{HOMFLY skein relation}
\label{HOMFLY}
\end{figure}

Given the data of an oriented link diagram $D$, one can compute the HOMFLY polynomial $P(D)$ (a polynomial over $\mathbb{Z}$ in the variables $a^{\pm 1}$ and $b^{\pm 1}$) using the local skein relations in Figure \ref{HOMFLY} (up to an arbitrary choice for the HOMFLY polynomial $P(U)$ of the unknot $U$).

%% Insert HOMFLY skein diagram here

This polynomial is invariant under the oriented Reidemeister moves and hence defines an invariant of oriented links.  If we specialize for $n \geq 2$ by defining

\[ P_n (q) = P(a=q^n, b = q + q^{-1}) \]

\noindent we obtain the $sl(n)$ quantum polynomial of the link given by the diagram $D$.  We will be using the normalization $P_n(U) = \frac{q^n - q^{-n}}{q - q^{-1}}$.  The polynomial $P_2(q)$ is known as the Jones polynomial.

In \cite{Khov1}, Khovanov associated to a link diagram $D$ a bigraded chain complex $CKh_2^{i,j}$ with differential

\[ d : CKh_2^{i,j}(D) \rightarrow CKh_2^{i+1,j}(D) \]

\noindent ($i$ is called the homological grading, $j$ is called the quantum grading).  To diagrams $D$ and $D'$ differing by a single Reidemeister move, Khovanov gave a quantum-degree $0$ chain homotopy equivalence between $CKh_2(D)$ and $CKh_2(D')$, thus showing that the homology groups $HKh_2^{i,j}(D)$ are knot invariants.  Furthermore these homology groups provided a \emph{categorification} of the Jones polynomial $P_2(D)$, by which is meant

\[ P_2(q) = \sum_{i,j} (-1)^i q^j dim(HKh_2^{i,j}) \]

A powerful facet of the $HKh_2$ theory was conjectured by Khovanov (and later proved by Jacobsson \cite{Jacobsson}), namely that the homology theory should be functorial for link cobordisms up to sign.  More explicitly, suppose we start with a smooth embedding of a surface $\Sigma$ with boundary 

\[ \Sigma \hookrightarrow [0,1] \times \mathbb{R}^3 \]

\[ \partial \Sigma = (L_0 \hookrightarrow \{0\} \times \mathbb{R}^3) \coprod (L_1 \hookrightarrow \{ 1 \} \times \mathbb{R}^3) \]

\begin{figure}
\centerline{
{
\psfrag{0handle}{$0$-handle}
\psfrag{1handle}{$1$-handle}
\psfrag{2handle}{$2$-handle}
\includegraphics[height=2.7in,width=3.2in]{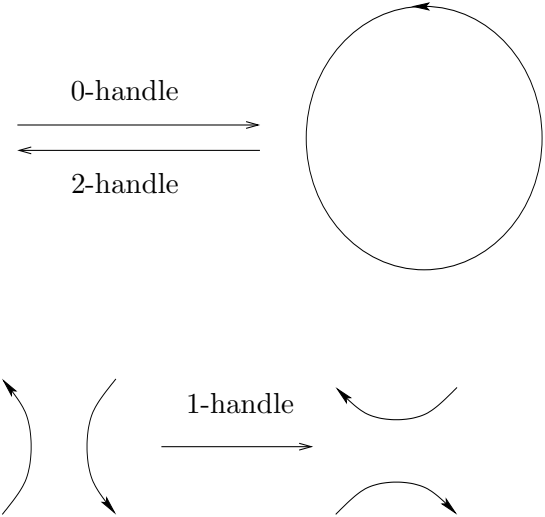}
}}
\caption{Morse moves for knot cobordism}
\label{morsemoves}
\end{figure}

\noindent Choose link diagrams $D_0$, $D_1$ of the links $L_0$, $L_1$.  Next we take a representation of the surface $\Sigma$ as a product of \emph{elementary cobordisms}.  Elementary cobordisms consist of before-and-after link diagrams where we have made one local change in the before diagram to get to the after diagram.  The local changes that we are allowed to make consist of each of the Reidemeister moves and the \emph{Morse moves} which correspond to adding a handle to the surface as shown in Figure \ref{morsemoves}.

%% Here show the Morse moves.

To diagrams $D$ and $D'$ differing by one elementary cobordism Khovanov associated a chain map $CKh_2(D) \rightarrow CKh_2(D')$, inducing a map on homology.  For the Reidemeister moves, these maps are just the chain homotopy equivalences referred to earlier.  For the $0$-handle and $2$-handle moves these maps are graded of quantum-degree $+1$ and for the $1$-handle move it is graded of degree $-1$.  By composing the maps corresponding to these elementary cobordisms we get a map graded of quantum-degree $\chi(\Sigma)$

\[ HKh_2(\Sigma) : HKh_2(L_0) \rightarrow HKh_2(L_1) \]

\noindent This map (up to multiplication by $-1$) is an invariant of $\Sigma$, as the notation suggests.  In other words, whichever decomposition of $\Sigma$ into elementary cobordisms is chosen, the induced map $HKh(\Sigma)$ stays the same up to sign.

\subsection{Perturbed $sl(2)$ theory}

In 2002, Lee gave a description of a perturbation $HKh'$ of the Khovanov homology $HKh_2$ \cite{Lee}.  This homology theory is no longer graded in the quantum grading direction, but filtered: $\ldots \subseteq \mathcal{F}^{j+1} HKh'^{i} \subseteq \mathcal{F}^j HKh'^{i} \subseteq \ldots$.  $HKh'(L)$ has a particularly simple form: when we take the base ring to be $\mathbb{C}$ it consists of $2^l$ copies of $\mathbb{C}$ where $l$ is the number of components of the link $L$.  Lee gave explicit chain representatives of these generators defined from any diagram presentation of $L$.  As a formal consequence of the properties of filtered chain complexes, there is a spectral sequence with $E_2$ page $HKh_2^{i,j}$ converging to $E_{\infty}$ page  $Gr^j HKh'^i := \mathcal{F}^j HKh'^i / \mathcal{F}^{j+1} HKh'^i$ the associated graded vector space of $HKh'$ with respect to the filtration $\mathcal{F}$.

Analysis of the behaviour of $HKh'$ under link cobordism was carried out by Rasmussen \cite{Ras}.  In this remarkable paper, Rasmussen showed that Lee's generators provide a lower bound for the slice genus of a knot.  Let us digress to define the slice genus.

Given a knot $K \hookrightarrow \mathbb{R}^3$, a classical knot invariant is the \emph{genus} $g \geq 0$ of the knot.  That is, the minimal genus of the surfaces-with-boundary $\Sigma$ for which there exists an embedding $\Sigma \hookrightarrow \mathbb{R}^3$ with $\partial \Sigma = K$.  

%%These $\Sigma$ are known as Seifert surfaces.  There exist classical invariants giving lower bounds on $g(K)$, and recently Heegaard-Floer knot homology has given a way to compute the genus precisely, even if the computations are impossible to do in practice at the moment for many knots.

If we believe that manifolds with boundary should rightfully live within other manifolds with boundary, we are motivated to consider surfaces $\Sigma$ smoothly embedded in $(-\infty,0] \times \mathbb{R}^3$ with $K = \partial \Sigma \hookrightarrow \{ 0 \} \times \mathbb{R}^3$.  The minimal genus of such surfaces we call the \emph{slice genus} $g_*(K)$ of the knot $K$.

By removing a neighbourhood of a point, such a surface provides a knot cobordism between $K$ and the unknot $U$.  Rasmussen showed that the map associated to any presentation (as a composition of elementary cobordisms) of a connected knot cobordism $\Sigma$ between two ($1$-component) knots $K_0$ and $K_1$ preserves the generators of Lee's homology $HKh'$.  This map is filtered of quantum-degree $\chi(\Sigma)$.

\[ \mathcal{F}^j HKh'^i \rightarrow \mathcal{F}^{j + \chi(\Sigma)} HKh'^i \]

The homology of the unknot is computable as $Gr^1 HKh'^0 = \mathbb{Q}$ and $Gr^{-1} HKh'^0 = \mathbb{Q}$ with no homology in any other bigrading.  Hence as an immediate corollary, $2g_* (K)$ is bounded below by one less than the highest degree in which $HKh'(K)$ is non-zero.

To compute $HKh'$ directly is difficult, but it can be carried out for positive knots (those knots that admit a diagram with only positive crossings).  In this case we get a tight bound on the slice genus since we can explicitly construct a Seifert surface $\Sigma$ of the same genus as our lower bound, and then just push this surface into $(-\infty,0] \times \mathbb{R}^3$ to get a slice surface.  In particular we can apply this to torus knots $K_{p,q}$ to get a combinatorial proof of Milnor's conjecture on the value of $g_*(K_{p,q})$.

\subsection{$sl(n)$ Khovanov homology}

In \cite{KhovRoz1}, Khovanov and Rozansky describe bigraded homology theories $HKh_n$ for $n \geq 3$ which categorify the $sl(n)$ polynomials $P_n$.  Again, given a knot diagram $D$ Khovanov and Rozansky associate to it a chain complex $CKh_n(D)$ (in a much more complicated way than for the case $n=2$) with differentials

\[ d_n : CKh_n^{i,j}(D) \rightarrow CKh_n^{i+1,j}(D) \]

\noindent And also as before, if $D$ and $D'$ differ by a Reidemeister move, they give a chain homotopy equivalence (graded of quantum-degree $0$) between $CKh_n(D)$ and $CKh_n(D')$, thus showing that $HKh_n^{i,j}(D)$ is an invariant of the knot and not just the diagram $D$.  This homology theory is also functorial for knot cobordisms.

Gornik has carried out for $HKh_n$ a direct analogue of Lee's work on $HKh_2$ \cite{Gornik}, describing a perturbation with filtered homology $HKh'_n$, which can be interpreted as the $E_{\infty}$ page of a spectral sequence with $E_2$ page $HKh_n$.  Given an $l$-component link diagram $D$, Gornik gives an explicit description of the $n^l$ generators of $HKh'_n(D)$.

\subsection{Statement of results}

Our intention in this paper is to do for Gornik's work something of what Rasmussen did for that of Lee's.  We start by generalizing Gornik's $HKh'_n$ to a theory that we shall denote $HKh_w$ where $w \in \mathbb{C}[x]$ is a monic polynomial of degree $n+1$ such that $\partial_x w = dw/dx$ is a product of distinct linear factors.  In this set up,

\[ HKh'_n = HKh_{x^{n+1} + (n+1)\beta^n x} \]

\noindent for $\beta \in \mathbb{C} - \{ 0 \}$.  $HKh_w$ is filtered in the quantum direction

\[ \ldots \subseteq \mathcal{F}^j HKh_w^{i}\subseteq \mathcal{F}^{j+1} HKh_w^i \subseteq \ldots \]

Our procedure will essentially follow that of \cite{Ras}.  For diagrams $D$, $D'$ differing by a Reidemeister move we wish construct a map $\mathcal{F}^j HKh_w^i(D) \rightarrow \mathcal{F}^j HKh_w^i(D')$ which preserves our analogues of the generators in \cite{Gornik}.  If $D$, $D'$ differ by a $0$- or $2$-handle attachment we wish to give maps $\mathcal{F}^j HKh_w^i(D) \rightarrow \mathcal{F}^{j - n +1} HKh_w^i(D')$ and if $D$, $D'$ differ by a $1$-handle we wish to give maps $\mathcal{F}^j HKh_w^i(D) \rightarrow \mathcal{F}^{j+n-1} HKh_w^i$.  We would aim to do this so that we could then compute that any representation of a connected knot cobordism $\Sigma : D \rightarrow D'$ as a product of elementary knot cobordisms preserves the generators of $HKh_w(D)$.

Due to the complexity of the Reidemeister $III$ move, we do not quite fulfill all of our wishes, however a topological argument (in which we only allow certain products of elementary cobordisms) is enough for us to get the desired analogue of the slice genus bound in \cite{Ras}.

We write $HKh^{i,j}_w(D) = \mathcal{F}^jHKh^i_w / \mathcal{F}^{j-1}HKh^i_w$ for the associated graded vector space to the quantum filtration.

Our main theorem is

\begin{theorem}
\label{maintheorem}

Let $D$ be a knot diagram of a knot $K$.  If $HKh_w^{0,j}(D) \not= 0$  and $HKh_w^{0,j'}(D) = 0$ for all $j'>j$ then

\[ (n-1)(2g^*(K) - 1) \geq - j \]

\end{theorem}

\noindent As a corollary we obtain another proof of Milnor's conjecture on the slice genus of torus knots.

It is worth noting here that this paper does not contain a proof that

\[ 1. \, HKh_w(D) = HKh_w(D') \]

\noindent for two knot diagrams differing by a single Reidemeister move, even for Gornik's $w$, nor that

\[ 2. \, HKh_w(D) = HKh_{w'}(D) \]

\noindent for two different degree $n+1$ polynomials $w$, $w'$.  However we expect both of these statements to be true (for example because there is a spectral sequence with $E_2$ pages independent of choice of $w$ or $D$ converging to $HKh_w(D)$).

In the final stages of the preparation of this paper a preprint appeared on the arXiv \cite{Wu}, which contains, amongst others, results similar to those in this paper.

\subsection{The Khovanov cube}

In this section we shall outline the construction of $HKh_w$.  This is similar to the construction of the original Khovanov-Rozansky homology $HKh_n$ - indeed, replacing any occurence of the word ``filtered'' by the word ``graded'' in this section will give an outline of the construction of $HKh_n$.

\begin{figure}
\centerline{
{
\psfrag{1}{$1$}
\psfrag{0}{$0$}
\psfrag{3}{$3$}
\psfrag{4}{$4$}
\psfrag{=}{$=$}
\psfrag{eta0}{$\eta_0$}
\psfrag{eta1}{$\eta_1$}
\psfrag{chi0}{$\chi_0$}
\psfrag{chi1}{$\chi_1$}
\includegraphics[height=2.5in,width=2.5in]{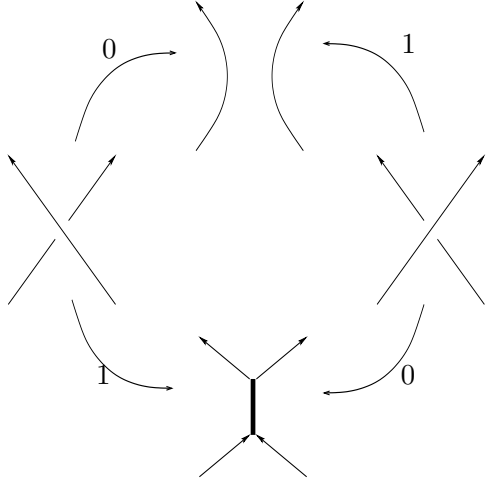}
}}
\caption{Forming resolutions}
\label{howtoresolve}
\end{figure}

To define $HKh_w$ of an oriented link diagram $D$ with $m$ crossings we start by forming the $2^m$ possible \emph{resolutions} of $D$ and decorating each vertex of the cube $[0,1]^n$ with one of these resolutions.  We form the resolutions as shown in Figure \ref{howtoresolve}.

Choosing an (arbitrary) ordering of the crossings allows us to associate each resolution $\Gamma$ of $D$ with one of the corners of the cube.  For a corner $v$ call the associated resolution $\Gamma_v$.  In the following subsections we shall see how to associate to $\Gamma_v$ a filtered vector space $H(\Gamma_v)$.  We shall often refer to the filtration as the \emph{quantum} filtration.

If two corners of the cube $v$, $v'$ are connected by an edge $e$ then we see that their associated resolutions $\Gamma_v$ and $\Gamma'_{v'}$ differ only in the resolution of a single crossing of $D$.  Suppose that $v$ is at the $0$-coordinate of the edge and $v'$ is at the $1$-coordinate of the edge.  In the following subsections we shall define a map $\Phi_e : H(\Gamma_v) \rightarrow H(\Gamma_{v'})$ of quantum filtered degree $1$.  If $e_1$, $e_2$, $e_3$, $e_4$ are edges bounding a face of the cube then $\Phi_{e_1}\Phi_{e_2} = \Phi_{e_3}\Phi_{e_4}$ if the edges are ordered so the composition makes sense.

\begin{figure}
\centerline{
{
\psfrag{plus}{$+1$}
\psfrag{minus}{$-1$}
\psfrag{3}{$3$}
\psfrag{4}{$4$}
\psfrag{=}{$=$}
\psfrag{eta0}{$\eta_0$}
\psfrag{eta1}{$\eta_1$}
\psfrag{chi0}{$\chi_0$}
\psfrag{chi1}{$\chi_1$}
\includegraphics[height=1.4in,width=3in]{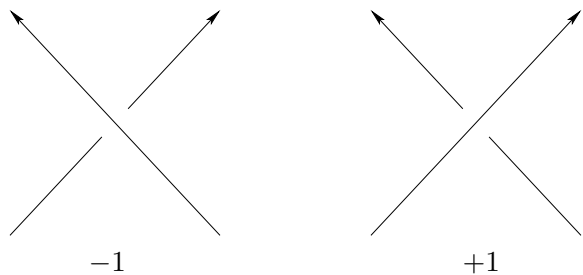}
}}
\caption{Definition of writhe}
\label{writhe}
\end{figure}

For a vertex $v$ write the sum of its coordinates as $s(v)$ (this will be an integer between $0$ and $n$).  We write $W$ for the \emph{writhe} of the knot diagram $D$: the signed number of crossings of $D$ as in Figure \ref{writhe}.  We define the chain groups of a chain complex by

\[ C^k = \sum_{s(v) = k + (1/2)(n+W)} H(\Gamma_v) \{(n-1)W - k\} \]

\noindent where by $\{ - \}$ we mean a shift in the quantum filtration.

The differentials of the chain complex are defined by

\[ d_k = \sum \pm \Phi_e : C^k \rightarrow C^{k+1} \]

\noindent where the sum is taken over all edges for which it makes sense.  The $\pm$ signs are chosen so that each face of the cube has exactly $1$ or $3$ of its edges decorated with a minus sign; the commutivity of the maps associated the edges of any face will then ensure that $d_{k+1}d_k = 0$.  The shift in the quantum filtrations of the $H(\Gamma_v)$'s ensures that the differentials are filtered of degree $0$.

Taking homology we obtain a vector space graded in the homological direction and filtered in the quantum direction.

\subsection{Introduction to matrix factorizations}

The definition of the $sl(n)$ Khovanov homology, standard or perturbed, makes use of the notion of a \emph{matrix factorization}.  Here we first define \emph{ungraded} matrix factorizations, in the next section we see how to associate \emph{filtered} matrix factorizations to trivalent graphs.

A polynomial $p(x) \in \mathbf{C}[x_1,...,x_n] = R$ may not admit a non-trivial factorization into polynomials.  However if $M^0 = M^1$ is a free module over $R$ then we may be able to find $R$-module maps $f_0, f_1$ 

\[ M^0 \stackrel{f_0}{\rightarrow} M^1 \stackrel{f_1}{\rightarrow} M^0 \]

\noindent such that $f_0 f_1 = f_1 f_0 = p(x)$.  Thus we would have factored $p(x)$ into a product of $R$-module maps - this is called a \emph{matrix factorization} of $p(x)$.

If $M$ and $\tilde{M}$ are both matrix factorizations of $p(x) \in \mathbf{C}[x_1,...,x_n] = R$ then a \emph{map of matrix factorizations} $G : M \rightarrow \tilde{M}$ is a pair of $R$-module maps

\[ G_0 : M^0 \rightarrow \tilde{M}^0, \, G_1 : M^1 \rightarrow \tilde{M}^1 \]

\noindent satisfying $\tilde{f}_0 G_0 = G_1 f_0$ and $\tilde{f}_1 G_1 = G_0 f_1$.  Note that this echoes the definition of a map of chain complexes.  Similarly, considering a matrix factorization as being something akin to a $2$-periodic chain complex, we can define the notion of a tensor product:

Suppose $M$ and $\tilde{M}$ are matrix factorizations of $p(x), q(x) \in R$ respectively.  Then we can define the matrix factorization $M \otimes \tilde{M}$ of $p(x) + q(x)$ by $(M \otimes \tilde{M})^0 = (M^0 \otimes \tilde{M}^0) \oplus (M^1 \otimes \tilde{M}^1)$, $(M \otimes \tilde{M})^1 = (M^1 \otimes \tilde{M}^0) \oplus (M^0 \otimes \tilde{M}^1)$ with maps

\[ \left( \begin{array}{cc} f_0 & \tilde{g}_1 \\ -\tilde{g}_0 & f_1 \end{array} \right) : (M \otimes \tilde{M})^0 \rightarrow (M \otimes \tilde{M})^1 \]

\[ \left( \begin{array}{cc} f_1 & -\tilde{g}_1 \\ \tilde{g}_0 & f_0 \end{array} \right) : (M \otimes \tilde{M})^1 \rightarrow (M \otimes \tilde{M})^0 \]

Since matrix factorizations have the look of $2$-periodic complexes (but with $d^2$ now being a polynomial, not necessarily zero), we can define the obvious notions of homotopic maps of matrix factorizations and homotopy equivalent matrix factorizations.  We will be working in the category of matrix factorizations and homotopy classes of maps.

\subsection{Trivalent graphs}

\begin{figure}
\centerline{
{
\psfrag{1}{$1$}
\psfrag{2}{$2$}
\psfrag{3}{$3$}
\psfrag{4}{$4$}
\psfrag{=}{$=$}
\psfrag{eta0}{$\eta_0$}
\psfrag{eta1}{$\eta_1$}
\psfrag{chi0}{$\chi_0$}
\psfrag{chi1}{$\chi_1$}
\includegraphics[height=2in,width=2in]{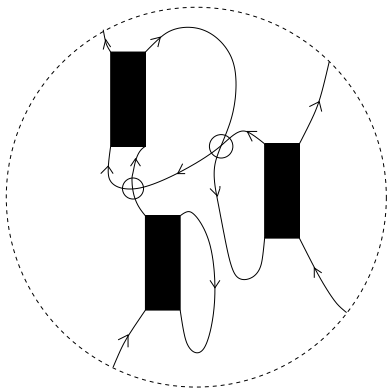}
}}
\caption{An example of a trivalent graph}
\label{exampletrivalent}
\end{figure}

Khovanov and Rozansky describe a way of associating a matrix factorization to any finite trivalent graph with thick edges and labelled boundary components.  An example of what we mean by a trivalent graph with thick edges is given in Figure \ref{exampletrivalent} - here we have indicated the boundary as a dotted circle which shall hereafter be omitted from our diagrams of trivalent graphs.

\begin{figure}
\centerline{
{
\psfrag{1}{$1$}
\psfrag{2}{$2$}
\psfrag{3}{$3$}
\psfrag{4}{$4$}
\psfrag{=}{$=$}
\psfrag{eta0}{$\eta_0$}
\psfrag{eta1}{$\eta_1$}
\psfrag{chi0}{$\chi_0$}
\psfrag{chi1}{$\chi_1$}
\includegraphics[height=1.5in,width=0.5in]{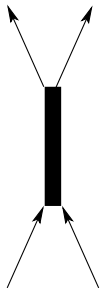}
}}
\caption{The neighbourhood of a thick edge}
\label{localthick}
\end{figure}

Note that the thin edges are oriented and at a thick edge the orientations of the incident thin edges look like Figure \ref{localthick}.  Each trivalent vertex has one thick edge and two thin edges incident to it.  Only thin edges are allowed to end on the boundary of the graph.  We allow closed thin loops.

Also note that we have included what appear to be crossings (which are circled) of thin edges.
But a trivalent graph need not come with an embedding into $\mathbb{R}^{2}$ or even $\mathbb{R}^{3}$,
and we are thinking of these thin edges as neither
intersecting nor giving rise to an ``overcrossing" or an
``undercrossing".

To associate a matrix factorization to such a trivalent graph will first involve choosing a \emph{potential}: a polynomial $w(x) \in \mathbb{C}[x]$ of degree $n+1$.  Khovanov's and Rozansky's original theory took $w = x^{n+1}$ and Gornik used the perturbation $w = x^{n+1} + (n+1) \beta^n x$ for $\beta \in \mathbb{C}-{0}$.  In the course of this paper we will be using a general polynomial $w$ of degree $n+1$, satisfying the condition that $\partial_x w  = dw/dx$ factors as a product of distinct linear factors.

Every variable $x_i$ that appears in the polynomial rings used to define matrix factorizations is taken to have quantum degree $2$ (we often use the word \emph{quantum} so as to distinguish from \emph{homological}).  We will be using \emph{filtered} matrix factorizations, which means each module $M$ is filtered $\ldots \subseteq \mathcal{F}^i M \subseteq \mathcal{F}^{i+1}M \subseteq \ldots$.  The filtration on $M$ (which will always be a free module over $R$) shall be induced by the grading on $R$ with possibly an overall shift.  To denote a shift of quantum degree $d$ in the filtration grading of a module $M$ we use the notation $ M \{ d \} $.  The two differentials in the matrix factorization should have filtered degree $n+1$ (note multiplication by $w(x)$ is filtered of degree $2(n+1)$).  Maps of matrix factorizations, unless stated otherwise, are of filtered degree $0$, homotopies are of filtered degree $-(n+1)$.

\begin{figure}
\centerline{
{
\psfrag{1}{$x_1$}
\psfrag{2}{$x_2$}
\psfrag{3}{$3$}
\psfrag{4}{$4$}
\psfrag{=}{$=$}
\psfrag{eta0}{$\eta_0$}
\psfrag{eta1}{$\eta_1$}
\psfrag{chi0}{$\chi_0$}
\psfrag{chi1}{$\chi_1$}
\includegraphics[height=1.4in,width=0.3in]{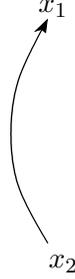}
}}
\caption{An oriented line segment with labelled boundary components}
\label{singlearc2}
\end{figure}

\begin{figure}
\centerline{
{
\psfrag{1}{$x_1$}
\psfrag{2}{$x_2$}
\psfrag{3}{$x_3$}
\psfrag{4}{$x_4$}
\psfrag{=}{$=$}
\psfrag{eta0}{$\eta_0$}
\psfrag{eta1}{$\eta_1$}
\psfrag{chi0}{$\chi_0$}
\psfrag{chi1}{$\chi_1$}
\includegraphics[height=1.4in,width=0.6in]{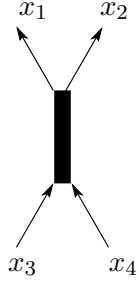}
}}
\caption{A thick edge with labelled boundary components}
\label{thickfactorization2}
\end{figure}

We define two fundamental factorizations - that of an oriented thin line segment (Figure \ref{singlearc2}) and that of a neighbourhood of a thick edge (Figure \ref{thickfactorization2}).  Then we can obtain the factorization corresponding to any trivalent graph by tensoring together these fundamental factorizations, identifying boundary variables where we have joined the fundamental factorizations together.

The fundamental factorization in Figure \ref{singlearc2} is the factorization of $w(x_1) - w(x_2)$ over
$\mathbb{C}[x_1,x_2]$ given by

\[ R \stackrel{\pi_{x_1 x_2}}{\longrightarrow} R\{1-n\} \stackrel{x_1 - x_2}{\longrightarrow}
R \]

\noindent where $M^0$ and $M^1$ are rank $1$ modules and $\pi_{x_1 x_2} =
\frac{w(x_1) - w(x_2)}{x_1 - x_2}$.  Remember the pair of curly brackets
denotes a grading shift so that each of our differentials is of
filtered degree $n+1$.

We now define the factorization in Figure \ref{thickfactorization2} over the ring $\mathbb{C}[x_1,x_2,x_3,x_4]$.  This will be a factorization of the polynomial $w(x_1) + w(x_2) - w(x_3) - w(x_4)$.

Consider the unique polynomial $p$ in two
variables such that $p(x+y, xy) = w(x) + w(y)$

Then
\begin{eqnarray*}
\lefteqn{w(x_1) + w(x_2) - w(x_3) - w(x_4)} \\
&& = p(x_1 + x_2, x_1 x_2) - p(x_3 +
x_4, x_3 x_4) \\
&& = p(x_1 + x_2, x_1 x_2) - p(x_3 + x_4, x_1 x_2) + p(x_3 + x_4, x_1
x_2) - p(x_3 + x_4, x_3 x_4) \\
&& = (x_1 + x_2 - x_3 - x_4)u_1 + (x_1 x_2 - x_3 x_4)u_2
\end{eqnarray*}
\noindent for some polynomials $u_1$, $u_2$.

We define the second fundamental factorization as the tensor product over $R(\Gamma)$ of the two factorizations

\[R\{ -1 \} \stackrel{u_{1}}{\longrightarrow} R\{ -n \} \stackrel{x_1 + x_2 -
  x_3 - x_4}{\longrightarrow} R\{ -1 \} \]

\[R \stackrel{u_{2}}{\longrightarrow} R\{3-n\} \stackrel{x_1 x_2 - x_3 x_4}
{\longrightarrow} R \]

\noindent Note that the differentials of the resulting factorization have filtered degree $n+1$.

Hence the factorization $C(\Gamma)$ (a tensor product of fundamental factorizations), associated to a trivalent graph $\Gamma$, also has differentials of this
filtered degree.  Let $R(\Gamma) = Q[{X_e \vert e \in {\rm{edge} \, \rm{endpoints}}}]$, then the
graph $\Gamma$ defines a
matrix factorization $C(\Gamma)$ over $R(\Gamma)$

\[ M^0 \stackrel{d^{0}}{\longrightarrow} M^1 \stackrel{d^{1}}{\longrightarrow} M^0 \]

\noindent where $M^0$ and $M^1$ are free modules over $R(\Gamma)$ and 

\begin{equation}
d^{1}d^{0}= d^{0}d^{1} = \mathop{\sum_{\rm{edge} \; \rm{endpoints} \; e}}_
{\rm{pointing} \; \rm{out \; of \; \Gamma}} w(X_{e}) - \mathop{\sum_{\rm{edge}\; \rm{endpoints} \; e}}_
{\rm{pointing} \; \rm{into \; \Gamma}} w(X_{e})
\label{inandoutburger}
\end{equation}

%\subsection{Choice of potential}

%We should be interested in potentials that correspond to equivariant
%cohomology rings of $\mathbb{CP}^{n-1}$.  We can specify an action of
%$S^1$ on $\mathbb{CP}^{n-1}$ by choosing a Lie group map $S^1 \rightarrow
%U(n)$.  Give this map by

%\[ \theta \mapsto diag(\lambda_1 \theta, \lambda_2 \theta, \ldots,
%\lambda_{n} \theta) \]

%The $S^1$-equivariant cohomology ring of $\mathbb{CP}^{n-1}$ with this
%action of $S^1$ is

%\[ H^*_{S^1}(\mathbb{CP}^{n-1}) = \mathbb{Q}[a,x]/w(a,x) \]

%\noindent where $w(a,x) = \prod_{i=1}^{n}(x + \lambda_i a)$.  (Here the
%cohomology ring $H^*(\mathbb{CP}^\infty) = H^*(BS^1) = \mathbb{Q}[a]$
%and both $a$ and $x$ have grading $2$).

%We shall be forgetting the homogenising variable $a$ and working just
%with the potential $w(x) = \prod_{i=1}^{n}(x + \lambda_i)$.  In
%fact, we will forget that the $\lambda_i$s are integers and let them
%take any values in $\mathbb{C}$, so long as $\partial_x w(x)$ is a
%product of distinct linear factors.  Thus
%we are thinking of our potential as a generic choice of tangent vector
%at $1$ along the maximal torus of $U(n)$.

%In ..., Khovanov gave a graded interpretation of Rasmussen's slice
%genus invariant.  Although we will give a lower bound on the slice
%genus arising from a filtered chain complex, we expect that the same
%bound can be recovered from a graded Khovanov-Rozansky homology by
%using the homogenised polynomial above.

\subsection{Method of attack}

We shall now give an overview of the organization of the remaining sections of this paper.

Section 2 contains definitions and theorems whose statements echo those to be found in \cite{Gornik} and \cite{KhovRoz1} but adapted to our choice of potential.  One purpose of the section is to show that the results of \cite{Gornik} hold for more general choices of potential than those considered in that paper.  As a consequence of the results in this section it will follow that for $D$ a diagram of an $l$-component link, $HKh_w$ consists of $n^l$ copies of $\mathbb{C}$.  We give explicit chain representatives of these $n^l$ generators of the homology.  When $l=1$, $HKh_w(D)$ is supported in homological degree $0$.

Although conventional practice dictates that one should deal with invariance under Reidemeister moves before other considerations, issues of logical dependence encourage us to postpone this until Section 4.

Section 3 deals with Morse moves for knot cobordisms as in Figure \ref{morsemoves}.  To link diagrams $D$, $D'$ differing by a Morse move we associate an isomorphism $HKh_w(D) \rightarrow HKh_w(D')$, induced by a chain map $CKh_w(D) \rightarrow CKh_w(D')$, that is filtered of quantum degree $1-n$ in the case of the $0$-handle and $2$-handle Morse moves and of degree $n-1$ in the case of the $1$-handle Morse move.  We see that the maps we define have good properties in terms of preserving the generators defined in Section 4.

Section 4 is also  computational.  The purpose of this section is to show that for link diagrams $D$ and $D'$ differing by a Reidemeister move I or II there is a map, filtered of quantum-degree $0$, $HKh(D) \rightarrow HKh(D')$, induced from a chain map $CKh_w(D) \rightarrow CKh_w(D)$. If we could show this also for Reidemeister move III, a consequence would be that the graded groups associated to the quantum filtration $HKh_w^{i,j}(D) = \mathcal{F}^jHKh^i_w / \mathcal{F}^{j-1}HKh^i_w$ would be invariants of the link represented by the diagram $D$.  Although we do not do this for Reidemeister III, a topological trick detailed in a later section means that we still have enough to prove our main result on the slice genus.

Sections 3 and 4 taken together tell us that if we have a surface $\Sigma$ with boundary

\[ \Sigma \hookrightarrow [0,1]\times\mathbb{R}^3 \]

\[ \partial \Sigma = (K \hookrightarrow \{ 0 \} \times \mathbb{R}^3) \coprod (U \hookrightarrow \{ 1 \} \times \mathbb{R}^3 ) \]

\noindent where $K$ is a $1$-component knot and $U$ is the unknot \emph{and} we have a representation of $\Sigma$ as a product of elementary cobordisms (section 1.1), starting with a diagram $D$ for $K$ and ending with the $0$-crossing diagram of the unknot, without using Reidemeister move III, then

\[ (n-1)(2\rm{genus}(\Sigma) + 3) > | j |  \]

\noindent where $HKh^{0,j}_w(D) \not= 0$.

Section 5 will detail a topological trick which shall mean that we can avoid doing for Reidemeister III a similar computation to those in Section 2 and yet still obtain our main theorem \ref{maintheorem}, the proof of which is completed in this section.

%% Section 6 is intended to put our results in the context of the recent progress made in understanding the relationship between the homology theories $HKh_n$ for varying $n$.

%There is inevitably an appendix into which we have shoehorned some necessary results which could not be made to fit into the beautiful inexorable flow and progression of idea upon idea that is this thesis.

\section{Homology of a trivalent graph}

% Still need to define \chi_0 and \chi_1 and check that their
% composition is what we want it to be.

We have required that $\partial_x w$ is a product of distinct linear
factors.  This condition is the main ingredient in many of results in \cite{Gornik}, the
proofs of which carry across very easily to our potential.  We shall
content ourselves, most often, with stating these results as they apply in our
case without proof, but in a few cases we shall give a proof where things are more difficult for our potential than for that of Gornik's or where we feel it is important for the rest of the paper to understand the proof.

Assume now that $\Gamma$ is a closed graph so that $C(\Gamma)$ (the matrix factorization associated to $\Gamma$)
is actually a 2-periodic complex (see Equation \ref{inandoutburger}), call the homology of this $2$-complex $H(\Gamma)$.  Since $C(\Gamma)$ is filtered $H(\Gamma)$ is also filtered $\ldots \subseteq \mathcal{F}^i H(\Gamma) \subseteq \mathcal{F}^{i+1} H(\Gamma) \subseteq \ldots$.  Let $R(\Gamma)$ be the polynomial ring over $\mathbb{C}$ generated by variables $X_e$ as $e$ runs over the thin edges of $\Gamma$ (note: not just the thin edges with boundary as in the previous section).  Certainly the homology $H(\Gamma)$ is a module over $R(\Gamma)$, but in fact it is also a
module over the ring $\bar{R}(\Gamma)$ obtained from $R(\Gamma)$ by quotienting
out any polynomial in the $X_e$'s appearing in the definition of the fundamental factorizations making up $C(\Gamma)$, (since these polynomials define
$0$-homotopic endomorphisms of $C(\Gamma)$).

\begin{proposition} (See \cite{Gornik} 2.4) The algebra $\bar{R}(\Gamma)$ is generated by
  elements $X_e$ (we are using $X_e$ here to mean the image of $X_e \in R(\Gamma)$ in $\bar{R}(\Gamma)$) where $e$ runs over the thin edges of $\Gamma$, and
  each $X_e$ satisfies

\[\partial_x w(X_e) = 0 \]

\end{proposition}

Writing $\Sigma_n$ for the roots of $\partial_x w$, we make an assignment

\[ \phi : e(\Gamma) \rightarrow \Sigma_n \]

\noindent of elements of $\Sigma_n$ to the thin edges $e(\Gamma)$ of $\Gamma$
is called a \emph{state} of $\Gamma$.

\begin{proposition}  (See \cite{Gornik} 2.5) Given a state $\phi$, define $Q_\phi \in
  \bar{R}(\Gamma)$ by

\[ Q_\phi = \prod_{e \in e(\Gamma)} \left( \prod_{\alpha \in \Sigma_n
  \setminus \phi(e)} \frac{1}{\phi(e) - \alpha} \right) \frac{1}{n+1} \frac{\partial_x w(X_e)}{X_e -
  \phi(e)} \]

\noindent where by $\partial_xw(X_e)/(X_e - \phi(e))$ we mean the result of substituting $X_e$ for $x$ in the polynomial $\partial_xw(x)/(x - \phi(e))$.

%% This /does/ actually agree precisely with Gornik's definition.
%  This makes the Q's idempotents.  Multiplying Q by X is the same as
%  multiplying Q by \alpha.  The fact that the Q's sum to 1 is a bit
%  more tricky - just check how they act (by multiplication) on a
%  C-basis of the space.

Then we have

\[ Q_{\phi_1}Q_{\phi_2} = \left\{ \begin{array}{cc}
Q_{\phi_1}, & \phi_1
=\phi_2 \\ 0, & \phi_1 \not= \phi_2
\end{array} \right. \]

\[ \sum_\phi Q_\phi = 1 \]

\end{proposition}

\begin{proof}

It suffices to check that the relations hold in the ring

\[\bigotimes_{e \in e(\Gamma)} \mathbb{C}[X_e]/\partial_x w (X_e) \]

\noindent since $\bar{R}(\Gamma)$ is a quotient of this ring.

\end{proof}

\begin{figure}
\centerline{
{
\psfrag{lam1}{$\lambda_1$}
\psfrag{lam2}{$\lambda_2$}
\psfrag{thing1}{$\lambda_1,\lambda_2 \in \Sigma_n$}
\psfrag{thing2}{$\lambda_1 \not= \lambda_2$}
\psfrag{type1}{$\rm{Type} \; \rm{I}$}
\psfrag{type2}{$\rm{Type} \; \rm{II}$}
\psfrag{6}{$6$}
\psfrag{id}{$id$}
\psfrag{chi0}{$\chi_0$}
\psfrag{chi1}{$\chi_1$}
\psfrag{phi}{$\phi$}
\includegraphics[height=2in,width=3in]{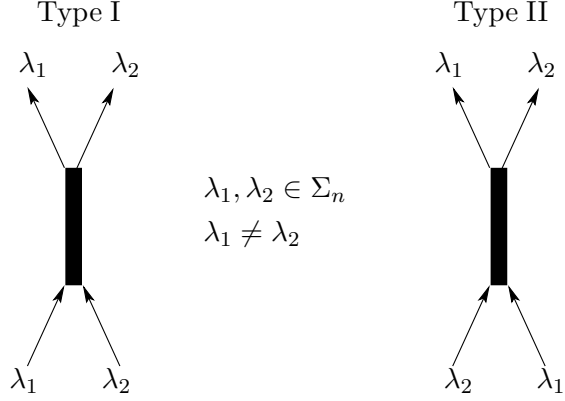}
}}
\caption{Type I and II admissable states}
\label{admissable}
\end{figure}

Given an embedding of $\Gamma$ into $\mathbb{R}^2$ we can define the notion of \emph{admissability} of states.  In the construction of the Khovanov homology of a knot, we start with a diagram of the knot and all $\Gamma$'s appearing in the construction will come with an embedding into $\mathbb{R}^2$ induced by the knot diagram.
A state $\phi$ is called \emph{admissable} if the assignment $\phi :
\Sigma_n \rightarrow e(\Gamma)$ looks like either of the possibilities
in Figure \ref{admissable} in a neighbourhood of each thick edge.  We call the set
of all states $S'(\Gamma)$ and the set of admissable states $S(\Gamma)$.

\begin{theorem} (See \cite{Gornik} 3) For non-admissable states $\phi$ we have

\[ Q_\phi = 0 \rm{.}\]

\noindent For admissable states $\phi$ we have

\[0 \not= \mathbb{C}Q_\phi = \bar{R}(\Gamma)Q_\phi \, \rm{and} \, \rm{so} \,
dim_\mathbb{C} \bar{R}(\Gamma)Q_\phi = 1 \rm{.}\]

\noindent We have a decomposition as a $\mathbb{C}$-algebra

\[ \bar{R}(\Gamma) = \oplus_{\phi \in S(\Gamma)} \mathbb{C}Q_\phi \]

\noindent and

\[ dim_\mathbb{C}\bar{R}(\Gamma) = P_n (\Gamma) |_{q=1} \rm{.} \]

\end{theorem}

%%% Need to say what boundary maps are...

\begin{proposition} Let $k \in \mathbb{Z}$.  Let $Kh^k (\Gamma)$ be
  the quantum degree $k$ piece of the classical $sl(n)$ Khovanov-Rozansky homology of
  the graph $\Gamma$.  Then there is an isomorphism of vector spaces

\[ \Phi : Kh^k (\Gamma) \rightarrow F^k H(\Gamma)/F^{k-1}H(\Gamma) \rm{.} \]

\end{proposition}

\begin{corollary}

Since the graded dimension of $Kh(\Gamma)$ is just $P_n(\Gamma)$ (the
$sl(n)$ polynomial of $\Gamma$), the filtered dimension of $H(\Gamma)$
is $P_n(\Gamma)$.  The number of admissable states is $P_n(\Gamma)(1)$
so that the dimension of $H(\Gamma)$ as a complex vector space agrees
with the number of admissable states.

\end{corollary}

Here since the translation of Gornik's proof is not completely
straightforward, we give some details.

\begin{proof} (of Proposition)  Khovanov and Rozansky have shown that
  the classical 2-periodic complex
  $CKh(\Gamma)$ has cohomology only in one of the two homological
  gradings.  Suppose without loss of generality that it lies in grading
  $1$.

$CKh(\Gamma)$ looks like

\[CKh^0 (\Gamma) \stackrel{d_{0}}{\longrightarrow} CKh^1 (\Gamma) \stackrel{d_{1}}{\longrightarrow}
CKh^0 (\Gamma) \]

\noindent where $d_0$ and $d_1$ are of graded degree $n+1$.  $C(\Gamma)$ looks like

\[CKh^0 (\Gamma) \stackrel{d'_{0}}{\longrightarrow} CKh^1 (\Gamma) \stackrel{d'_{1}}{\longrightarrow}
CKh^0(\Gamma) \]

\noindent where $d'_0$ and $d'_1$ are filtered of degree $n+1$.  In fact we can
decompose $d'_i$ as

\[d'_i = d_i^0 + d_i^1 + d_i^2 + \ldots \]

\noindent where $d_i^j$ is graded of degree $n + 1 - 2j$ and $d_i^0 =
d_i$ for $i=1,2$.
To define the $\Phi$ mentioned in the Proposition we give first a
lift of $\Phi$

\[\phi : (\ker d_1)^k \rightarrow F^k H(\Gamma) / F^{k-1} H(\Gamma) \]

\[\phi : \alpha \mapsto \alpha + \alpha^1 + \alpha^2 + \ldots \]

\noindent where $\alpha^i$ has degree $k-2i$.  We define each $\alpha^i$
inductively.  We
wish to define each $\alpha_i$ so that

\begin{equation}
\label{joycer}
\sum_{l=0}^k d_1^i (\alpha^{j-i}) = 0 \; \forall j
\end{equation}

Writing $\alpha = \alpha^0$ gives us the root case.  Suppose that we
know (\ref{joycer}) for $k \leq K$; we wish to define $\alpha^{K+1}$ so
that we have (\ref{joycer}) for $k \leq K+1$.

If we can show that the following holds:

\begin{equation}
\label{dewar}
d_0^0 \left( \sum_{i=0}^K d_1^{i+1} (\alpha^{K-i}) \right)= 0
\end{equation}

\noindent then we will be done since we know that $\ker d_0^0 = \mathrm{im} d_1^0 $ so we can find an
$\alpha^{K+1}$ satisfying

\[ d_1^0 ( - \alpha^{K+1} ) = \sum_{i=0}^K d_1^{i+1} (\alpha^{K-i}) \rm{.} \]

Now,

\begin{eqnarray*}
\lefteqn{d_0^0 \left( \sum_{i=0}^{K} d_1^{i+1}(\alpha^{K-i})  \right)  } \\
&& = \sum_{i=0}^{K} d_0^0 d_1^{i+1}(\alpha^{K-i}) = -\sum_i^K \sum_{j=0}^i d_0^{j+1} d_1^{i-j}(\alpha^{K-i}) \\
&& = -\sum_{j=0}^K \sum_{i=j}^K d_0^{j+1} d_1^{i-j}(\alpha^{K-i}) = -\sum_{j=0}^K d_0^{j+1} \sum_{i=0}^{K-j} d_1^i (\alpha^{(K-j)-i}) = 0
\end{eqnarray*}

It is easy to check that the map $\phi$ so defined gives a
well-defined isomorphism $\Phi$.

%%Let me just check that...  Yes, it does.

\end{proof}

\begin{proposition}
\label{HGammastructure}

\[ H(\Gamma) = \bigoplus_{\phi \in S(\Gamma)} Q_\phi H(\Gamma) \]

\[ \rm{dim}_{\mathbb{C}} Q_\phi H(\Gamma) = 1 \]

\end{proposition}

\begin{proof}

For dimensional reasons, it is enough to show that for any $\phi \in
S(\Gamma)$ we can find a non-zero element of $Q_\phi H(\Gamma)$.
Such a non-zero element is described explicitly below.

\begin{figure}
\centerline{
{
\psfrag{lam1}{$\lambda_1$}
\psfrag{lam2}{$\lambda_2$}
\psfrag{thing1}{$\lambda_1,\lambda_2 \in \Sigma_n$}
\psfrag{thing2}{$\lambda_1 \not= \lambda_2$}
\psfrag{type1}{$\rm{Type} \; \rm{I}$}
\psfrag{type2}{$\rm{Type} \; \rm{II}$}
\psfrag{6}{$6$}
\psfrag{id}{$id$}
\psfrag{chi0}{$\chi_0$}
\psfrag{chi1}{$\chi_1$}
\psfrag{phi}{$\phi$}
\includegraphics[height=1.5in,width=2in]{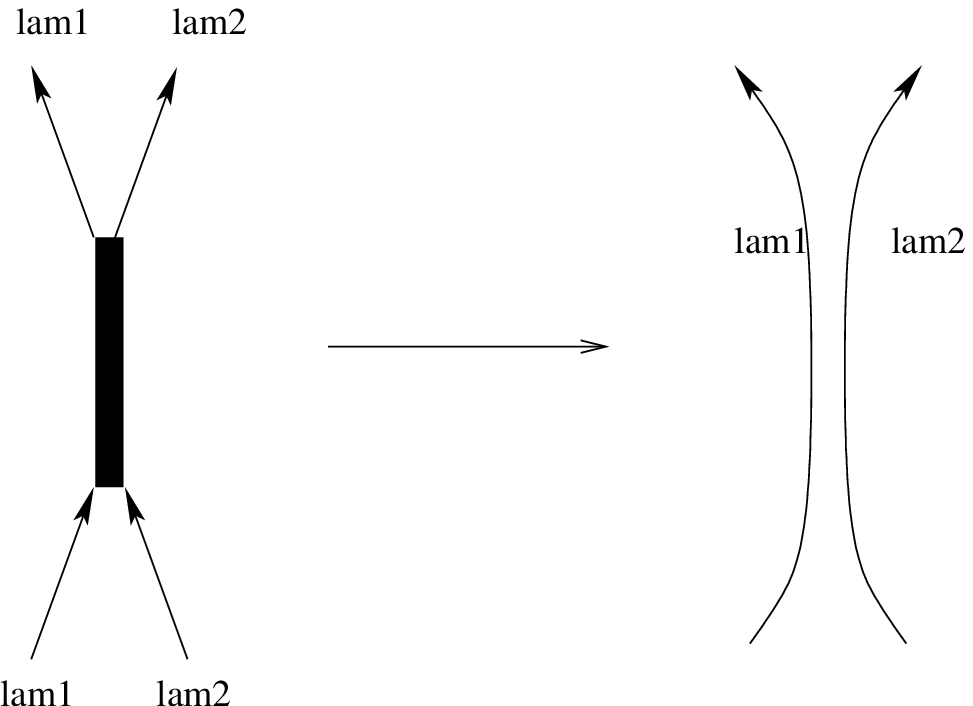}
}}
\caption{Type I}
\label{typeone}
\end{figure}

\begin{figure}
\centerline{
{
\psfrag{lam1}{$\lambda_1$}
\psfrag{lam2}{$\lambda_2$}
\psfrag{thing1}{$\lambda_1,\lambda_2 \in \Sigma_n$}
\psfrag{thing2}{$\lambda_1 \not= \lambda_2$}
\psfrag{type1}{$\rm{Type} \; \rm{I}$}
\psfrag{type2}{$\rm{Type} \; \rm{II}$}
\psfrag{6}{$6$}
\psfrag{id}{$id$}
\psfrag{chi0}{$\chi_0$}
\psfrag{chi1}{$\chi_1$}
\psfrag{phi}{$\phi$}
\includegraphics[height=1.5in,width=2in]{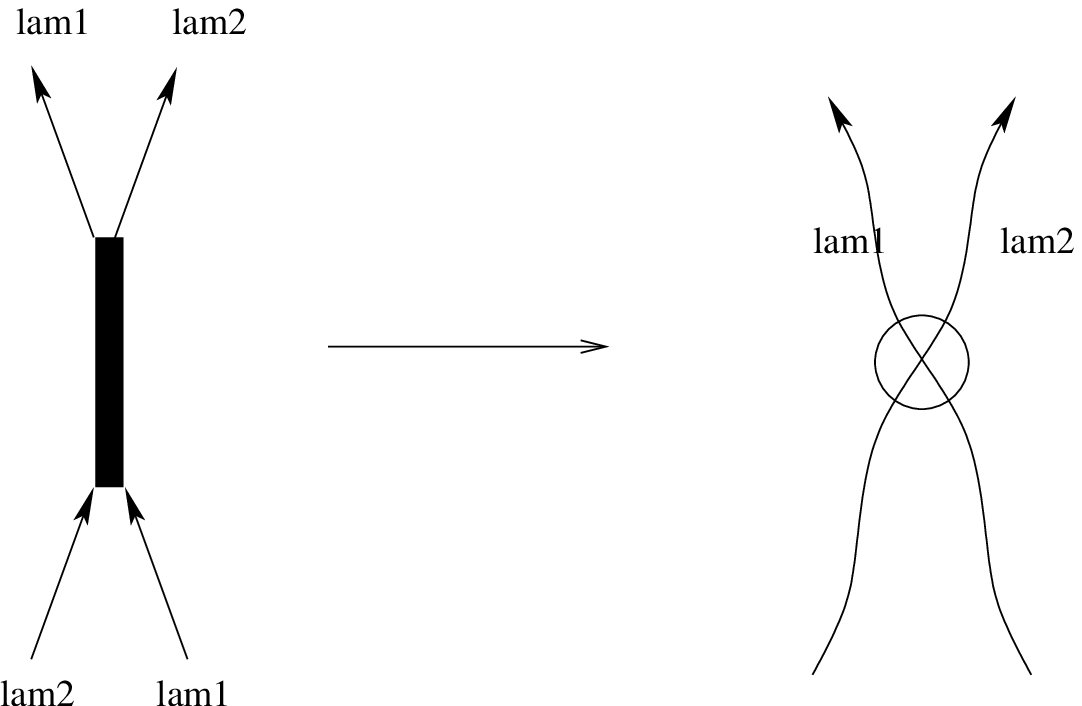}
}}
\caption{Type II}
\label{typetwo}
\end{figure}

% Insert diagram here.

Given $\phi \in S(\Gamma)$ form a new graph $\Gamma'$, with each component decorated with an element $\lambda \in \Sigma$ by following the local rule in Figures \ref{typeone} and \ref{typetwo} at all thick edges.  We remind the reader that a trivalent graph describing a matrix need not come with an embedding into $\mathbb{R}^2$.  The circled crossings in $\Gamma'$ should not be thought of as undercrossings or overcrossings but just as the two strands not intersecting.

% Insert another fucking diagram.

The state $\phi$ gives a corresponding state $\phi'$ of
$\Gamma'$.  Now since $\Gamma'$ is a union of circles, $H(\Gamma') =
R(\Gamma')$ naturally so we can pick unambiguously a $Q_{\phi'} \in H(\Gamma)$.

We apply $\chi_0$ and $\eta_0$ repeatedly to $H(\Gamma')$ until we arrive in
$H(\Gamma)$.  ($\chi_0$ is defined in the remainder of this section and $\eta_0$ is defined in Appendix \ref{the eta map}).  We define $Q_\phi$ to be the image of $Q_{\phi'}$ under
these map.  Gornik shows that this element is non-zero.

\end{proof}

To define $HKh_w(D)$ for a link diagram $D$ we have to
give two maps of matrix factorizations

\[ \chi_0 : C(\Gamma_0) \rightarrow C(\Gamma_1) \]

\[ \chi_1 : C(\Gamma_1) \rightarrow C(\Gamma_0) \]

\begin{figure}
\centerline{
{
\psfrag{x1}{$x_1$}
\psfrag{x2}{$x_2$}
\psfrag{x3}{$x_3$}
\psfrag{x4}{$x_4$}
\psfrag{type1}{$\rm{Type} \; \rm{I}$}
\psfrag{type2}{$\rm{Type} \; \rm{II}$}
\psfrag{6}{$6$}
\psfrag{id}{$id$}
\psfrag{chi0}{$\chi_0$}
\psfrag{chi1}{$\chi_1$}
\psfrag{phi}{$\phi$}
\includegraphics[height=1.5in,width=0.5in]{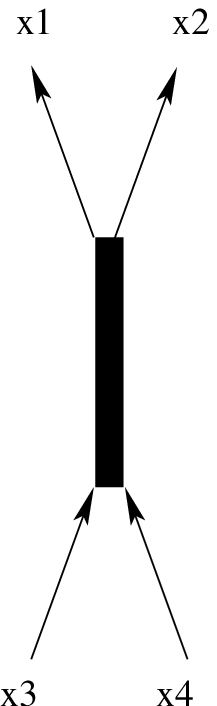}
}}
\caption{$\Gamma_1$}
\label{Gamma1}
\end{figure}

\begin{figure}
\centerline{
{
\psfrag{x1}{$x_1$}
\psfrag{x2}{$x_2$}
\psfrag{x3}{$x_3$}
\psfrag{x4}{$x_4$}
\psfrag{id}{$id$}
\psfrag{chi0}{$\chi_0$}
\psfrag{chi1}{$\chi_1$}
\psfrag{phi}{$\phi$}
\includegraphics[height=1.5in,width=0.5in]{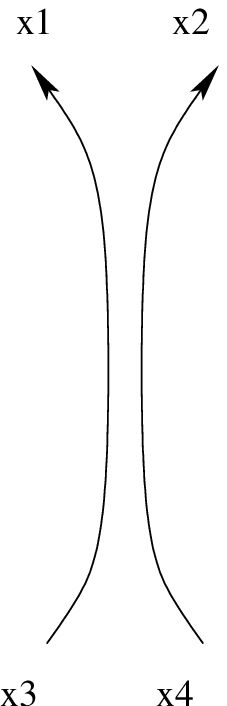}
}}
\caption{$\Gamma_0$}
\label{Gamma0}
\end{figure}

% Have to give pictures of $\Gamma_0$ and $\Gamma_1$ here - or
% whatever it is that Khovanov calls them.

From these maps we build up the differentials in the chain complex of
matrix factorization $CKh_w (D)$.  Given a knot diagram $D$, the Khovanov cube with closed trivalent graphs (giving matrix factorizations of $0$ i.e. $2$-periodic complexes) at the vertices and maps of matrix factorizations on each edge is given by taking the tensor product of

\[ 0 \rightarrow C(\Gamma_0)\{1-n\} \stackrel{\chi_0}{\rightarrow} C(\Gamma_1)\{-n\} \rightarrow 0 \]

\noindent for every positive crossing (the left of Figure \ref{writhe}) and

\[ 0 \rightarrow C(\Gamma_1) \{ n \} \stackrel{\chi_1}{\rightarrow} C(\Gamma_0)\{ n-1 \} \rightarrow 0 \]

\noindent for every negative crossing (the right of Figure \ref{writhe}).  This definition clearly agrees with that of Subsection 1.5.

We give an explicit description ($R$ is the ring
$\mathbb{C}[x_1,x_2,x_3,x_4]$). The factorization $C(\Gamma_0)$ is

\[ \left( \begin{array}{c} R \\ R\{2-2n\} \end{array} \right)
\stackrel{P_0}{\longrightarrow} \left( \begin{array}{c} R\{1-n\} \\ R\{1-n\}
\end{array} \right) \stackrel{P_1}{\longrightarrow} \left( \begin{array}{c} R \\
  R\{2-2n\} \end{array} \right) \]

\[ P_0 = \left( \begin{array}{cc} \pi_{13} & x_2 - x_4 \\ \pi_{24} & x_3
  - x_1 \end{array} \right) \]

\[ P_1 = \left( \begin{array}{cc} x_1 - x_3 & x_2 - x_4 \\ \pi_{24} & -\pi_{13} \end{array} \right) \]

\noindent and $C(\Gamma_1)$ is the factorization

\[ \left( \begin{array}{c} R\{-1\} \\ R\{3-2n\} \end{array} \right)
\underrightarrow{Q_0} \left( \begin{array}{c} R\{-n\} \\ R\{2-n\}
\end{array} \right) \underrightarrow{Q_1} \left( \begin{array}{c} R\{-1\} \\
  R\{3-2n\} \end{array} \right) \]

\[ Q_0 = \left( \begin{array}{cc} u_1 & x_1x_2 - x_3x_4 \\ u_2 & x_4 +
  x_3 - x-1 - x_2 \end{array} \right) \]

\[ Q_1 = \left( \begin{array}{cc} x_1 + x_2 - x_3 - x_4 & x_1x_2 - x_3x_4 \\ u_2 & -u_1 \end{array} \right) \]

% Insert a couple of factorizations here.

% \mu = 0, \lambda = 1

The map $\chi_0 : C(\Gamma_0) \rightarrow C(\Gamma_1)$ is defined by the pair of maps:

\[ U_0 = \left( \begin{array}{cc} x_3 - x_2 & 0 \\ a & 1 \end{array} \right) \]

\[ U_1 = \left( \begin{array}{cc} x_3 & -x_2 \\ -1 & 1 \end{array} \right) \]

\noindent Here

\[ a = \frac{u_1 + x_3u_2 - \pi_{24}}{x_1 - x_3} \rm{.}\]

The map $\chi_1 : C(\Gamma_1) \rightarrow C(\Gamma_0)$ is defined by the pair of maps:

\[ V_0 = \left( \begin{array}{cc} 1 & 0 \\ -a & x_3 - x_2 \end{array} \right) \]

\[ V_1 = \left( \begin{array}{cc} 1 & x_2 \\ 1 & x_3 \end{array} \right) \rm{.}\]

Note that $\chi_0 \chi_1$ and $\chi_1 \chi_0$ are homotopic to multiplication by $x_2 - x_3$ (we can just take the zero homotopies in both cases).

\begin{proposition}

The Khovanov homology groups of an $l$-component link diagram $D$ for
our potential $w$, $Kh_w(D)$, have dimension

\[ \dim_\mathbb{C} = n^l \]

There is a canonical basis of generators, one generator for each
assignment

\[ \psi : components(D) \rightarrow \Sigma_n \].

Each such $\psi$ defines in an obvious way a \emph{type II} admissable
state of a resolution $\Gamma_\psi$.  The generator corresponding to
$\psi$ lives in the chain group summand $H(\Gamma_\psi)$, and comes
from the state $\psi$ as described in (\ref{HGammastructure}).

\end{proposition}

It is these generators that will play the same roles in this paper as Lee's
generators \cite{Lee} played in Rasmussen's slice genus result \cite{Ras}

\section{Morse moves}

A link cobordism can be written as a finite sequence of link diagrams,
where successive diagrams differ either by a Reidemeister
move or a Morse move.  The Morse moves correspond to adding $0$-,
$1$-, and $2$-handles.

In this section we will assign to each Morse move a filtered chain map
between the complexes of the diagrams that it connects.  Then we shall
compute what each Morse moves does to the canonical generators of the
homology of a link diagram.

\subsection{$0$-handle move}

\begin{figure}
\centerline{
{
\psfrag{0handle}{$0$-handle}
\includegraphics[height=1.4in,width=3in]{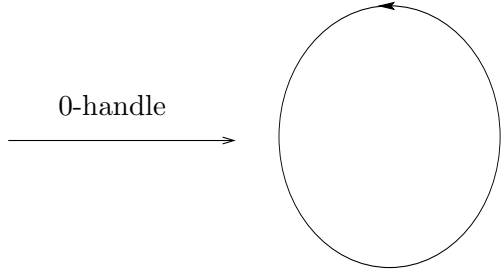}
}}
\caption{$0$-handle addition}
\label{0-handle}
\end{figure}

A $0$-handle move is the creation of a simple loop as in Figure \ref{0-handle}.

%% Insert diagram.

The $2$-periodic complex associated to a circle is

\[ A<1> = (0 \rightarrow A \rightarrow 0) \]

\noindent where $A = (\mathbb{C}[x]/\partial_x w(x)) \{ 1-n \}$ (the curly brackets indicate a shift in the filtration).

The unit map $i : \mathbb{C} \rightarrow A$ has filtered degree
$1-n$.  To the $0$-handle move above we associate the map of complexes

\begin{eqnarray*}
1 \otimes i : CKh_w (D) \rightarrow CKh_w(D \sqcup S^1)
\end{eqnarray*}

The canonical generators of the cohomology of $A<1>$ are the chain
elements

\[ q_\beta =  \frac{1}{n+1} \left( \prod_{\alpha \in \Sigma_n
  \setminus \beta} \frac{1}{\beta - \alpha} \right) \frac{\partial_x w(x)}{x -
  \beta} \]

\noindent and they satisfy

\[ 1 = \sum_{\beta \in \Sigma_n} q_\beta \]

So our map $1 \otimes i$ takes a canonical generator $g \in CKh_w(D)$ to
the sum of canonical generators $\sum_{\beta \in \Sigma_n} g \otimes q_\beta$.

\subsection{$2$-handle move}

\begin{figure}
\centerline{
{
\psfrag{2handle}{$2$-handle}
\includegraphics[height=1.4in,width=3in]{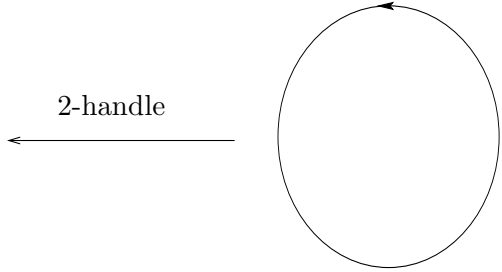}
}}
\caption{$2$-handle addition}
\label{2-handle}
\end{figure}

The $2$-handle move is the removal of a simple closed loop as in Figure \ref{2-handle}.

%Insert diagram.

The trace map $ \epsilon : A \rightarrow \mathbb{C} $ is defined by

\[ \epsilon (x^i) = \left\{ \begin{array}{cc}
1, & i = n-1 \\ 0, & i < n-1
\end{array} \right. \]

To the $2$-handle move above we associate the map of complexes

\[ 1 \otimes \epsilon : CKh_w (D \sqcup S^1) \rightarrow CKh_w(D) \]

\noindent and this has filtered degree $1-n$. Since the coefficient of $x^{n-1}$ in $q_\beta$ is just $\prod_{\alpha \in \Sigma_n
  \setminus \beta} \frac{1}{\beta - \alpha}$ our map takes a canonical
  generator $g \otimes q_\beta \in CKh_w(D \sqcup S^1)$ to

\[ \left( \prod_{\alpha \in \Sigma_n
  \setminus \beta} \frac{1}{\beta - \alpha} \right) g \]

\subsection{$1$-handle move}

\begin{figure}
\centerline{
{
\psfrag{1handle}{$1$-handle}
\psfrag{x1}{$x_1$} 
\psfrag{x2}{$x_2$} 
\psfrag{x3}{$x_3$} 
\psfrag{x4}{$x_4$} 
\psfrag{Gamma0}{$\Gamma_0$} 
\psfrag{Gamma1}{$\Gamma_1$} 
\includegraphics[height=1.2in,width=3in]{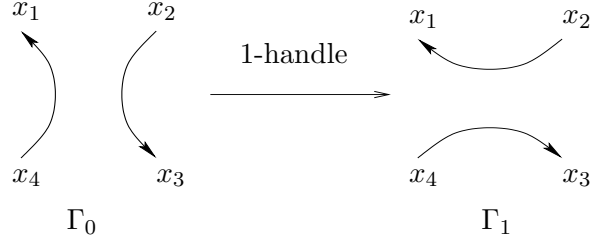}
}}
\caption{$1$-handle addition}
\label{1-handle}
\end{figure}

The $1$-handle move is the addition of a saddle as in Figure \ref{1-handle}.

% Draw picture.

Here $C(\Gamma_0)$ is the factorization

\[ \left( \begin{array}{c} R \\ R\{2-2n\} \end{array} \right)
\underrightarrow{P_0} \left( \begin{array}{c} R\{1-n\} \\ R\{1-n\}
\end{array} \right) \underrightarrow{P_1} \left( \begin{array}{c} R \\
  R\{2-2n\} \end{array} \right) \]

\[ P_0 = \left( \begin{array}{cc} \pi_{14} & -(x_3 - x_2) \\ \pi_{23} & x_1
  - x_4 \end{array} \right) \]

\[ P_1 = \left( \begin{array}{cc} x_1 - x_4 & x_3 - x_2 \\ -\pi_{23} & \pi_{14} \end{array} \right) \]

And $C(\Gamma_1)<1>$ is the factorization

% Should just relabel 2 <---> 4 from above factorization.  Let's do
% that...

\[ \left( \begin{array}{c} R\{1-n\} \\ R\{1-n\} \end{array} \right)
\underrightarrow{Q_0} \left( \begin{array}{c} R \\ R\{2-2n\}
\end{array} \right) \underrightarrow{Q_1} \left( \begin{array}{c} R\{1-n\} \\
  R\{1-n\} \end{array} \right) \]

\[ Q_0 = \left( \begin{array}{cc} x_1 - x_2 & x_3 - x_4 \\ -\pi_{34} & \pi_{12} \end{array} \right) \]

\[ Q_1 = \left( \begin{array}{cc} \pi_{12} & -(x_3 - x_4) \\ \pi_{34} & x_1
  - x_2 \end{array} \right) \]

%Does this work with our choice below of F?  Let me check...
%We need either F_0 ---> -F_0 or F_1 ---> F_1, otherwise
%we do indeed have a matrix factorization morphism.

Ignoring the grading shift, the filtered $\mathbb{C}[x_1, x_2, x_3, x_4]$-module of filtered
  matrix factorization maps $C(\Gamma_0) \rightarrow C(\Gamma_1)<1>$ is isomorphic to
  the quotient  $\mathbb{C}[x_1, x_2, x_3, x_4]/(x_1 = x_2 = x_3 =
  x_4, \partial_x w(x_1)=0)$.  A generator is given by the pair of
  matrices:

%%\footnote{Our definition of $F$ corrects an omission of a
%%  factor of $1/2$ missing from the '$e$' entries of the [..] $F$.}

\[ F_0 = \left( \begin{array}{cc} e_{123}
  & -1 \\ -e_{134} & -1 \end{array}
  \right) \]

\[ F_1 = \left( \begin{array}{cc} -1 & 1 \\ e_{123} & e_{134}
  \end{array} \right) \]

\noindent where

\[ e_{ijk} = \frac{\pi_{ik} - \pi_{jk}}{x_i - x_j} \]

The map $F$ is filtered of degree $n-1$.

Suppose we have two knot diagrams $D_1$ and $D_2$ (of links $L_1$
$L_2$) such that $D_2$ is
got from $D_1$ by the addition of a $1$-handle.  We will compute the
map induced by $F$

\[ HKh_w (D_1) \rightarrow HKh_w (D_2) \]

\noindent in terms of the generators of $HKh_w (D_1)$ and $HKh_w (D_2)$ that were constructed in Section $1$.

Recall that for each choice of assignment $\phi : \{components \: of
  \: L_1\} \rightarrow \Sigma_n$ we got a corresponding resolution
  $\Gamma_\phi$ of $D_1$ and a state $\phi : e(\Gamma_\phi)
  \rightarrow \Sigma_n$ giving us an element $Q_\phi \in
  H(\Gamma_\phi)$ which represents one of our basis elements of $CKh_w
  (D_1)$.

To construct $Q_\phi$ we followed the recipe which converted the
resolution $\Gamma_\phi$ into a disjoint union of circles
$\Gamma'_\phi$ and pushing forward an element (again determined by
$\phi$) of $H(\Gamma'_\phi)$ to $H(\Gamma_\phi)$ via repeated use of
$\eta_0$.  Since $F$ commutes with $\eta_0$, it is enough for us to
determine the map induced by $F$ acting on the homology of disjoint circles.

There are thus two possibilities: either the $1$-handle move joins two
components of the link as in Figure \ref{1-handleI} or it splits a single component into two as in Figure \ref{1-handleII}.

\subsubsection{First case}

% Insert diagram of first case.  \Gamma'_0 etc.

\begin{figure}
\centerline{
{
\psfrag{1handle}{$1$-handle}
\psfrag{1234}{$1$-handle}
\psfrag{x1}{$x_1$} 
\psfrag{x2}{$x_2$} 
\psfrag{x3}{$x_3$} 
\psfrag{x4}{$x_4$} 
\psfrag{Gamma0}{$\Gamma'_0$} 
\psfrag{Gamma1}{$\Gamma'_1$} 
\includegraphics[height=1.4in,width=5in]{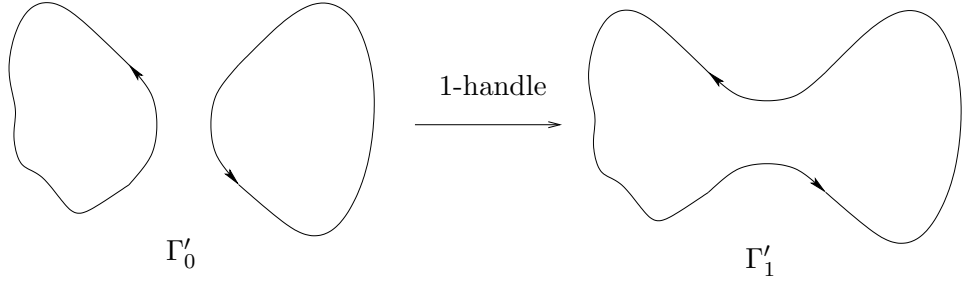}
}}
\caption{$1$-handle addition, first case}
\label{1-handleI}
\end{figure}

% 1 = 4 , 2 = 3.

Performing the endpoint indentifications in Figure \ref{1-handleI}, we have

\[ C(\Gamma'_0) = \left( \begin{array}{c} R \\ R\{2-2n\} \end{array} \right)
\underrightarrow{P'_0} \left( \begin{array}{c} R\{1-n\} \\ R\{1-n\}
\end{array} \right) \underrightarrow{P'_1} \left( \begin{array}{c} R \\
  R\{2-2n\} \end{array} \right) \]

\[ P'_0 = \left( \begin{array}{cc} \partial_x w(x_1) & 0 \\ \partial_x
  w(x_2) & 0 \end{array} \right) \]

\[ P'_1 = \left( \begin{array}{cc} 0 & 0 \\ -\partial_x w(x_2) &
  \partial_x w(x_1) \end{array} \right) \]

\noindent where $R = \mathbb{C}[x_1,x_2]$ thought of as an infinite
dimensional $\mathbb{C}$-module.

\[ C(\Gamma'_1)<1> = \left( \begin{array}{c} R\{1-n\} \\ R\{1-n\} \end{array} \right)
\underrightarrow{Q'_0} \left( \begin{array}{c} R \\ R\{2-2n\}
\end{array} \right) \underrightarrow{Q'_1} \left( \begin{array}{c} R\{1-n\} \\
  R\{1-n\} \end{array} \right) \]

\[ Q'_0 = \left( \begin{array}{cc} (x_1 - x_2) & -(x_1 - x_2) \\ -\pi_{12} & \pi_{12} \end{array} \right) \]

\[ Q'_1 = \left( \begin{array}{cc} \pi_{12} & (x_1 - x_2) \\ \pi_{12} & (x_1
  - x_2) \end{array} \right) \]

And the map $F'$ is given by

\[ F'_0 = \left( \begin{array}{cc} e' & -1 \\ -f' & -1 \end{array} \right)
\]

\[ F'_1 = \left( \begin{array}{cc} -1 & 1 \\ e' & f' \end{array}
\right) \]

\noindent where

\[ e' = \frac{\pi_{12} - \partial_x w(x_2)}{x_1 - x_2}, 
\; f' = \frac{\partial_x w (x_1) - \pi_{12}}{x_1 - x_2} \]

%homology lies in the zeroeth bit.

We observe that $H(\Gamma'_0) = ker(P'_0)/im(P'_1)$.  The canonical
generators indexed by pairs of roots $\alpha, \beta \in \Sigma_n$
have representatives in the zeroeth chain group of $C(\Gamma_0)$:

\[ \left( \begin{array}{c} 0 \\ q^1_{\alpha} q^2_{\beta} \end{array}
\right) \in \left( \begin{array}{c} R \\ R\{2-2n\} \end{array} \right)
\]

Also, $H(\Gamma'_1<1>) = ker(Q'_0)/im(Q'_1)$, and the canonical
generators indexed by roots $\alpha \in \Sigma_n$ have representatives
in the $0$th chain group of $C(\Gamma_1)$:

\[ \left( \begin{array}{c} q^1_{\alpha} \\ q^1_{\alpha} \end{array}
\right) \in \left( \begin{array}{c} R\{1-n\} \\ R\{1-n\} \end{array} \right)
\]

\noindent where

\[ q^i_{\alpha} = \frac{1}{n+1} \left( \prod_{\beta \in \Sigma_n
  \setminus \alpha} \frac{1}{\alpha - \beta} \right) \frac{\partial_x w(x_i)}{x_i -
  \alpha} \; \rm{for} \; i = 1,2 \]

Now $H ( \Gamma'_1 < 1 > )$ is isomorphic to $\mathbb{C}[x_1, x_2]/(x_1=x_2, \partial_x w(x_1) = 0)$, and in this module

\[ q^1_\alpha q^2_\beta = \left\{ \begin{array}{cc} 0, & \alpha \not= \beta \\ q^1_\alpha, & \alpha = \beta \end{array} \right. \]

So it is easy to see what $\Psi'$ does to the generators of $H (\Gamma'_0)$:

\begin{eqnarray*}
\left( \begin{array}{cc} e' & -1 \\ -f' & -1 \end{array} \right) \left( \begin{array}{c} 0 \\ q^1_{\alpha} q^2_{\beta} \end{array} \right) &=& \left( \begin{array}{c}  -q^1_{\alpha} q^2_{\beta} \\  -q^1_{\alpha} q^2_{\beta} \end{array}
\right) \\
&=& \left( \begin{array}{c} -q^1_\alpha \\ -q^1_\alpha  \end{array} \right) \, \rm{or} \, \left( \begin{array}{c} 0 \\ 0 \end{array} \right)
\end{eqnarray*}

\noindent depending as $\alpha = \beta$ or $\alpha \not= \beta$.

\subsubsection{Second case}

% 1 <---> 2, 3 <---> 4 $

\begin{figure}
\centerline{
{
\psfrag{1handle}{$1$-handle}
\psfrag{x1}{$x_1$} 
\psfrag{x2}{$x_2$} 
\psfrag{x3}{$x_3$} 
\psfrag{x4}{$x_4$} 
\psfrag{Gamma0}{$\Gamma''_0$} 
\psfrag{Gamma1}{$\Gamma''_1$} 
\psfrag{1handle}{$1$-handle}
\includegraphics[height=2.5in,width=4.2in]{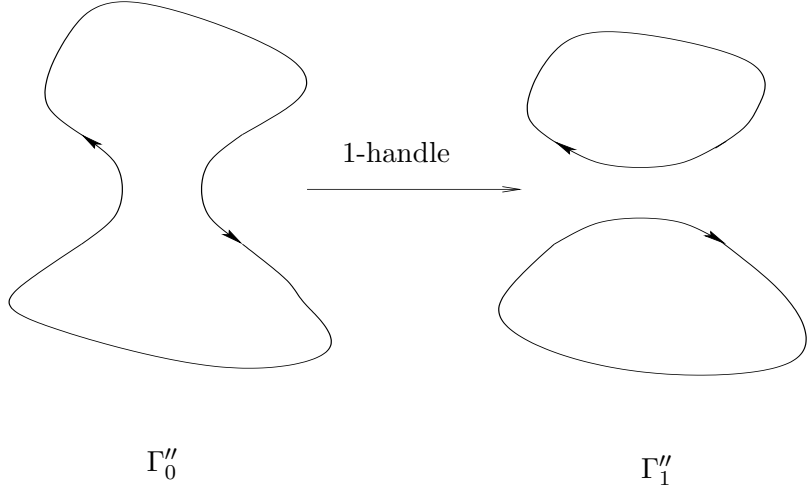}
}}
\caption{$1$-handle addition, second case}
\label{1-handleII}
\end{figure}

Under the identifications in Figure \ref{1-handleII}, we have

\[ C(\Gamma''_0) = \left( \begin{array}{c} R \\ R\{2-2n\} \end{array} \right)
\underrightarrow{P''_0} \left( \begin{array}{c} R\{1-n\} \\ R\{1-n\}
\end{array} \right) \underrightarrow{P''_1} \left( \begin{array}{c} R \\
  R\{2-2n\} \end{array} \right) \]

\[ P''_0 = \left( \begin{array}{cc} \pi_{13} & -(x_3 - x_1) \\ \pi_{13} &
  x_1 - x_3 \end{array} \right) \]

\[ P''_1 = \left( \begin{array}{cc} x_1 - x_3 & x_3 - x_1 \\ -\pi_{13} &
  \pi_{13} \end{array} \right) \]

\noindent where $R = \mathbb{C}[x_1, x_3]$

\[ C(\Gamma''_1)<1> = \left( \begin{array}{c} R\{1-n\} \\ R\{1-n\} \end{array} \right)
\underrightarrow{Q''_0} \left( \begin{array}{c} R \\ R\{2-2n\}
\end{array} \right) \underrightarrow{Q''_1} \left( \begin{array}{c} R\{1-n\} \\
  R\{1-n\} \end{array} \right) \]

\[ Q''_0 = \left( \begin{array}{cc} 0 & 0 \\ -\partial_x w(x_3) &
  \partial_x w(x_1) \end{array} \right) \]

\[ Q''_1 = \left( \begin{array}{cc} \partial_x w(x_1) & 0 \\
  \partial_x w(x_3) & 0 \end{array} \right) \]

The map $F''$ is given by

% eta{123} = \pi.

\[ F''_0 = \left( \begin{array}{cc} e''
  & -1 \\ -f'' & -1 \end{array}
  \right) \]

\[ F''_1 = \left( \begin{array}{cc} -1 & 1 \\ e'' & f''
  \end{array} \right) \]

\noindent where

\[ e'' = \frac{\partial_x w(x_1)  - \pi_{13}}{x_1 - x_3}, 
\; f'' = \frac{\pi_{13} - \partial_x w(x_3)}{x_1 - x_3} \]

We observe that $H(\Gamma''_0) = ker(P''_1)/im(P''_0)$.  The canonical
generators indexed by roots $\alpha \in \Sigma_n$
have representatives in the $0$th chain group of $C(\Gamma''_0)$

\[ \left( \begin{array}{c} q^1_{\alpha} \\ q^1_{\alpha} \end{array}
\right) \in \left( \begin{array}{c} R\{1-n\} \\ R\{1-n\} \end{array} \right)
\]

Also $H(\Gamma''_1<1>) = ker(Q''_1)/im(Q''_0)$, and the canonical
generators indexed by roots $\alpha, \beta \in \Sigma_n$ have representatives
in the $0$th chain group of $C(\Gamma''_1)<1>$

\[ \left( \begin{array}{c} 0 \\ q^1_{\alpha} q^3_{\beta} \end{array}
\right) \in \left( \begin{array}{c} R \\ R\{2-2n\} \end{array} \right)
\]

\noindent where

\[ q^i_{\alpha} = \frac{1}{n+1} \left( \prod_{\beta \in \Sigma_n
  \setminus \alpha} \frac{1}{\alpha - \beta} \right) \frac{\partial_x w(x_i)}{x_i -
  \alpha} \; \rm{for} \; i = 1,3 \]

Now $H(\Gamma''_1 <1>)$ is isomorphic to $\mathbb{C}[x_1,x_3]/(\partial_x w (x_1), \partial_x w (x_3))$ and in this module

\[ \left( \frac{\partial_x w (x_1) - \partial_x w (x_3)}{x_1 - x_3} \right) \left( \frac{ \partial_x w(x_1)}{x_1 - \alpha} \right) = \left( \frac{ \partial_x w (x_1)}{x_1 - \alpha} \right) \left( \frac{\partial_x w (x_3)}{x_3 - \alpha} \right) \]

\noindent since

\[ \frac{\partial_x w(x_1) - \partial_x w (x_3)}{x_1 - x_3} - \frac{\partial_x w(x_3)}{x_3 - \alpha} \]

\noindent has a factor of $(x_1 - \alpha)$.

Now it is easy to see what $F''$ does to the generators:

\begin{eqnarray*}
\left( \begin{array}{cc} -1 & 1 \\
e'' & f'' \end{array} \right) \left( \begin{array}{c} q^1_{\alpha}
\\ q^1_{\alpha} \end{array}
\right) &=& \left( \begin{array}{c} 0 \\ (e'' + f'')q^1_{\alpha} 
\end{array} \right) \\ 
&=& \left( \begin{array}{c} 0 \\  \frac{\partial_x w (x_1) - \partial_x w (x_3)}{x_1 - x_3} q^1_\alpha  \end{array} \right) \\
&=& \left( \prod_{\beta \in \Sigma_n \setminus \alpha} (\alpha - \beta) \right)\left( \begin{array}{c} 0 \\  q^1_\alpha q^3_\alpha   \end{array} \right)
\end{eqnarray*}

%%\subsection{Summary}

\section{Reidemeister moves I and II}

In this section we prove invariance of $HKh_w(D)$ under changing the (closed) link diagram $D$ by an oriented Reidemeister I or II move.

The diagrams on either side of an oriented Reidemeister move each give a chain complex of matrix factorizations.  We wish to define two chain maps (one in each direction) between these two chain complexes.  Each chain map shall be quantum-graded of degree $0$.  These chain maps will then induce degree $0$ chain maps on the chain complexes associated to \emph{closed} link diagrams differing locally by the oriented Reidemester move.  We have a description of chain representatives for the generators of the homology $HKh_w$ of each link diagram and we show that the maps preserve the generators up to multiplication by a non-zero number.  Since the chain maps were of degree $0$ this tells us that the graded vector space associated to the quantum filtration of $HKh_w$ is invariant under changing the input diagram by the Reidemeister move.

Every matrix in this section describes a map between free modules over polynomial rings.  Each of these modules, as per the definition of $CKh_w$, is possibly subject to a shift in quantum filtration.  In what is below, we suppress mention of the shifts in quantum filtration.  For those readers wishing to check the validity of the following computations, this omission should not encumber them with any great difficulties, while explicitly to describe the shifts would thicken somewhat our exposition.  In this section we write $\partial_i w$ for $\partial_x w(x_i)$.

%% In this section we will associate to Reidemeister moves I and II
%%quantum filtered degree $0$ maps
%%between the complexes of matrix factorizations described by the start
%%and end diagrams of the moves.  Although we believe that these maps
%%are chain homotopy equivalences, we will be more interested in studying
%%the maps induced on the homology $HKh_w (D_1) \rightarrow HKh_w (D_2)$
%%where $D_1$ and $D_2$ are (closed) link diagrams differing by a Reidemeister
%%move.  By this approach we shall see that we have at least quasi-isomorphisms.

\subsection{Reidemeister move I}

\subsubsection{Reidemeister move I.1}

\begin{figure}
\centerline{
{
\psfrag{0}{$x_1$}
\psfrag{x2}{$x_2$}
\psfrag{x3}{$x_3$}
\psfrag{x4}{$x_4$}
\psfrag{id}{$id$}
\psfrag{chi0}{$\chi_0$}
\psfrag{chi1}{$\chi_1$}
\psfrag{phi}{$\phi$}
\includegraphics[height=1.5in,width=3in]{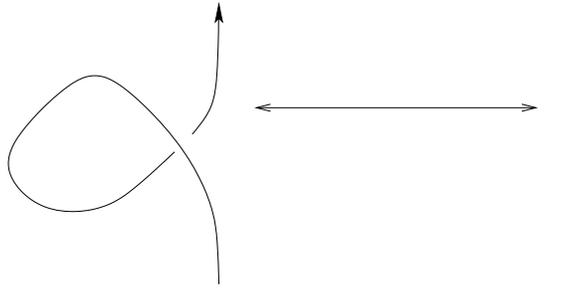}
}}
\caption{Reidemeister move I.1}
\label{R11}
\end{figure}

The first case of the  Reidemeister I move is shown in Figure \ref{R11}.  Either side of the move corresponds to a chain complex of matrix facorizations.  Our first task is to define a chain maps between the left chain complex to and the right complex of degree $0$.  Then we wish to see that these chain maps induce maps on the homologies of two closed link diagrams, differing locally by the move, which preserve the generators of the homology.

\begin{figure}
\centerline{
{
\psfrag{0}{$0$}
\psfrag{F}{$F$}
\psfrag{G}{$G$}
\psfrag{chi0}{$\chi_0$}
\psfrag{M}{$M$}
\psfrag{N}{$N$}
\psfrag{P}{$P$}
\psfrag{phi}{$\phi$}
\psfrag{x1}{$x_1$}
\psfrag{x2}{$x_2$}
\psfrag{x3}{$x_3$}
\includegraphics[height=1.5in,width=3in]{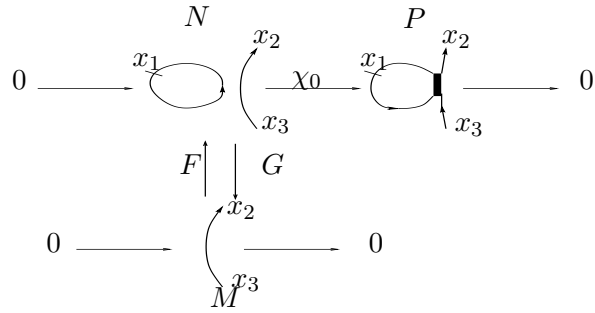}
}}
\caption{Reidemeister move I.1 chain maps}
\label{R11chainmaps}
\end{figure}

The set-up is shown in Figure \ref{R11chainmaps}.  The chain complex from the left of Figure \ref{R11} is above, and the one from the right is below.  To give the two chain maps we need to define $F$ and $G$.  First we shall write down the three factorizations $M$, $N$, and $P$.

The factorizations M, N, and P are

\[ M : R\{ 1-n \} \stackrel{x_2 - x_3}{\rightarrow} R \stackrel{\pi_{23}}{\rightarrow} R \{ 1-n \} \]

\[ N : \left( \begin{array}{c} S \{ 2-2n \} \\ S \{ 2-2n \} \end{array} \right) \underrightarrow{A^0}  \left( \begin{array}{c} S \{ 1-n \} \\ S \{ 3-3n \} \end{array} \right) \underrightarrow{A^1} \left( \begin{array}{c} S \{ 2-2n \} \\ S \{ 2-2n \} \end{array} \right) \]

 \[ A^0 = \left( \begin{array}{cc} 0  &  x_2 - x_3 \\ - \pi_{23} & \partial w_1  \end{array} \right), A^1 = \left( \begin{array}{cc} \partial w_1  &  -(x_2 - x_3) \\ \pi_{23} & 0  \end{array} \right) \]

\[ P : \left( \begin{array}{c} S \{ -2n \} \\ S \{ 2-2n  \} \end{array} \right) \underrightarrow{B^0} \left( \begin{array}{c} S \{ -1-n \} \\ S \{ 3-3n \} \end{array} \right) \underrightarrow{B^1}  \left( \begin{array}{c} S \{ -2n \} \\ S \{ 2-2n \} \end{array} \right) \]

\[ B^0 = \left( \begin{array}{cc} x_2 - x_3 & x_1(x_2 - x_3) \\  -u_2 & u_1  \end{array} \right) , B^1 = \left( \begin{array}{cc} u_1 & -x_1(x_2 - x_3) \\ u_2 & x_2 - x_3  \end{array} \right) \]

\noindent where $R = \mathbb{C}[x_2, x_3]$ and $S = \mathbb{C}[x_1, x_2, x_3]$ considered as a module over $R$.  Also

\[ \pi_{23} = \frac{w(x_2) - w(x_3)}{x_2 - x_3}, \partial w_1 = w'(x_1) \]

\[ u_1 = \frac{p(x_1 + x_2, x_1 x_2) - p(x_1 + x_3, x_1 x_2)}{x_2 - x_3} \]

\[ u_2 = \frac{p(x_1 + x_3, x_1 x_2) - p(x_1 + x_3, x_1 x_3)}{x_1 (x_2 - x_3)} \]

\noindent As usual, $p$ denotes the unique two-variable polynomial such that $p(x+y, xy) = w(x) + w(y)$.

The eagle-eyed reader may spot that there has been a shift in the usual $\mathbb{Z}/2$-grading of the matrix factorizations $N$ and $P$.  This is due to the introduction of an extra component in the oriented resolutions of the diagram defining these factorizations.  We shall not mention this kind of shift in the remainder of this paper where it may occur, but content ourselves most often with giving the factorizations explicitly on which we are focussed.

The map $\chi_0$ is given by the pair of matrices

\[ \left( \begin{array}{cc} x_1 & -x_2 \\ -1 & 1 \end{array} \right) \]

\[ \left( \begin{array}{cc} x_1 - x_2 & 0 \\ -a & 1 \end{array} \right) \]

Finally, we are in a position to define the first of our chain maps: $F$.  This comes as the pair of matrices

\[ \left( \begin{array}{c} \alpha \\ 0 \end{array} \right) , \left( \begin{array}{c} 0 \\ -\alpha  \end{array} \right)  \]

\noindent where $\alpha = (p(x_1+ x_2, x_1 x_2) - p(x_1 + x_2, x_1 x_3))/(x_1 (x_2 - x_3))$.

It is straightforward to check that this a map of matrix factorizations, and it is clearly filtered of degree $0$.  However, it is less straightforward to check that this gives a map of chain complexes; to do so we need to see that $\chi_0 F$ is homotopic to the zero map of chain complexes.  The map $\chi_0 F$ is given by the pair of matrices:

\[ \left( \begin{array}{c} x_1 \alpha \\ -\alpha \end{array} \right) ,  \left( \begin{array}{c} 0 \\ -\alpha \end{array} \right) \]

We define a homotopy $H^0 : M_0 \rightarrow P_1$, $H^1 : M_1 \rightarrow P_0$:

\[ H^0 = \left( \begin{array}{c} -1 \\ \frac{u_2 - \alpha}{x_2 - x_3} \end{array} \right),   H^1 = \left( \begin{array}{c} 1 \\ 0 \end{array} \right) \]

\noindent (a glance at the definitions of $\alpha$ and $u_2$ assures us that the second entry of $H^0$ is a polynomial).

To see that $H$ is a homotopy between $\chi_0 F$ and $0$ we need to compute that

\[ \left( \begin{array}{c} x_1 \alpha \\ -\alpha \end{array} \right) = H^1 (\pi_{23}) + B^1 H^0,  \left( \begin{array}{c} 0 \\ -\alpha \end{array} \right) = H^0 (x_2 - x_3) + B^0 H^1 \]

\noindent which is left as an exercise.

\begin{figure}
\centerline{
{
\psfrag{0}{$0$}
\psfrag{F}{$\tilde{F}$}
\psfrag{G}{$\tilde{G}$}
\psfrag{chi0}{$\chi_0$}
\psfrag{M}{$M$}
\psfrag{N}{$N$}
\psfrag{P}{$P$}
\psfrag{phi}{$\phi$}
\psfrag{x1}{$x_1$}
\psfrag{x2}{$x_2$}
\psfrag{x3}{$x_3$}
\includegraphics[height=1.5in,width=2.3in]{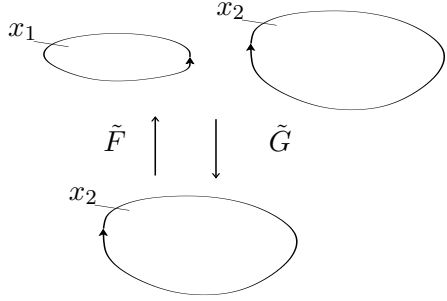}
}}
\caption{$\tilde{F}$ induced by $F$}
\label{R11closedup}
\end{figure}

We have now shown that $F$ is a degree $0$ chain map.  To complete the proof of one direction of the Reidemeister I.1 move, it remains to check that $F$ induces a map that preserves the generators of the homology of a closed link diagram.  Because of the way that the generators are defined (as the image of various elements of the homology of some disjoint circles under $\eta_0$ maps), it is enough to check that the generators are preserved by the map $\tilde{F}$ that $F$ induces when $x_2$ is identified with $x_3$ as in Figure \ref{R11closedup}.

We write $p = p(s,t)$ where $s = x_1 + x_2$ and $t = x_1 x_2$.  Looking above at the definition of the polynomial $\alpha$ that appears in the entries of the matrix $F$, we see that result of substituting $x_2$ for $x_3$ in $\alpha$ is

\begin{eqnarray*}
\alpha_{3 \rightarrow 2} &=& \frac{\partial}{\partial t} p(s,t) = \frac{\partial}{\partial t} (w_1 + w_2) \\
&=& \frac{\partial x_1}{\partial t} \partial w_1 + \frac{\partial x_2}{\partial t} \partial w_2 \\
&=& (x_2 - x_1)^{-1} \partial w_1 - (x_2 - x_1)^{-1} \partial w_2 \\
&=& - \frac{\partial w_1 - \partial w_2}{x_1 - x_2}
\end{eqnarray*}

Taking homology of the graphs in Figure \ref{R11closedup}, $\tilde{F}$ becomes multiplication by $\alpha_{3 \rightarrow 2}$

\[ \tilde{F} = \alpha_{3 \rightarrow 2} : \mathbb{C}[x_2]/(\partial w_2) \rightarrow \mathbb{C}[x_1, x_2]/(\partial w_1 , \partial w_2) \]

A general basis element, as constructed in Section 2, of the homology of the lower circle in Figure \ref{R11closedup} is written (up to non-zero multiple) as $\partial w_2 / ( x_2 - \psi )$ where $\psi$ is a root of $\partial w$.  Under $\tilde{F}$ this gets mapped to

\begin{eqnarray*}
\alpha_{3 \rightarrow 2} \frac{\partial w_2}{x_2 - \psi} &=& - \left(\frac{\partial w_1 - \partial w_2}{x_1 - x_2}\right) \left( \frac{\partial w_2}{x_2  - \psi} \right)\\
&=& - \left( \frac{\partial w_1}{x_1 - \psi} \right) \left( \frac{ \partial w_2 }{x_2 - \psi} \right)
\end{eqnarray*}

\noindent since

\[ \frac{\partial w_1 - \partial w_2}{x_1 - x_2} - \frac{\partial w_1}{x_1 - \psi} \]

\noindent has factor of $(x_2 - \psi)$ and $\partial w_2 = 0$ in the image of $\tilde{F}$.

%%Bingo.

Write $D$, $D'$ for closed link diagrams which are the same except as they differ locally as the right and left parts of Figure \ref{R11} respectively.  We have shown that if $g$ is a standard basis element of $HKh_w(D)$ then it gets taken by $\tilde{F}$ to (a non-zero multiple of) the basis element of $HKh_w(D')$ which is obtained by the same decoration of components of $D'$ with elements of $\Sigma_n$ as is $g$.

Next we need to define the map $G$ to give us a degree $0$ chain map from the upper chain complex of matrix factorizations to the lower.  In fact, we can see that $G$ just has to be a map of matrix factorizations of degree-$0$ and it will define a chain map between the two chain complexes since the only commuting square that we need to worry about will be automatically $0$ in both directions.

As a $\mathbb{C}[x_2,x_3]$-module, $\mathbb{C}[x_1,x_2,x_3]$ is isomorphic to the following direct sum

\[ \mathbb{C}[x_1 , x_2 , x_3] = \bigoplus_{i=0}^{n-2} x_1^i \mathbb{C}[x_2, x_3] \oplus \bigoplus_{i \geq 0} x_1^{i-1} (\partial w_1 - c) \mathbb{C}[x_2, x_3] \]

\noindent where we have written $c$ for the constant term in the polynomial $\partial w$.  Let $\beta$ be the $\mathbb{C}[x_2,x_3]$-module map

\[ \beta : \mathbb{C}[x_1,x_2,x_3] \rightarrow \mathbb{C}[x_2,x_3] \]

\[ 1, x_1, x_1^2, \ldots, x_1^{n-2} \mapsto 0, x_1^{i-1} (\partial w_1 -c) \mapsto (x_2 - x_3)^{i} \forall i \geq 0 \]

And let $1_{1 \rightarrow 2,3}$ be the map

\[ 1_{1 \rightarrow 2,3} : \mathbb{C}[x_1,x_2,x_3] \rightarrow \mathbb{C}[x_2,x_3] \]

\[ x_1^i \mapsto (x_2 - x_3)^i \forall i \geq 0 \]

We define the map of matrix factorizations $G$ by the pair of matrices

\[ \left( \begin{array}{cc} \beta & 0 \end{array} \right) , \left( \begin{array}{cc} 1_{1 \rightarrow 2,3} & -\beta \end{array} \right) \]

\noindent (it is straightforward to check that this defines a map of matrix factorizations filtered of degree $0$).

The map $G$ automatically defines a map of chain complexes.  Next, as we did for $\tilde{F}$, we need to compute what the induced map $\tilde{G}$ does on the homologies of the closed-up graphs in Figure \ref{R11closedup}.  Looking at the definition of $\beta$ we see that $\tilde{G}$ will take

\[ \left( \frac{\partial w_1}{x_1 - \psi} \right) \left( \frac{\partial w_2}{x_2 - \psi} \right) \in \mathbb{C}[x_1, x_2]/(\partial w_1, \partial w_2) \]

to

\[ \frac{\partial w_2}{x_2 - \psi} \in \mathbb{C}[x_2]/(\partial w_2) \]

%%Bingo.

Again, write $D$, $D'$ for closed link diagrams which are the same except as they differ locally as the right and left parts of Figure \ref{R11} respectively.  We have shown that if $g$ is a standard basis element of $HKh_w(D')$ then it gets taken by $\tilde{G}$ to (a non-zero multiple of) the basis element of $HKh_w(D)$ which is obtained by the same decoration of components of $D$ with elements of $\Sigma_n$ as is $g$.

\subsubsection{Reidemeister move I.2}

\begin{figure}
\centerline{
{
\psfrag{R11}{Reidemeister I.1}
\psfrag{02}{$0$-handle or $2$-handle}
\psfrag{1}{$1$-handle}
\psfrag{2}{Reidemeister move II.2}
\psfrag{R12}{Reidemeister move I.2}
\psfrag{N}{$N$}
\psfrag{P}{$P$}
\psfrag{phi}{$\phi$}
\psfrag{x1}{$x_1$}
\psfrag{x2}{$x_2$}
\psfrag{x3}{$x_3$}
\includegraphics[height=3in,width=3in]{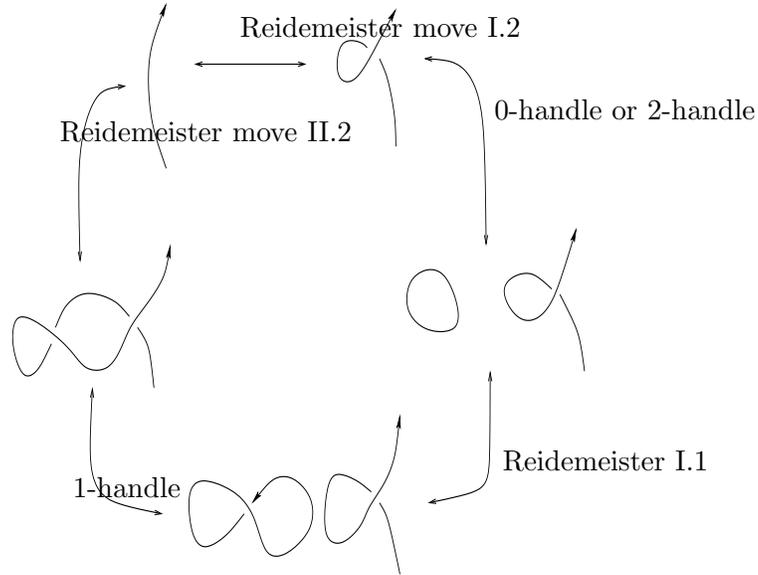}
}}
\caption{Reidemeister I.2}
\label{R12}
\end{figure}

Figure \ref{R12} decomposes the Reidemeister I.2 move into other elementary cobordisms.  From other parts of this paper it follows that this gives degree $0$ chain maps (by composition of chain maps corresponding to the elementary cobordisms used) between $CKh_w(D)$ and $CKh_w(D')$ where $D$ and $D'$ are closed link diagrams differing locally by the Reidemeister I.2 move.  These chain maps will preserve generators of $HKh_w(D)$ and $HKh_w(D)$ up to non-zero multiples just as those corresponding to the Reidemeister I.1 move did.

\subsection{Reidemeister move II}

\begin{figure}
\centerline{
{
\psfrag{0}{$0$}
\psfrag{1}{$1$}
\psfrag{2}{$2$}
\psfrag{3}{$3$}
\psfrag{4}{$4$}
\psfrag{5}{$5$}
\psfrag{6}{$6$}
\psfrag{id}{$id$}
\psfrag{chi0}{$\chi_0$}
\psfrag{chi1}{$\chi_1$}
\psfrag{phi}{$\phi$}
\includegraphics[height=2in,width=2.5in]{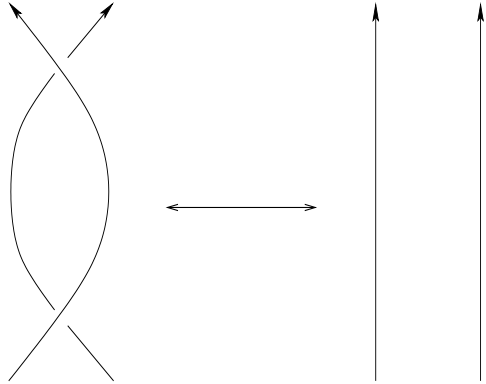}
}}
\caption{Reidemeister II.1 move}
\label{ronen}
\end{figure}

Since our link comes with an orientation, there are two cases to
compute, the first is in Figure \ref{ronen}.

\subsubsection{Reidemeister move II.1}

\begin{figure}
\centerline{
{
\psfrag{0}{$0$}
\psfrag{1}{$1$}
\psfrag{2}{$2$}
\psfrag{3}{$3$}
\psfrag{4}{$4$}
\psfrag{5}{$5$}
\psfrag{6}{$6$}
\psfrag{id}{$id$}
\psfrag{chi0}{$\chi_0$}
\psfrag{chi1}{$\chi_1$}
\psfrag{phi}{$\phi$}
\includegraphics[height=3in,width=3in]{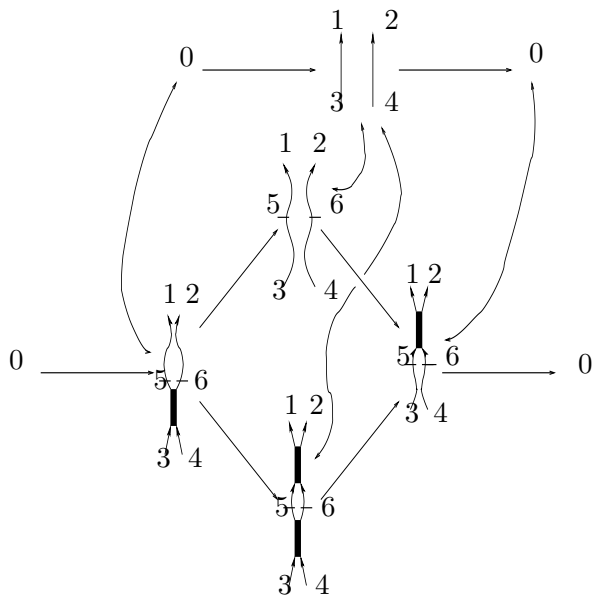}
}}
\caption{Reidemeister II.1 chain maps}
\label{R22}
\end{figure}

We start with a note about ordering of bases.  The trivalent graphs in Figure \ref{Gamma1} and Figure \ref{Gamma0} are associated to factorizations defined explicitly in Section 2.  We write either of these factorizations as

\[ \left( \begin{array}{c} M_{00} \\ M_{11} \end{array} \right) \rightarrow \left( \begin{array}{c} M_{10} \\ M_{01} \end{array} \right) \rightarrow \left( \begin{array}{c} M_{00} \\ M_{11} \end{array} \right) \]

\noindent where each $M_{ij}$ is a free module of rank 1 over the ring $\mathbb{C}[x_1,x_2,x_3,x_4]$.  The basis that we use for the tensor product of two of these factorisations, as appears in the lower half of Figure \ref{R22} for example, is

\[  \left( \begin{array}{c} M_{0000} \\ M_{0011} \\ M_{0101} \\ M_{0110} \\ M_{1001} \\ M_{1010} \\ M_{1100} \\ M_{1111} \end{array} \right) \rightarrow \left( \begin{array}{c} M_{0001} \\ M_{0010} \\ M_{0100} \\ M_{1000} \\ M_{1110} \\ M_{1101} \\ M_{1011} \\ M_{0111} \end{array} \right) \rightarrow \left( \begin{array}{c} M_{0000} \\ M_{0011} \\ M_{0101} \\ M_{0110} \\ M_{1001} \\ M_{1010} \\ M_{1100} \\ M_{1111} \end{array} \right) \]

\noindent where $M_{ijkl} = M_{ij} \otimes M_{kl}$ and $M_{ij}$ is a summand of the higher matrix factorization, $M_{kl}$ a summand of the lower.

Each $M_{ijkl}$ is a $\mathbb{C}[x_1,x_2,x_3,x_4,x_5,x_6]$ considered as a module over the ring $\mathbb{C}[x_1,x_2,x_3,x_4]$.  In what follows we have many $8 \times 8$ matrices.  Each matrix entry should be understood as a $\mathbb{C}[x_1,x_2,x_3,x_4]$-module map $\mathbb{C}[x_1,x_2,x_3,x_4,x_5,x_6] \rightarrow \mathbb{C}[x_1,x_2,x_3,x_4,x_5,x_6]$.  We will often abuse notation - for example $x_5$ appearing as a matrix entry will mean the ``multiply by $x_5$'' map, even though $x_5$ technically lies in the module and not the ring.

%% Insert drawing here.

The complex of matrix factorizations coming from the left-hand side of
Reidemeister move II.1 is in the lower half of Figure \ref{R22}.  In this diagram, as in the others in this section, where a point would normally be labelled with a variable $x_i$ we have suppressed the $x$ and just labelled with $i$.

Figure \ref{R22} shows the two complexes of matrix factorizations coming from either side of Figure \ref{ronen}.  There are also some double-headed arrows.  Each double-headed arrow is two maps of matrix factorizations (one in either direction); these maps of matrix factorizations are the component maps of two chain maps of matrix factorizations (one going each way between the higher and lower complexes).  It is our first task in this section to define these maps of matrix factorizations so that we indeed \emph{do} have two chain maps.

We also aim to define these maps so that the induced chain maps between two complexes coming from closed link diagrams differing by a single Reidemeister II.1 move preserve the generators of the homology $HKh_w$.

Two of the double-headed arrows are automatically the zero map.  The double-headed arrow going between the two factorizations with no thick edges shall just be homotopy equivalences of matrix factorizations as described in the Appendix.  For the double-headed arrow which goes between the top factorization and the doubly thick-edged factorization, we shall write down two explicit maps of matrix factorizations.

\begin{figure}
\centerline{
{
\psfrag{0}{$0$}
\psfrag{1}{$1$}
\psfrag{2}{$2$}
\psfrag{3}{$3$}
\psfrag{4}{$4$}
\psfrag{5}{$5$}
\psfrag{6}{$6$}
\psfrag{id}{$id$}
\psfrag{chi0}{$\chi_0$}
\psfrag{chi1}{$\chi_1$}
\psfrag{phi}{$\phi$}
\includegraphics[height=3in,width=3in]{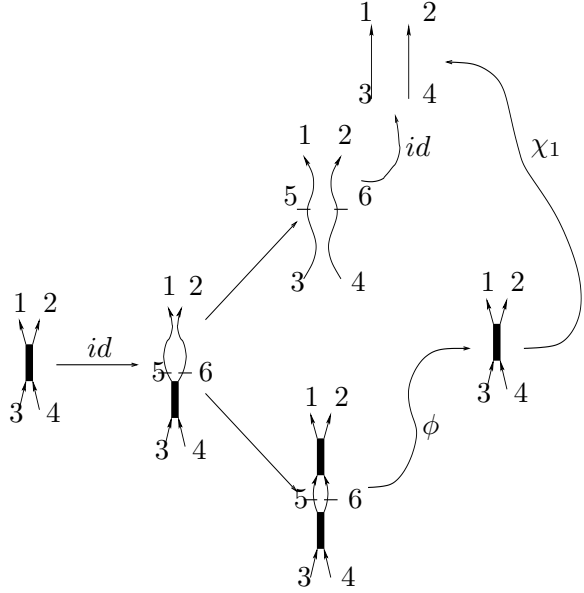}
}}
\caption{Factoring the upwards chain map}
\label{R21upwards}
\end{figure}

First we shall consider the chain map from the lower chain complex of matrix factorizations to the higher.  In Figure \ref{R21upwards} we have isolated the relevant part of Figure \ref{R22} and factored the map from the doubly thick-edged factorization into the composition of two maps.  The arrows marked $id$ are the identity in the homotopy category of matrix factorizations.  We will have a chain map upwards if we can define the map $\phi$ so that Figure \ref{R21upwards} is anti-commutative.

Using the basis conventions discussed earlier we define $\phi$ to be

\[ \Pi  = \phi_0 \circ \left( \begin{array}{cccccccc} 1 & 0 & 0 & 0 & 0 & 0 & 0 & 0 \\ 0 & 0 & 0 & 0 & 0 & 0 & 1 & 0 \end{array} \right)  \]

\[ \Pi \circ = \phi_1 \circ \left( \begin{array}{cccccccc} 0 & 0 & 0 & 1 & 0 & 0 & 0 & 0 \\ 0 & 0 & 1 & 0 & 0 & 0 & 0 & 0 \end{array} \right) \]

\noindent where $\Pi$ is the map which on each module summand $\mathbb{C}[x_1,x_2,x_3,x_4,x_5,x_6]$ looks like

\begin{eqnarray*}
\mathbb{C}[x_1,x_2,x_3,x_4,x_5,x_6] &\rightarrow&  \mathbb{C}[x_1,x_2,x_3,x_4,x_5,x_6]/(x_5+x_6-x_3-x_4,x_5x_6-x_3x_4) \\
&=& \mathbb{C}[x_1,x_2,x_3,x_4]{1} \oplus \mathbb{C}[x_1,x_2,x_3,x_4]{x_5} \\
&\rightarrow& \mathbb{C}[x_1,x_2,x_3,x_4]
\end{eqnarray*}

\noindent (the last map is projection onto the \emph{second} module summand).  The map $\phi$ is easily checked to be a map of matrix factorizations.  Next we want to see that we have defined a chain map by seeing that Figure \ref{R21upwards} commutes.

Let us compute the map in Figure \ref{R21upwards} which runs from the leftmost singly thick-edged factorization to the singly thick-edged factorization which is the target of $\phi$ (this is a composition of three maps of matrix factorizations).  We will denote this map by the pair of matrices

\[ \left( \begin{array}{cc} b^{01} & b^{02} \\ b^{03} & b^{04} \end{array} \right), \left( \begin{array}{cc} b^{11} & b^{12} \\ b^{13} & b^{14} \end{array} \right) \]

\begin{figure}
\centerline{
{
\psfrag{0}{$0$}
\psfrag{1}{$1$}
\psfrag{2}{$2$}
\psfrag{3}{$3$}
\psfrag{4}{$4$}
\psfrag{5}{$5$}
\psfrag{6}{$6$}
\psfrag{id}{$id$}
\psfrag{chi0}{$\chi_0$}
\psfrag{chi1}{$\chi_1$}
\psfrag{phi}{$\phi$}
\psfrag{f1,g1}{$f_1,g_1$}
\psfrag{f2,g2}{$f_2,g_2$}
\psfrag{f5,g5}{$f_5,g_5$}
\psfrag{f6.g6}{$f_6,g_6$}
\includegraphics[height=2in,width=2in]{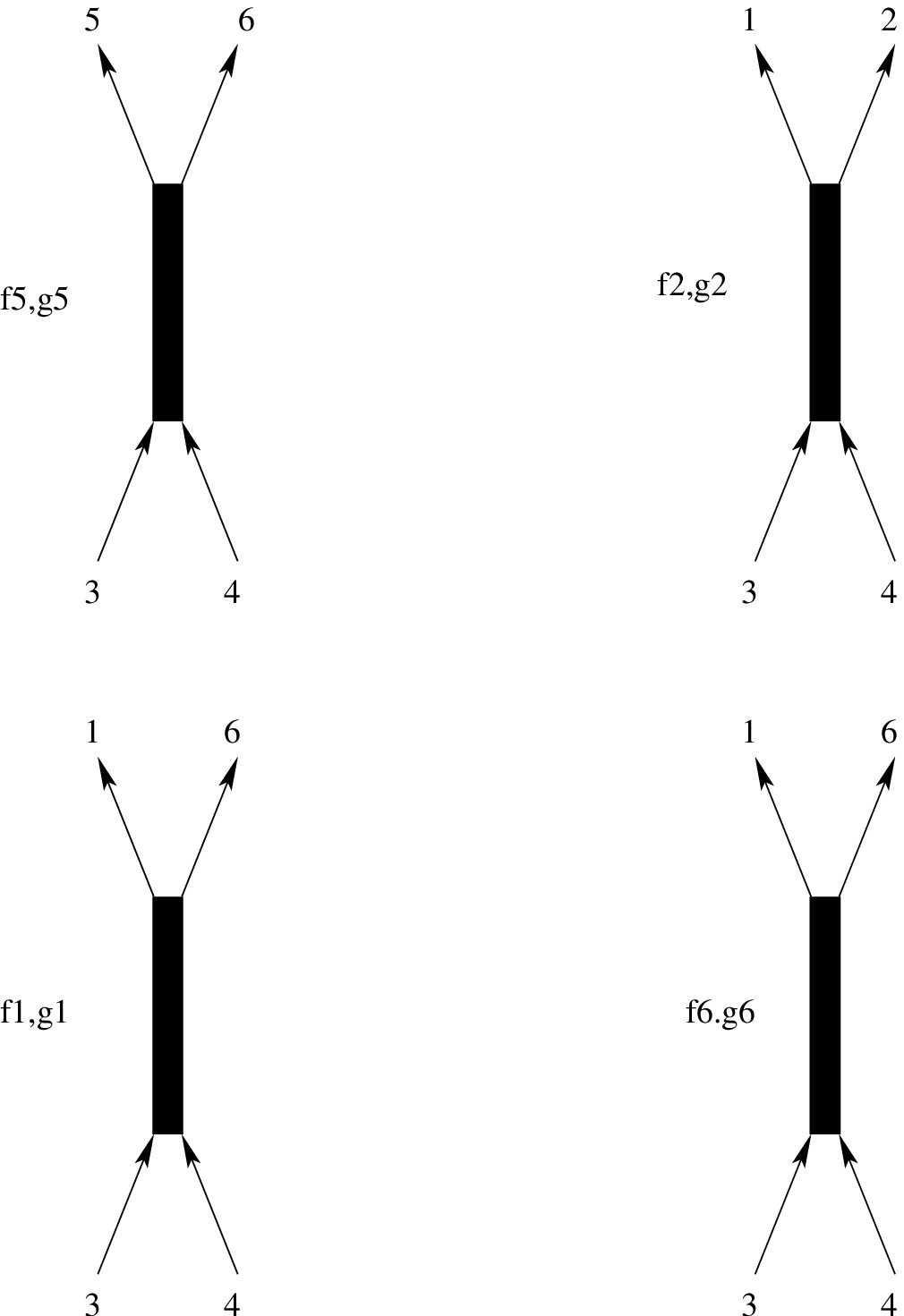}
}}
\caption{}
\label{singleboxes}
\end{figure}

We shall write $f:M^0 \rightarrow M^1$ and $g:M^1 \rightarrow M^0$ with suitable subscripts to stand for the matrix factorizations in Figure \ref{singleboxes}

Let

\[ f_{15} = (f_1 - f_5)/(x_1 - x_5), g_{15} = (g_1 - g_5)/(x_1 - x_5) \]

\[ f_{26} = (f_2 - f_6)/(x_2 - x_6), g_{26} = (g_2 - g_6)/(x_2 -x_6) \]

\begin{figure}
\centerline{
{
\psfrag{0}{$0$}
\psfrag{1}{$1$}
\psfrag{2}{$2$}
\psfrag{3}{$3$}
\psfrag{4}{$4$}
\psfrag{5}{$5$}
\psfrag{6}{$6$}
\psfrag{id}{$id$}
\psfrag{chi0}{$\chi_0$}
\psfrag{chi1}{$\chi_1$}
\psfrag{phi}{$\phi$}
\psfrag{f1,g1}{$f_1,g_1$}
\psfrag{f2,g2}{$f_2,g_2$}
\psfrag{f5,g5}{$f_5,g_5$}
\psfrag{f6.g6}{$f_6,g_6$}
\includegraphics[height=1.4in,width=2in]{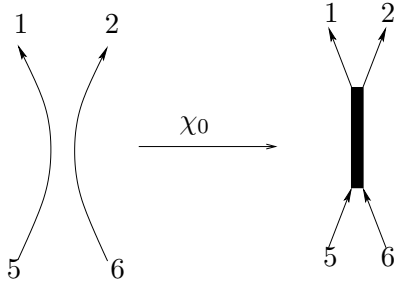}
}}
\caption{$\chi_0$ map}
\label{defn_of_a1}
\end{figure}

Then we have

\begin{eqnarray*}
 \left( \begin{array}{cc} b^{01} & b^{02} \\ b^{03} & b^{04} \end{array} \right) &=& \phi_0 \circ \left( \begin{array}{cccccccc} x_5 - x_2 & 0 & 0 & 0 & 0 & 0 & 0 & 0 \\ 0 & x_5 - x_2 & 0 & 0 & 0 & 0 & 0 & 0 \\ 0 & 0 & 1 & 0  & -1 & 0 & 0 & 0 \\ 0  & 0 & 0 & 1 & 0 & -1 & 0 & 0 \\ 0  & 0 & -x_2 & 0 & x_5 & 0 & 0 & 0 \\ 0 & 0 & 0 & -x_2 & 0 & x_5 & 0 & 0 \\ -a_1 & 0 & 0 & 0 & 0 & 0 & 1 & 0 \\ 0 & -a_1 & 0 & 0 & 0 & 0 & 0 & 1 \end{array} \right) \\
&\circ& \left( \begin{array}{cccccccc} 1 & 0 & 0 & 0 & 0 & 0 & 0 & 0 \\ 0 & 1 & 0 & 0 & 0 & 0 & 0 & 0 \\ 0 & 0 & 0 & 1 & 0 & 0 & 0 & 0 \\ 0 & 0 & 1 & 0 & 0 & 0 & 0 & 0 \\ 0 & 0 & 0 & 0 & 0 & 0 & 0 & 1 \\ 0 & 0 & 0 & 0 & 0 & 0 & 1 & 0 \\ 0 & 0 & 0 & 0 & 1 & 0 & 0 & 0 \\ 0 & 0 & 0 & 0 & 0 & 1 & 0 & 0 \end{array} \right) \left( \begin{array}{c} id_{R^2} \\ f_{26} \\ -g_{15} \circ f_{26} \\ f_{15} \end{array} \right) 
\end{eqnarray*}

\begin{eqnarray*}
\left( \begin{array}{cc} b^{11} & b^{12} \\ b^{13} & b^{14} \end{array} \right) &=& \phi_1 \circ \left( \begin{array}{cccccccc} x_5 - x_2 & 0 & 0 & 0 & 0 & 0 & 0 & 0 \\ 0 & x_5 - x_2 & 0 & 0 & 0 & 0 & 0 & 0 \\ 0 & 0 & 1 & -1 & 0 & 0 & 0 & 0 \\ 0 & 0 & -x_2 & x_5 & 0 & 0 & 0 & 0 \\ 0 & -a_1 & 0 & 0 & 1 & 0 & 0 & 0 \\ -a_1 & 0 & 0 & 0 & 0 & 1 & 0 & 0 \\ 0 & 0 & 0 & 0 & 0 & 0 & x_5 & -x_2 \\ 0 & 0 & 0 & 0 & 0 & 0 & -1 & 1 \end{array} \right) \\
&\circ& \left( \begin{array}{cccccccc} 0 & 0 & 0 & 0 & 0 & 0 & 0 & 1 \\ 0 & 0 & 0 & 0 & 0 & 0 & 1 & 0 \\ 0 & 0 & 0 & 0 & 1 & 0 & 0 & 0 \\ 1 & 0 & 0 & 0 & 0 & 0 & 0 & 0 \\ 0 & 0 & 1 & 0 & 0 & 0 & 0 & 0 \\ 0 & 0 & 0 & 1 & 0 & 0 & 0 & 0 \\ 0 & 1 & 0 & 0 & 0 & 0 & 0 & 0 \\ 0 & 0 & 0 & 0 & 0 & 1 & 0 & 0 \end{array} \right) \left( \begin{array}{c} g_{15} \\ -f_{15} \circ g_{26} \\ g_{26} \\ id_{R^2} \end{array} \right)
\end{eqnarray*}

\noindent where $a_1$ is the usual polynomial appearing in the $\chi_0$ map in Figure \ref{defn_of_a1}.  In our previous notation (with now $x_5 + x_6$ and $x_5x_6$ replacing $x_3 + x_4$ and $x_3x_4$ respectively in the definitions of $u_1$ and $u_2$) we have

\[ a_1 = \frac{u_1 + x_5u_2 - \pi_{26}}{x_1 - x_5} = - \frac{u_1 + x_2u_2 - \pi_{15}}{x_2 - x_6} \]

Multiplying out the matrices:

\begin{eqnarray*}
\left( \begin{array}{cc} b^{01} & b^{02} \\ b^{03} & b^{04} \end{array} \right) &=& \Pi \circ \left( \begin{array}{cccccccc} x_5 - x_2 & 0 & 0 &  0 & 0 & 0 & 0 & 0 \\ -a_1 & 0 & 0 & 0 & 1 & 0 & 0 & 0 \end{array} \right) \left( \begin{array}{c} id_{R^2} \\ f_{26} \\ - g_{15}f_{26} \\ f_{15} \end{array} \right) \\
&=& \Pi \circ \left( \left( \begin{array}{cc} x_5 - x_2 & 0 \\ -a_1 & 0 \end{array} \right) - \left( \begin{array}{cc} 0 & 0 \\ 1 & 0 \end{array} \right) g_{15}f_{26} \right)
\end{eqnarray*}

Let us write

\[ f_2 = \left( \begin{array}{cc} v_{12} & x_3x_4 - x_1x_2 \\ v_{22} & x_1 + x_2 - x_3 - x_4 \end{array} \right), f_6 = \left( \begin{array}{cc} v_{16} & x_3x_4 - x_1x_2 \\ v_{26} & x_1 + x_2 - x_3 - x_4 \end{array} \right) \]

\noindent so

\[ f_{26} = \left( \begin{array}{cc} \frac{v_{12} - v_{16}}{x_2 - x_6} & -x_1 \\ \frac{v_{22} - v_{26}}{x_2 - x_6} & 1 \end{array} \right) \]

Now

\[ g_{15} = \left( \begin{array}{cc} 1 & x_6 \\ . & . \end{array} \right) \]

\noindent so that we see immediately

\[ b^{01} = 1, b^{02} = 0, b^{04} = 1 \]

It remains to compute $b^{03}$:

\begin{eqnarray*}
b^{03} &=& \Pi\frac{1}{x_2 - x_6}[(u_1 + x_2u_2 - \pi_{15}) - (v_{12} - v_{16}) - (x_6(v_{22} - v_{26}))] \\
&=& \Pi \left( \frac{1}{x_2 - x_6} (u_1 + x_2 u_2 - \pi_{15} - v_{12} + v_{16} - x_2 v_{22} + x_6 v_{26}) + v_{22} \right) \\
&=& \Pi \left( \frac{1}{x_2 - x_6} (\pi_{15} + v_{16} + x_6v_{26}) + v_{22} \right)
\end{eqnarray*}

\noindent(the last equality follows from the action of $\Pi$)

Now $v_{22}$ certainly gets killed by $\Pi$ and

\[ \frac{-\pi_{15} + v_{16} + x_6v_{26}}{x_2 - x_6} \]

\noindent is a polynomial with no term involving $x_2$ in the numerator and hence also gets killed by $\Pi$.  Thus we have shown that

\[b^{03} = 0\]

Now we work on $b^{11},b^{12},b^{13},b^{14}$

\begin{eqnarray*}
\left( \begin{array}{cc} b^{11} & b^{12} \\ b^{13} & b^{14} \end{array} \right) &=& \Pi \circ \left( \begin{array}{cccccccc} x_5 & 0 & 0 & 0 & -x_2 & 0 & 0 & 0 \\ -1 & 0 & 0 & 0 & 1 & 0 & 0 & 0 \end{array} \right) \left( \begin{array}{c} g_{15} \\ -f_{15}g_{26} \\ g_{26} \\ id_{R^2} \end{array} \right) \\
&=& \Pi \circ \left( \left( \begin{array}{cc} x_5 & 0 \\ -1 & 0 \end{array} \right) g_{15} + \left( \begin{array}{cc} -x_2 & 0 \\ 1 & 0 \end{array} \right) g_{26} \right)
\end{eqnarray*}

We compute that

\[ g_{15} = \left( \begin{array}{cc} 1 & x_6 \\ . & . \end{array} \right) , g_{26} = \left( \begin{array}{cc} 1 & x_1 \\ . & . \end{array} \right) \]

Hence

\[ \left( \begin{array}{cc} b^{11} & b^{12} \\ b^{13} & b^{14} \end{array} \right) = \Pi \circ \left( \begin{array}{cc} x_5 - x_2 & x_5x_6 - x_1x_2 \\ -1 + 1 & -x_6 + x_1 \end{array} \right) = \left( \begin{array}{cc} 1 & 0 \\ 0 & 1 \end{array} \right) \]

Now it is clear that Figure \ref{R21upwards} commutes (see theorem \ref{losingmarks} in the appendix if you need more justification).  Hence we have defined a chain map.

It remains to check that the chain map we have defined preserves the generators of the homology.  Recall that the first step to define a generator involves decorating the link components with roots of $\partial w$.  We then resolve the diagram by adding a thick edge at each crossing where the roots disagree and adding the oriented resolution at each crossing where the roots agree.

If the decoration of the two strands of the Reidemeister II move is by the same root then we have nothing to show since the relevant component of the chain map is homotopic to the identity.

%%(And because $\chi_1 \phi$ commutes with the action of, say, $x_1$ and $x_2$, the induced map on the homology of the correctly resolved link diagram shall preserve the eigenspaces of $x_1$ and $x_2$.  The $\alpha_i$-eigenspace is the $1$-dimensional vector space generated by the basis element corresponding to the root $\alpha_i$.)

\begin{figure}
\centerline{
{
\psfrag{0}{$0$}
\psfrag{1}{$1$}
\psfrag{2}{$2$}
\psfrag{3}{$3$}
\psfrag{4}{$4$}
\psfrag{5}{$5$}
\psfrag{6}{$6$}
\psfrag{id}{$id$}
\psfrag{chi0}{$\chi_0$}
\psfrag{chi1}{$\chi_1$}
\psfrag{phi}{$\phi$}
\psfrag{chi1phi}{$\chi_1 \phi$}
\includegraphics[height=2in,width=2in]{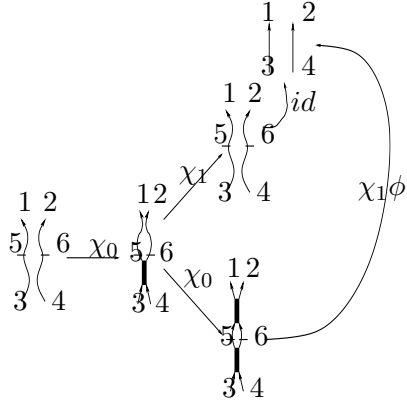}
}}
\caption{Preservation of generators under Reidemeister II.1}
\label{R21upwardspreserve}
\end{figure}

In the case that the decoration is by two different roots then we use the result in \ref{the eta map} and refer to Figure \ref{R21upwardspreserve}.  Suppose that the strand with left endpoints is decorated by the root $\alpha_1$ and the strand with right endpoints is decorated by the root $\alpha_2$.

Recall that to produce the generator we have to start with a resolution which has the homology of a set of disjoint circles and then push forward a generator (specified by our choice of root decoration) of that homology by $\eta_0$ maps.  The appendix \ref{the eta map} tells us that this pushing forward by $\eta_0$ maps gives us the same generator, up to sign, as pushing forward by two $\chi_0$ maps at the site of the Reidemeister II move and by $\eta_0$ maps elsewhere.

In other words, we can create the generator by pushing forward the relevant generator of the homology of the diagram looking locally like the far left resolution in Figure \ref{R21upwardspreserve} by the two $\chi_0$ maps indicated.  Since we have shown that this is a commutative diagram, the image of the generator (in the homology of the diagram looking locally like the top of Figure \ref{R21upwardspreserve}) under the chain map that we gave is the same as first applying $\chi_0$ and then $\chi_1$ to the corresponding generator of the homology of the diagram looking locally like the far left resolution.  Earlier, we have shown that this is the same up to homotopy as the map induced by multiplication by $x_6 - x_3$.  This is the same as multiplying our generator by the non-zero number $\alpha_2 - \alpha_1$.  Hence our generator is preserved up to non-zero multiple.

\begin{figure}
\centerline{
{
\psfrag{0}{$0$}
\psfrag{1}{$1$}
\psfrag{2}{$2$}
\psfrag{3}{$3$}
\psfrag{4}{$4$}
\psfrag{5}{$5$}
\psfrag{6}{$6$}
\psfrag{id}{$id$}
\psfrag{chi0}{$\chi_0$}
\psfrag{chi1}{$\chi_1$}
\psfrag{phi}{$\phi$}
\psfrag{f1,g1}{$f_1,g_1$}
\psfrag{f2,g2}{$f_2,g_2$}
\psfrag{f5,g5}{$f_5,g_5$}
\psfrag{f6.g6}{$f_6,g_6$}
\includegraphics[height=2in,width=2in]{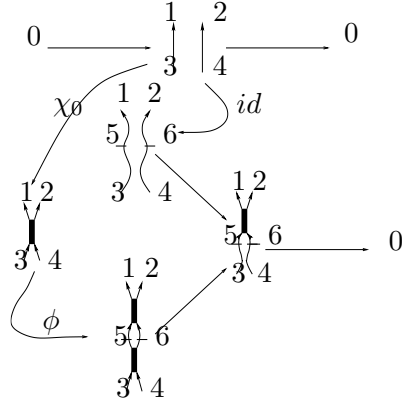}
}}
\caption{Factoring the downwards chain map}
\label{R21downwards}
\end{figure}

Now we shall consider the chain map going from the higher chain complex to the lower.  In Figure \ref{R21downwards} we have isolated the relevant part of the diagram and factored the map to the doubly thick-edged factorization into the compositon of two maps.   The arrows marked $id$ are the identity in the homotopy category of matrix factorizations.  We shall have a chain map if we can give the map $\phi$ so that the diagram anti-commutes.

Let $F$ and $G$ be the $2 \times 2$ matrices of the factorization in Figure \ref{thickfactorization2}.  The entries of $F$ and $G$ are polynomials in $x_1+x_2, x_3+x_4,x_1x_2,x_3x_4$, and we shall write $F_{5+6}$ and $G_{5+6}$ for the matrices obtained from $F$ and $G$ by replacing $x_3 + x_4$ with $x_5 + x_6$ and $F_{5+6}^{56}$ and $G_{5+6}^{56}$ for the matrices obtained from $F$ and $G$ by replacing $x_3 + x_4$ and $x_3x_4$ with $x_5 + x_6$ and $x_5x_6$ respectively.

Writing $P(x+y,xy) = w(x) + w(y)$ we define

\[ u_1(x,y,z) = \frac{P(x,z) - P(y,z)}{x-y} \]

Using the basis-ordering conventions given at the start of this section we define $\phi$ to be

\begin{eqnarray*}
\phi_0 &=& \left( \begin{array}{cccccccc} 1 & 0 & 0 & 0 & 0 & 0 & 0 & 0 \\ 0 &0 & 0 & 0 & 0 & 0 & 1 & 0 \\ 0 & 0 & 0 & 0 & 0 & 1 & 0 & 0 \\ 0 & 0 & 0 & 1 & 0 & 0 & 0 & 0 \\ 0 & 0 & 0 & 0 & 1 & 0 & 0 & 0 \\ 0 & 0 & 1 & 0 & 0 & 0 & 0 & 0 \\ 0 & 1 & 0 & 0 & 0 & 0 & 0 & 0 \\ 0 & 0 & 0 & 0 & 0 & 0 & 0 & 1 \end{array} \right)  \left( \begin{array}{c} Id_{R^2} \\ X_1 \\ X_2 \\ X_3 + X_4 \end{array} \right)
\end{eqnarray*}

\begin{eqnarray*}
\phi_1 &=& \left( \begin{array}{cccccccc} 0 & 0 & 0 & 0 & 1 & 0 & 0 & 0 \\ 0 & 0 & 1 & 0 & 0 & 0 & 0 & 0 \\ 0 & 1 & 0 & 0 & 0 & 0 & 0 & 0 \\ 1 & 0 & 0 & 0 & 0 & 0 & 0 & 0 \\ 0 & 0 & 0 & 1 & 0 & 0 & 0 & 0 \\ 0 & 0 & 0 & 0 & 0 & 1 & 0 & 0 \\ 0 & 0 & 0 & 0 & 0 & 0 & 1 & 0 \\ 0 & 0 & 0 & 0 & 0 & 0 & 0 & 1 \end{array} \right) \left( \begin{array}{c} Id_{R^2} \\ Y_1 \\ Y_2 \\ Y_3 + Y_4 \end{array} \right)
\end{eqnarray*}

\noindent where

\begin{eqnarray*}
X_1 &=&  \frac{F-F_{5+6}}{x_3 + x_4 - x_5 - x_6} , X_2 =  \frac{F_{5+6} - F_{5+6}^{56}}{x_3x_4 - x_5x_6},\\
X_3 &=&  \frac{u_1(x_5+x_6, x_3 + x_4, x_3x_4) - u_1(x_5 + x_6,x_3 + x_4,x_5 x_6)}{x_3x_4 - x_5x_6}id_{R^2} \rm{,}\\
X_4 &=&  \frac{G_{5+6} - G_{5+6}^{56}}{x_3x_4 - x_5x_6} \frac{F - F_{5+6}}{x_3 + x_4 - x_5 - x_6} \\
Y_1 &=& \frac{G - G_{5+6}}{x_3 + x_4 - x_5 - x_6}, Y_2 =  \frac{G_{5+6} - G_{5+6}^{56}}{x_3x_4 - x_5x_6},\\
Y_3 &=&  -\frac{u_1(x_5 + x_6, x_3 + x_4, x_3x_4) - u_1(x_5 + x_6, x_3 + x_4, x_5 x_6)}{x_3x_4 - x_5x_6}id_{R^2} \rm{,} \\
Y_4 &=&  \frac{F_{5+6} - F_{5+6}^{56}}{x_3x_4 - x_5x_6}\frac{G - G_{5+6}}{x_3 + x_4 - x_5 - x_6}
\end{eqnarray*}

\noindent It is a simple, if somewhat laborious, matter to check that this is a well-defined map of matrix factorizations.  We now need to check that we have a chain map by seeing that Figure \ref{R21downwards} commutes.

If we let $\Pi$ be the map which on each module summand looks like

\begin{eqnarray*}
\Pi &:& \mathbb{C}[x_1,x_2,x_3,x_4,x_5,x_6] \rightarrow \mathbb{C}[x_1,x_2,x_3,x_4,x_5,x_6]/(x_5 - x_3, x_6 - x_4) \\
&=& \mathbb{C}[x_1,x_2,x_3,x_4]
\end{eqnarray*}

\noindent then the components of the map from the leftmost singly thick-edged factorization to the rightmost singly thick-edged factorization are

\begin{eqnarray*}
\Pi &\circ& \left( \begin{array}{cccccccc} 1 & 0 & 0 & 0 & 0 & 0 & 0 & 0 \\ 0 & 0 & 0 & 0 & 0 & 0 & 1 & 0 \end{array} \right) \left( \begin{array}{cccccccc} 1 & 0 & 0 & 0 & 0 & 0 & 0 & 0 \\ -a_2 & a_3 & 0 & 0 & 0 & 0 & 0 & 0 \\ 0 & 0 & x_3 & 1 & 0 & 0 & 0 & 0 \\ 0 & 0 & x_6 & 1 & 0 & 0 & 0 & 0 \\ 0 & 0 & 0 & 0 & x_3 & 1 & 0 & 0 \\ 0 & 0 & 0 & 0 & x_6 & 1 & 0 & 0 \\ 0 & 0 & 0 & 0 & 0 & 0 & 1 & 0 \\ 0 & 0 & 0 & 0 & 0 & 0 & -a_2 & a_3 \end{array} \right) \\
&\circ& \left( \begin{array}{cccccccc} 1 & 0 & 0 & 0 & 0 & 0 & 0 & 0 \\ 0 &0 & 0 & 0 & 0 & 0 & 1 & 0 \\ 0 & 0 & 0 & 0 & 0 & 1 & 0 & 0 \\ 0 & 0 & 0 & 1 & 0 & 0 & 0 & 0 \\ 0 & 0 & 0 & 0 & 1 & 0 & 0 & 0 \\ 0 & 0 & 1 & 0 & 0 & 0 & 0 & 0 \\ 0 & 1 & 0 & 0 & 0 & 0 & 0 & 0 \\ 0 & 0 & 0 & 0 & 0 & 0 & 0 & 1 \end{array} \right)  \left( \begin{array}{c} Id_{R^2} \\ X_1 \\ X_2 \\ X_3 + X_4 \end{array} \right) = \left( \begin{array}{cc} 1 & 0 \\ 0 & 1 \end{array} \right)
\end{eqnarray*}

\noindent and

\begin{eqnarray*}
\Pi &\circ& \left( \begin{array}{cccccccc} 0 & 0 & 0 & 1 & 0 & 0 & 0 & 0 \\ 0 & 0 & 1 & 0 & 0 & 0 & 0 & 0 \end{array} \right) \left( \begin{array}{cccccccc} x_3 & 1 & 0 & 0 & 0 & 0 & 0 & 0 \\ x_6 & 1 & 0 & 0 & 0 & 0 & 0 & 0 \\ 0 & 0 & 1 & 0 & 0 & 0 & 0 & 0 \\ 0 & 0 & 0 & 1 & 0 & 0 & 0 & 0 \\ 0 & 0 & 0 & 0 & 1 & x_6 & 0 & 0 \\ 0 & 0 & 0 & 0 & 1 & x_3 & 0 & 0 \\ 0 & 0 & 0 & -a_2 & 0 & 0 & a_3 & 0 \\ 0 & 0 & -a_2 & 0 & 0 & 0 & 0 & a_3 \end{array} \right) \\
&\circ& \left( \begin{array}{cccccccc} 0 & 0 & 0 & 0 & 1 & 0 & 0 & 0 \\ 0 & 0 & 1 & 0 & 0 & 0 & 0 & 0 \\ 0 & 1 & 0 & 0 & 0 & 0 & 0 & 0 \\ 1 & 0 & 0 & 0 & 0 & 0 & 0 & 0 \\ 0 & 0 & 0 & 1 & 0 & 0 & 0 & 0 \\ 0 & 0 & 0 & 0 & 0 & 1 & 0 & 0 \\ 0 & 0 & 0 & 0 & 0 & 0 & 1 & 0 \\ 0 & 0 & 0 & 0 & 0 & 0 & 0 & 1 \end{array} \right)  \left( \begin{array}{c} Id_{R^2} \\ Y_1 \\ Y_2 \\ Y_3 + Y_4 \end{array} \right) = \left( \begin{array}{cc} 1 & 0 \\ 0 & 1 \end{array} \right)
\end{eqnarray*}

Hence we see, just as before, that we have indeed defined a chain map.

We need now to check that this chain map, so defined, preserves the generators of the homology.  Again, the only case that is not immediate is that when the roots of $\partial w$ decorating the two strands of the Reidemeister move are distinct.  Let us suppose that the strand with ends on the left is decorated by the root $\alpha_1$ and the strand with ends on the right is decorated by the root $\alpha_2$.

\begin{figure}
\centerline{
{
\psfrag{0}{$0$}
\psfrag{1}{$1$}
\psfrag{2}{$2$}
\psfrag{3}{$3$}
\psfrag{4}{$4$}
\psfrag{5}{$5$}
\psfrag{6}{$6$}
\psfrag{id}{$id$}
\psfrag{chi0}{$\chi_0$}
\psfrag{chi1}{$\chi_1$}
\psfrag{phi}{$\phi$}
\psfrag{phichi0}{$\phi \chi_0$}
\includegraphics[height=2in,width=2in]{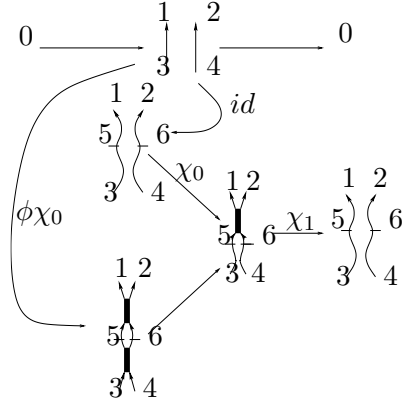}
}}
\caption{Preservation of generators under Reidemeister $2.1$}
\label{R21downwardspreserve}
\end{figure}

Suppose that $g$ is the corresponding generator in $H(\Gamma)$ where $\Gamma$ looks locally like the lowest trivalent graph in Figure \ref{R21downwardspreserve} and $H(\Gamma)$ is a chain group summand of $CKh_w$.  The generator $g$ was obtained by pushing forward $g' \in H(\Gamma')$ by $\eta_0$ maps where $\Gamma$ is a union of disjoint circles.  In \cite{Gornik}, Gornik shows that applying $\eta_1$ to each thick edge of $\Gamma$ results in a non-zero multiple of $g'$.

In Section \ref{the eta map} it is shown that applying two $\eta_1$ maps at the site of the Reidemeister II.1 move is the same up to sign as applying two $\chi_1$ maps.  Since Figure \ref{R21downwardspreserve} commutes we see that applying two $\eta_1$ maps to $g$ is the same as applying $\chi_1 \chi_0$ to $g'$.  And $\chi_1 \chi_0$ is homotopic to the map induced by multiplication by $x_2 - x_5$.  This acts the same on $g'$ as multiplication by $\alpha_2 - \alpha_1$ which is a non-zero number.  Hence our generator is preserved up to non-zero multiple.

\subsubsection{Reidemeister move II.2}

\begin{figure}
\centerline{
{
\psfrag{0handle}{$0$-handle}
\psfrag{1handle}{$1$-handle}
\psfrag{2handle}{$2$-handle}
\psfrag{Rii}{Reidemeister II.2}
\psfrag{Rii.i}{Reidemeister II.1}
\includegraphics[height=2.9in,width=4in]{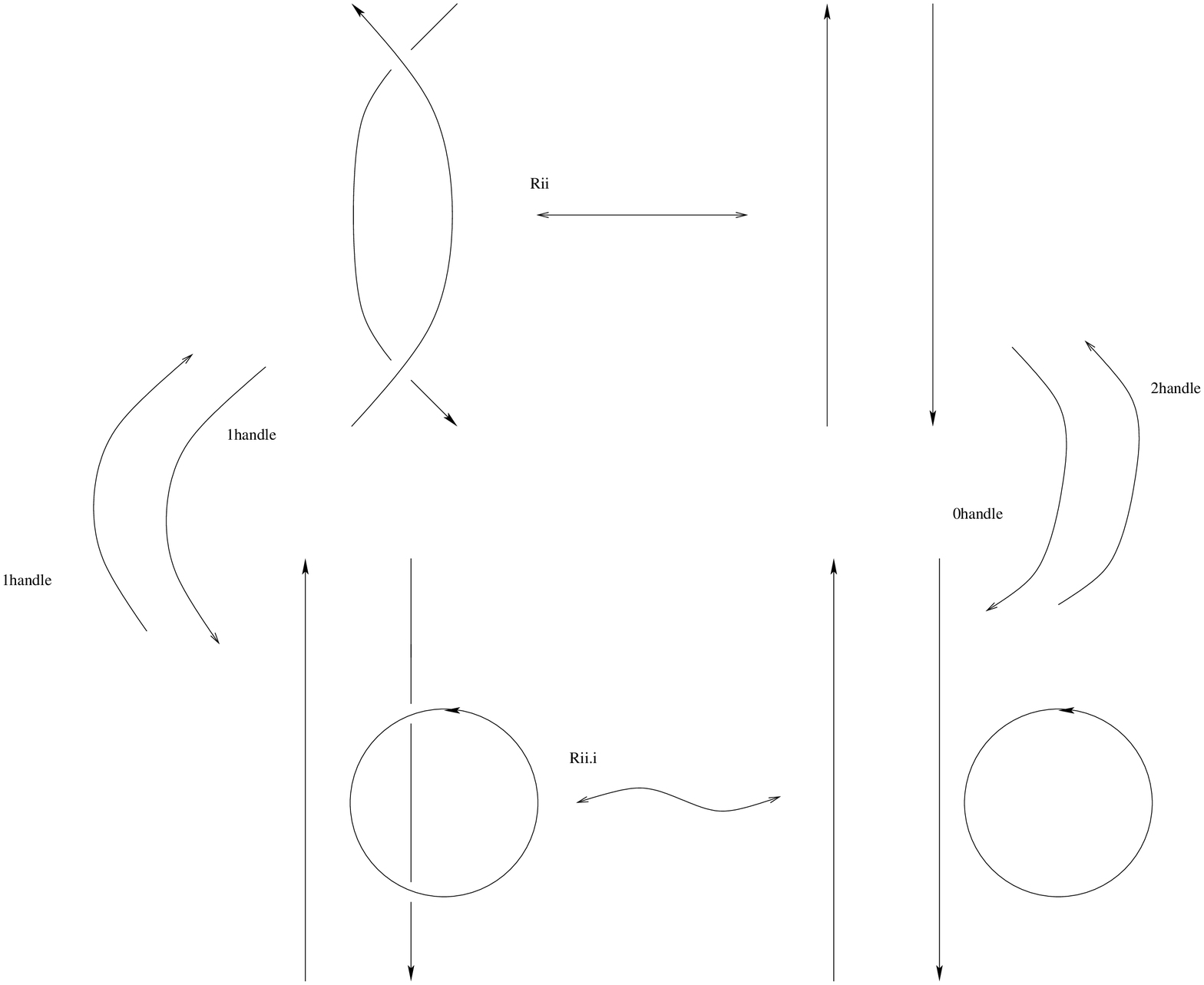}
}}
\caption{Decomposition of Reidemeister II.2 move}
\label{R22dodge}
\end{figure}

Figure \ref{R22dodge} decomposes the Reidemeister II.2 move into other elementary cobordisms.  It follows that this gives degree $0$ chain maps (by composition of chain maps corresponding to the elementary cobordisms used) between $CKh_w(D)$ and $CKh_w(D')$ where $D$ and $D'$ are closed link diagrams differing locally by the Reidemeister II.2 move.  These chain maps will preserve generators of $HKh_w(D)$ and $HKh_w(D)$ up to non-zero multiples.

%%In Figure \ref{R22dodge} we have decomposed the Reidemeister II.2 move into a product of three other elementary cobordisms.  Taking the product of the corresponding chain maps gives degree $0$ chain maps between the two chain complexes corresponding to the two sides of the Reidemeister II.2.  This chain map preserves the generators of the homology.

\section{Computation for a slice surface of a knot}

We start with a knot diagram and a presentation of a slice surface for the knot which begins with the given diagram.  We would like to see that the map induced by the presentation between the homologies of the knot and of the unknot is an isomorphism of a filtered degree depending on the genus on the surface.  This will then give us our main theorem.

This is not precisely what we do, however.  First we perturb our presentation so that we avoid having to prove similar theorems for the Reidemeister III move as we have already done for the I and II moves and for the handle attachments.  Essentially we show that it is possible to delay any Reidemeister III moves until the end of the presentation, which makes our degree computations very easy.

\subsection{Avoiding Reidemeister III}

In this section we detail a topological argument that we use to avoid having to do similar calculations to those in the previous sections.  The idea in a nutshell is that in describing a smooth surface $\Sigma \hookrightarrow \mathbb{R}^4$ as a composition of elementary cobordisms we can postpone the Reidemeister III moves until after all the $0$- and $1$-handles have been added.  This is not the whole story; we shall also be postponing some of the $1$-handle moves, and some of Reidemeister I and II moves as well.

Suppose that $S$ is a composition of elementary cobordisms representing a connected surface $\Sigma$

\[ \Sigma \hookrightarrow \mathbb{R}^3 \times [0,1] \]

\[ \partial \Sigma = K \cup U , \, K \hookrightarrow \mathbb{R}^3 \times \{ 0 \} , \, U \hookrightarrow \mathbb{R}^3 \times \{ 1 \} \]

\noindent where $S$ begins with a diagram $D$ of the knot $K$ and ends with the $0$-crossing diagram of the unknot $U$.

We can assume that the $0$-handle additions of $S$ take place before any other elementary cobordism (since a $0$-handle addition at any point of the cobordism can be replaced by a $0$-handle addition at a distance from the rest of the diagram at the start of the cobordism and then later sliding it over by Reidemeister II moves to the point where it is required).  Similarly, we can assume that the $2$-handle additions take place after every other elementary cobordism.

Suppose that $S$ has $h$ $0$-handle additions at the beginning, so that immediately after the $0$-handle additions we have a diagram of a $(h+1)$-component link.  Since $\Sigma$ is connected we can make a list of $h$ ``fusing'' $1$-handle additions in $S$ that connect all the components together.  We shall be wanting to postpone all the Reidemeister III moves and all the other $1$-handle additions until after these $h$ $1$-handles have been added, giving us a ($1$-component) knot diagram.

\begin{figure}
\centerline{
{
\psfrag{=}{$=$}
\psfrag{Di+1}{$D_{i+1}$}
\psfrag{tDi}{$\tilde{D_i}$}
\psfrag{tDi+1}{$\tilde{D_{i+1}}$}
\psfrag{Rii.i}{Reidemeister II.1}
\includegraphics[height=2in,width=4in]{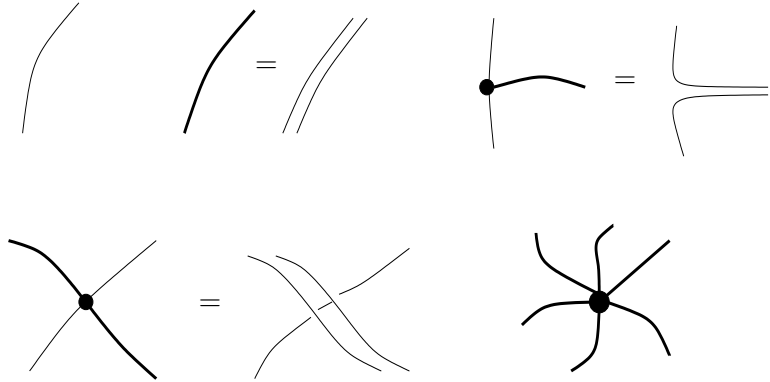}
}}
\caption{Local bestiary}
\label{localbestiary}
\end{figure}

The composition  of elementary cobordisms $S$ consists of a list of diagrams $D_1, D_2, D_3, \ldots, D_d$ where each successive diagram differs from the previous diagram by a specified elementary cobordism.  Forgetting about the $0$- and $2$-handle attachments, we can assume that $D_1$ is $D$ union $h$ $0$-crossing unknots, and $D_d$ is the union of $k+1$ $0$-crossing unknots.

\begin{figure}
\centerline{
{
\psfrag{Di}{$D_i$}
\psfrag{Di+1}{$D_{i+1}$}
\psfrag{tDi}{$\tilde{D_i}$}
\psfrag{tDi+1}{$\tilde{D_{i+1}}$}
\psfrag{Rii.i}{Reidemeister II.1}
\includegraphics[height=5in,width=4.5in]{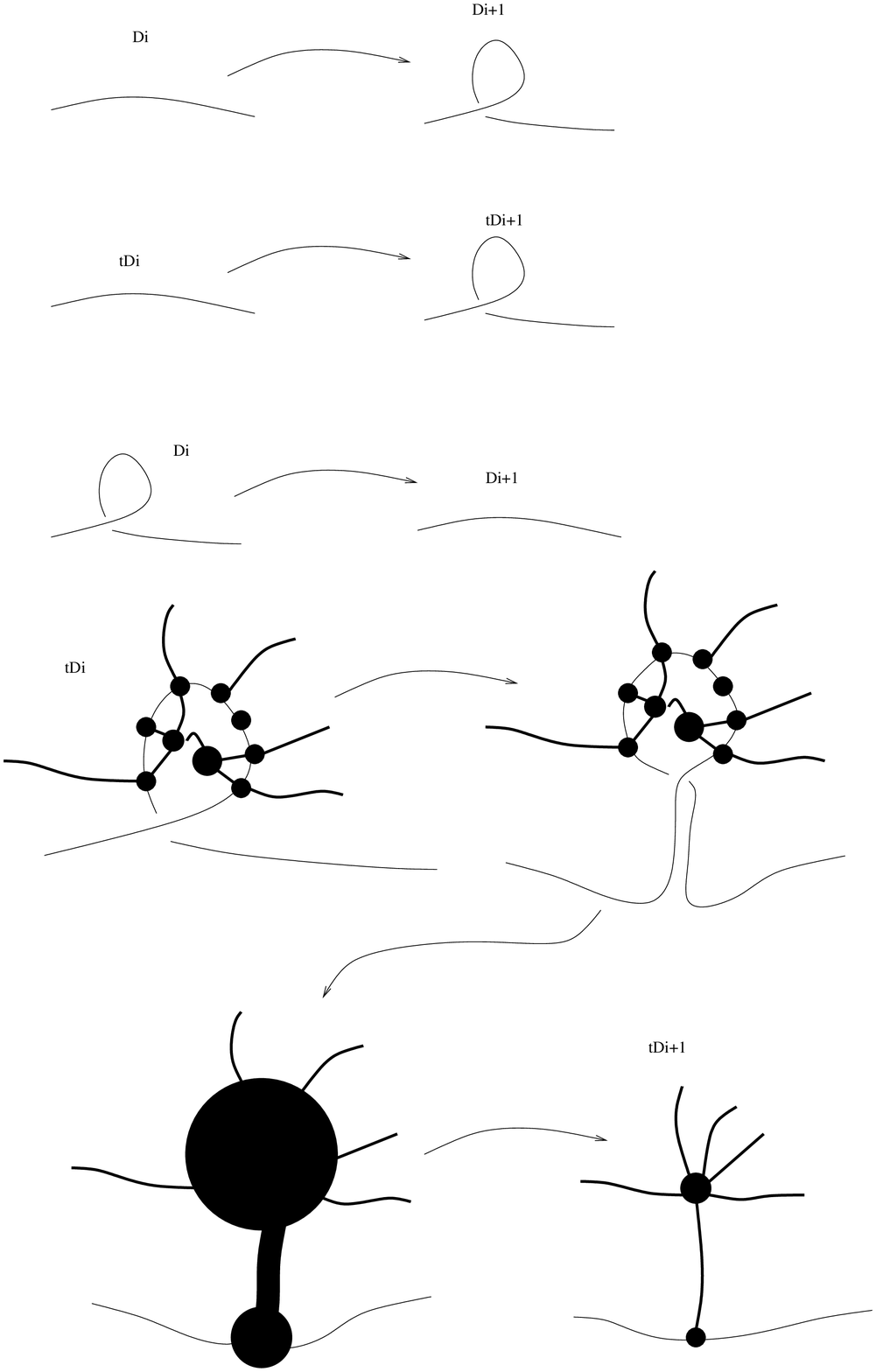}
}}
\caption{Reidemeister I}
\label{spaghettiR1}
\end{figure}

\begin{figure}
\centerline{
{
\psfrag{Di}{$D_i$}
\psfrag{Di+1}{$D_{i+1}$}
\psfrag{tDi}{$\tilde{D_i}$}
\psfrag{tDi+1}{$\tilde{D_{i+1}}$}
\psfrag{0handle}{$0$-handle}
\psfrag{1handle}{$1$-handle}
\psfrag{2handle}{$2$-handle}
\psfrag{Rii}{Reidemeister II.2}
\psfrag{Rii.i}{Reidemeister II.1}
\includegraphics[height=7in,width=4.5in]{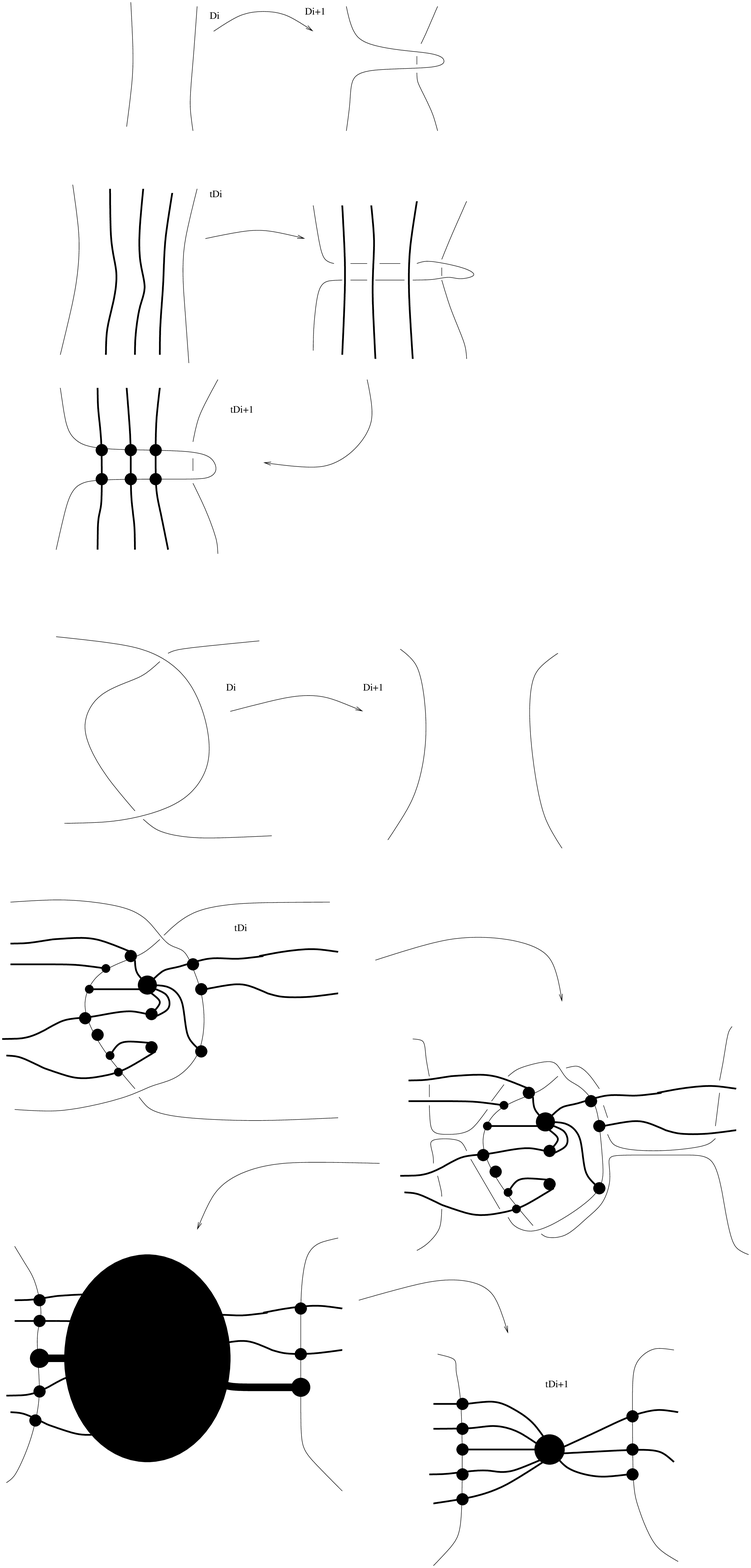}
}}
\caption{Reidemeister II}
\label{spaghettiR2}
\end{figure}

\begin{figure}
\centerline{
{
\psfrag{Di}{$D_i$}
\psfrag{Di+1}{$D_{i+1}$}
\psfrag{tDi}{$\tilde{D_i}$}
\psfrag{tDi+1}{$\tilde{D_{i+1}}$}
\psfrag{0handle}{$0$-handle}
\psfrag{1handle}{$1$-handle}
\psfrag{2handle}{$2$-handle}
\psfrag{Rii}{Reidemeister II.2}
\psfrag{Rii.i}{Reidemeister II.1}
\includegraphics[height=6.5in,width=4.5in]{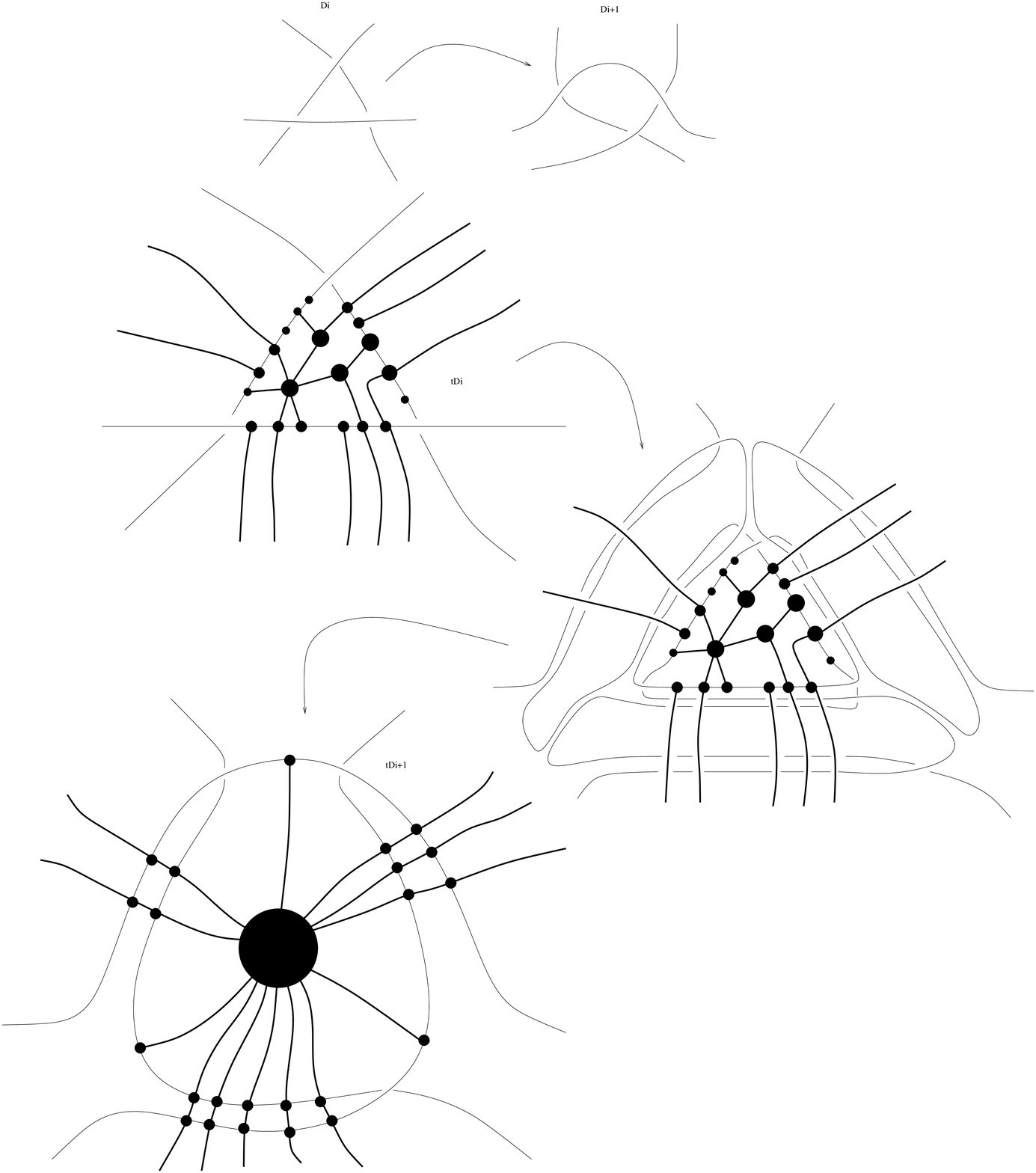}
}}
\caption{Reidemeister III}
\label{spaghettiR3}
\end{figure}

\begin{figure}
\centerline{
{
\psfrag{Di}{$D_i$}
\psfrag{Di+1}{$D_{i+1}$}
\psfrag{tDi}{$\tilde{D_i}$}
\psfrag{tDi+1}{$\tilde{D_{i+1}}$}
\psfrag{0handle}{$0$-handle}
\psfrag{1handle}{$1$-handle}
\psfrag{2handle}{$2$-handle}
\psfrag{Rii}{Reidemeister II.2}
\psfrag{Rii.i}{Reidemeister II.1}
\includegraphics[height=2.9in,width=4in]{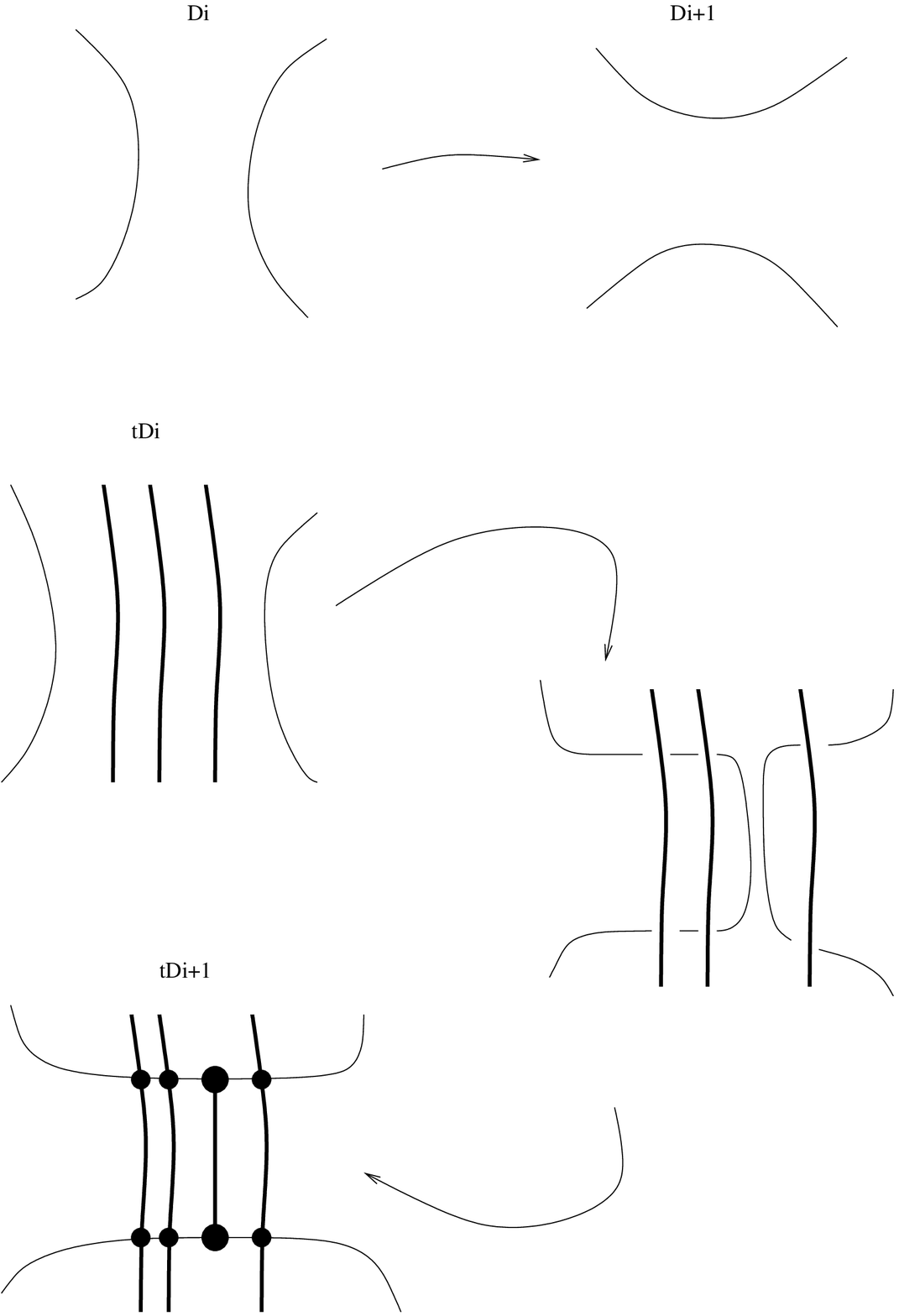}
}}
\caption{$1$-handle}
\label{1handlespaghetti}
\end{figure}

We shall specify a presentation of a surface which starts with the diagram $\tilde{D_1} = D_1$, and to get from diagram $\tilde{D_i}$ to $\tilde{D_{i+1}}$ we perform a sequence of Reidemeister II moves (each of which introduces two new crossings) or a Reidemeister I move (introducing a crossing) or no moves.  Note that we never add any handles in moving from $\tilde{D_i}$ to $\tilde{D_{i+1}}$.  Each $\tilde{D_i}$ is a link diagram (of $h+1$ components), but we shall find it convenient to draw $\tilde{D_i}$ as a diagram consisting of normal edges (\emph{thin edges}), thick edges (\emph{pipes}), and \emph{balls}.  In Figure \ref{localbestiary} we show what each $\tilde{D_i}$ can look like locally.   Each pipe represents two uncrossed strands, each ball with three or more incident pipes represents some arbitrary tangle.  If we delete all pipes and balls from the diagram $\tilde{D_i}$ we get the diagram $D_i$.   Now we give some diagrams showing with what it is that we replace each elementary cobordism of the presentation $S$ (Figures \ref{spaghettiR1}, \ref{spaghettiR2}, \ref{spaghettiR3}, \ref{1handlespaghetti}).

%%{spaghettiR1} {R2} {R3}

Now that we have reached diagram $\tilde{D_d}$, the next thing that we want to do in our surface presentation is to add some $1$-handles, joining the strands together in certain pipes.  Certain of the cobordisms described already have amalgamated previous pipes and balls into bigger balls.  If we imagine that each ball is translucent, we can see all the pipes that we have introduced.  We add the following $1$-handles:

\begin{itemize}

\item In one of the Reidemeister I moves (Figure \ref{spaghettiR1}), a new pipe is introduced.  Add a $1$-handle joining the two strands in this pipe.

\item In one of the Reidemeister II moves (Figure \ref{spaghettiR2}), two new pipes are introduced.  Add two $1$-handles, one in each pipe, joining the two strands together.

\item In the Reidemeister III move (Figure \ref{spaghettiR3}), three new pipes are introduced.  Add three $1$-handles, one in each pipe, joining the two strands together.

\item In the $1$-handle move (Figure \ref{1handlespaghetti}), a new pipe is introduced.  Add a $1$-handle joining the two strands in this pipe.

\end{itemize}

The only condition that we impose for the order of these $1$-handle additions is that we want to add the $h$ chosen $1$-handles that fuse together the $h+1$ components of $\tilde{D_d}$ before we add the rest of the $1$-handles.  After doing all the $1$-handle additions we arrive at diagram $\tilde{D_{d+1}}$, in which some of the pipes have had $1$-handles added.  The idea now is to isotope the strands in each pipe away down the pipe, moving away from the site of the $1$-handle addition.  We are going to do this in the reverse order to that in which the pipes appeared in the cobordism.

\begin{figure}
\centerline{
{
\psfrag{Ei}{$E_i$}
\psfrag{Ei+1}{$E_{i+1}$}
\psfrag{tDi}{$\tilde{D_i}$}
\psfrag{tDi+1}{$\tilde{D_{i+1}}$}
\psfrag{0handle}{$0$-handle}
\psfrag{1handle}{Sites of $1$-handle addition}
\psfrag{2handle}{$2$-handle}
\psfrag{Rii}{Reidemeister II.2}
\psfrag{Rii.i}{Reidemeister II.1}
\includegraphics[height=5in,width=5in]{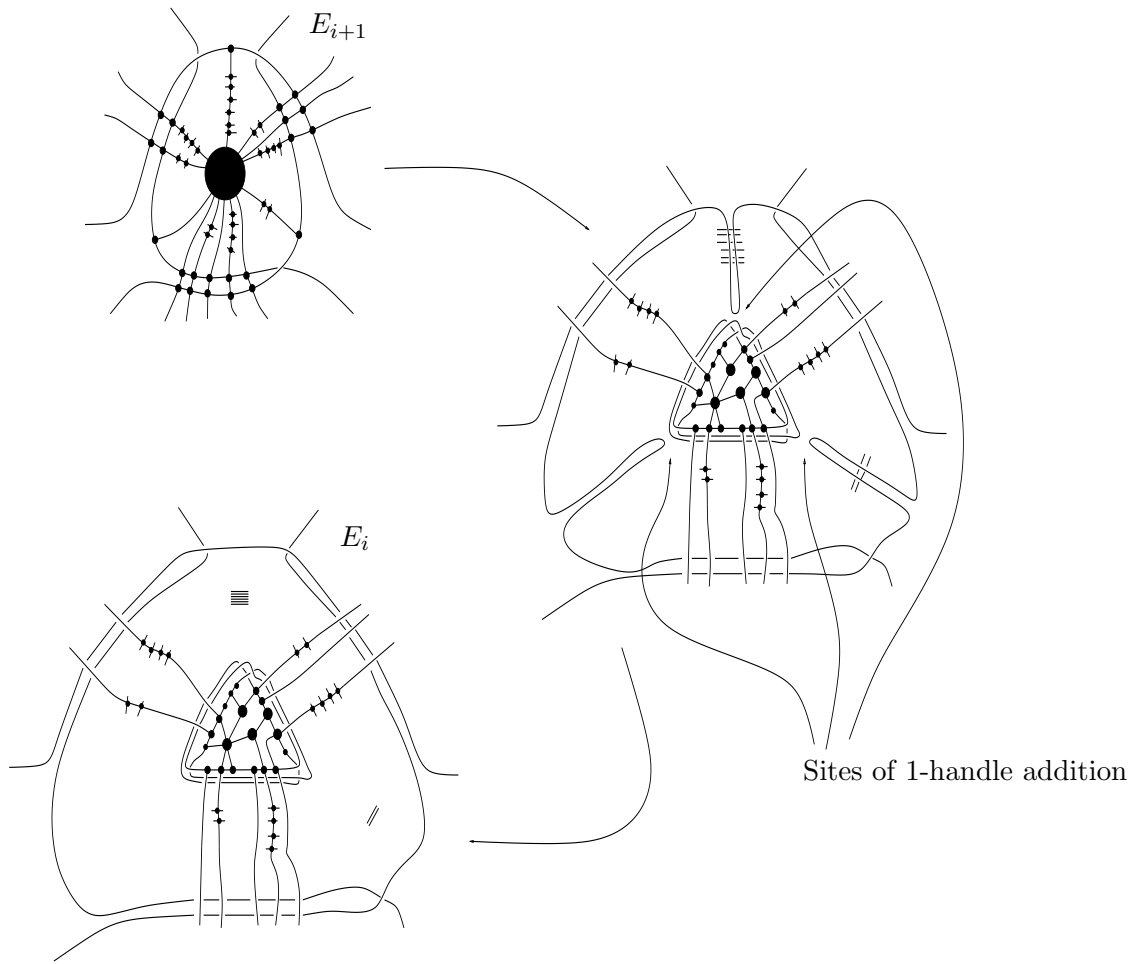}
}}
\caption{Retracting pipes}
\label{suckback}
\end{figure}

We shall write $E_d = \tilde{D_{d+1}}$.  Our cobordism will now pass through intermediary diagrams $E_d, E_{d-1}, \ldots , E_1$.  We describe the cobordism that takes us from $E_{i+1}$ to $E_i$:

\begin{itemize}

\item If $D_i$ was taken to $D_{i+1}$ by a Reidemeister move I or II that increased the total number of crossings then $E_i = E_{i+1}$ with the trivial cobordism.

\item If $D_i$ was taken to $D_{i+1}$ by a Reidemeister move I or II that decreased the total number of crossings, by a Reidemeister III move, or by a $1$-handle addition, then $E_{i+1}$ is taken to $E_i$ by a number of Reidemeister II moves, each reducing the total number of crossings, as illustrated in Figure \ref{suckback}. 

\end{itemize}

Although only the Reidemeister III case is shown in Figure \ref{suckback}, the reader should be able to deduce what is intended for the other cases.  The effect intended is that associated to trying eat a plate of spaghetti by putting a few strands in one's mouth and sucking.  The extra thin edge strands crossing the three pipes in Figure \ref{suckback} which are not present in Figure \ref{spaghettiR3} represent the possibility of a Reidemeister II or III move or $1$-handle $D_j \rightarrow D_{j+1}$ for $i < j < d$ introducing new balls in $\tilde{D_{j+1}}$ where a thin edge has to cross some pipes as in the intermediate stages of Figures \ref{spaghettiR2} and \ref{spaghettiR3}.

\begin{figure}
\centerline{
{
\psfrag{Ei}{$E_i$}
\psfrag{Ei+1}{$E_{i+1}$}
\psfrag{phi}{$\phi$}
\psfrag{tDi+1}{$\tilde{D_{i+1}}$}
\psfrag{0handle}{$0$-handle}
\psfrag{1handle}{Sites of $1$-handle addition}
\psfrag{2handle}{$2$-handle}
\psfrag{Rii}{Reidemeister II.2}
\psfrag{Rii.i}{Reidemeister II.1}
\includegraphics[height=3in,width=3in]{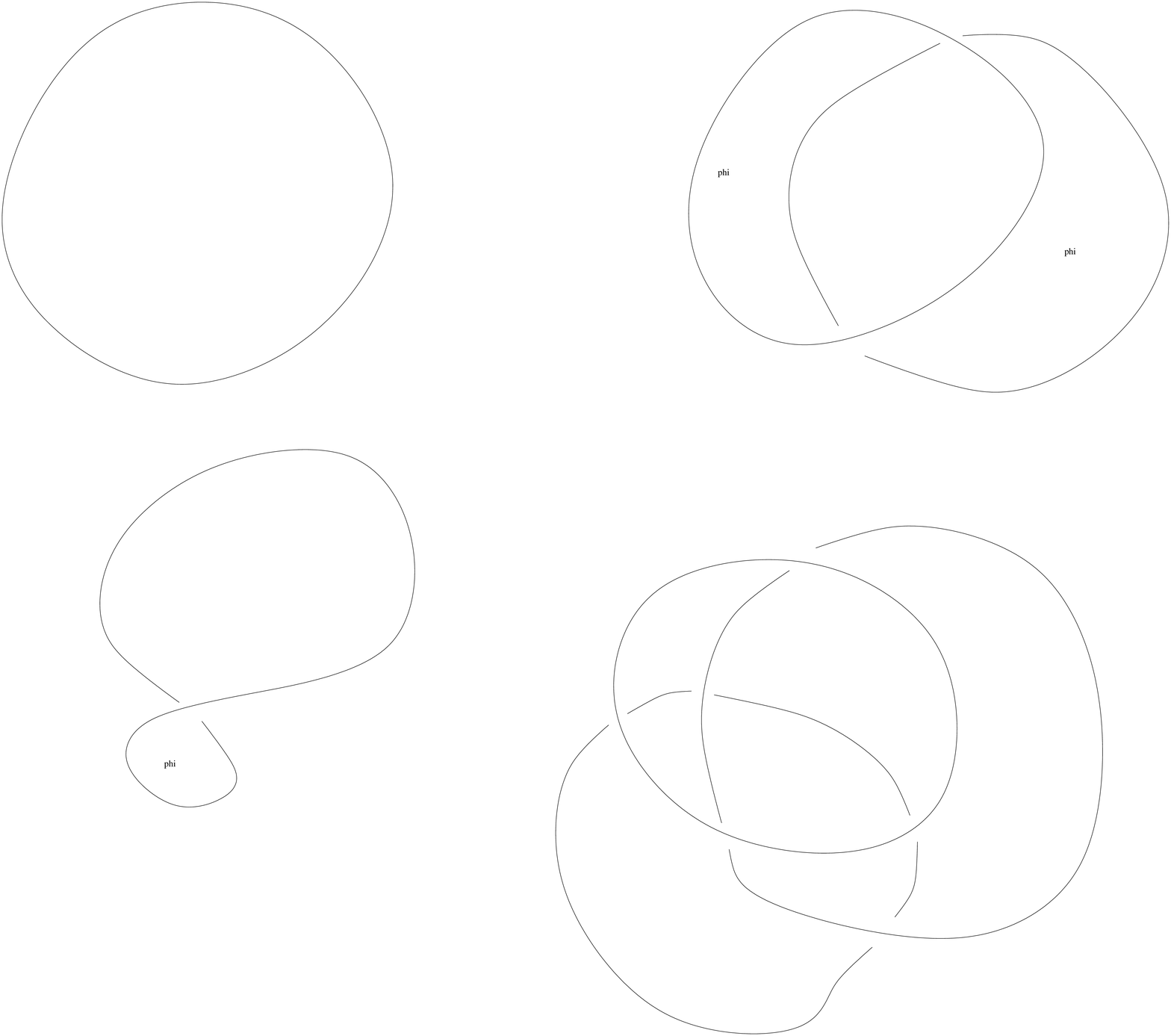}
}}
\caption{Components of $E_1$}
\label{complicatedunlink}
\end{figure}

We arrive at diagram $E_1$ which is now free of pipes and balls.  The components of $E_1$ are as shown in Figure \ref{complicatedunlink}, the areas of the components labelled by $\phi$ are, by construction, empty of other components.  By Reidemeister I and II moves we arrive at diagram $E_0$, a diagram of the unlink consisting of $0$-crossing unknots and $3$-component unlinks as in the bottom right of Figure \ref{complicatedunlink}.

It can be shown that if we perform Reidemeister moves to convert our diagram to the $0$-crossing diagram of the unlink and then add $2$-handles to cap off all but one of the components then what we have described is a new presentation $S'$ of the surface $\Sigma$.  However, we do not strictly need to know this for the proof of our result; but the reader should convince herself that $S'$ is a presentation of a surface of the same genus as $\Sigma$ (by counting $1$-handles and the number of components of our final unlink).

\subsection{The image of a generator}

Associated to the presentation $S'$ of the previous section is a map

\[ HKh_w(S') : HKh_w(D) \rightarrow HKh_w(E_0) \]

\noindent where $E_0$ is the final diagram in the presentation $S'$ and is an unlink of, say, $e_0$ components.

By the work of the previous sections we know that $HKh_w(S')$ is injective.  If $\alpha_\xi \in HKh_w(D)$ is the generator associated to decorating $K$ with the root $\xi$ of $\partial w$, then we know that

\[ HKh_w(S')(\alpha_\xi) = \alpha_\xi \in HKh_w(E_0) \]

\noindent where we write $\alpha_\xi$ also for the element of $HKh_w(E_0)$ associated to decorating each component of $E_0$ with $\xi$ (here is where it is important to remember that we added the $h$ ``fusing'' $1$-handles first.

\begin{figure}
\centerline{
{
\psfrag{B}{$B$}
\psfrag{Bo}{$B^o$}
\psfrag{phi}{$\phi$}
\psfrag{tDi+1}{$\tilde{D_{i+1}}$}
\psfrag{0handle}{$0$-handle}
\psfrag{1handle}{Sites of $1$-handle addition}
\psfrag{2handle}{$2$-handle}
\psfrag{Rii}{Reidemeister II.2}
\psfrag{Rii.i}{Reidemeister II.1}
\includegraphics[height=2in,width=4in]{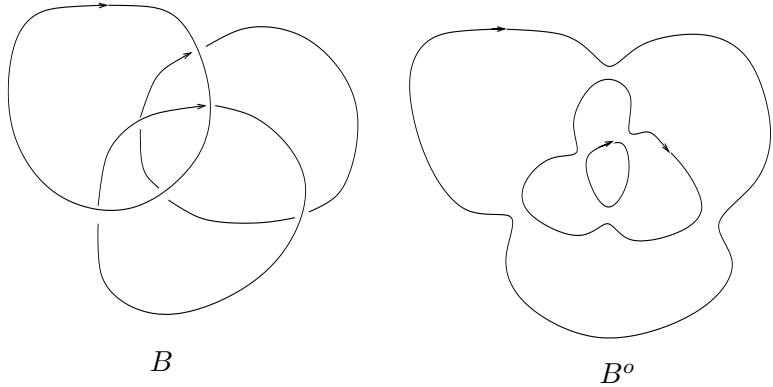}
}}
\caption{Faux Borromean Rings}
\label{orientedFBR}
\end{figure}

We would like know in which graded degree of the grading associated to the quantum filtration of $HKh_w(E_0)$ $\alpha_\xi$ lies.  Write $B$ for the \emph{fauxs} Borromean rings unlink that appears in the bottom right of Figure \ref{complicatedunlink} and in Figure \ref{orientedFBR}.  First we would like to compute in which degree $\alpha_\xi$ lies in $HKh_w(B)$.  Figure \ref{orientedFBR} shows $B$ and the oriented resolution $B^o$ of $B$.  The chain representative of $\alpha_\xi$ lies in $H(B^o)$ which is a summand of the $0$th homological degree chain group.

The chain representative of $\alpha_\xi$ in $H(B^o)$ is, up to non-zero scalar multiplication, the module element

\[ \prod_{i=1,2,3} \frac{ \partial w_i }{ x_i - \xi } \in (\mathbb{C}[x_1,x_2,x_3]/(\partial w_1, \partial w_2, \partial w_3)) \{ 3-3n \} \]

\noindent which is of top filtration grading $3n-3$.  \emph{A priori}, of course, it does not follow that $[ \alpha_\xi ] \in HKh_w^{0, 3n-3}$, since we have not yet seen that $\alpha_\xi$ is not homologous to an element of $CKh_w$ which is of a lower filtered degree.  We shall show this now.

The diagram $B$ is a diagram of the unlink, so for dimensional reasons the spectral sequence converging from $E_2 = HKh_n^{i,j}$ to $E_\infty = HKh_w^{i,j}$ collapses immediately ($E_2 = E_\infty$).

The reduction $\tilde{\alpha}_\xi \in \mathcal{F}^{3n-3}H(B^o) / \mathcal{F}^{3n-4}H(B^o)$ is a cycle in the page $E_1 = CKh_n^{i,j}$.  If we can show that $\tilde{\alpha}_\xi$ is not a boundary with respect to the $E_1$ differential (which is just the standard Khovanov-Rozansky $sl(n)$ differential) then, since $E_2 = E_\infty$, $\tilde{\alpha}_\xi$ will represent a non-zero class on the $E_\infty$ page.  The grading of this class (which is necessarily $j=3n-3$) will be the grading in which $[\alpha_\xi]$ lies in $HKh_w^{0,j}$.

So, in order to show that $ [ \alpha_\xi ] \in HKh_w^{0,3n-3} $, we need to see that $\tilde{\alpha}_\xi$ represents a non-zero class in $HKh_n$.  Khovanov and Rozansky \cite{KhovRoz1} have provided quantum-degree $0$ chain homotopy equivalences between the chain complexes $CKh_n$ corresponding to tangle diagrams differing by a single Reidemeister move.  We shall change $B$ by certain of these Reidemeister moves to arrive at the $0$-crossing $3$-component unlink.  The chain maps of Khovanov-Rozansky's chain homotopy equivalences will act on $\tilde{\alpha}_\xi$, mapping it to a cycle of the chain complex corresponding to the $0$-crossing $3$-component unlink.  This cycle represents a non-zero class in the homology, so $\tilde{\alpha}_\xi$ represents a non-zero class.  This will complete our argument that $[ \alpha_\xi ] \in HKh_w^{0,3n-3}$.

\begin{figure}
\centerline{
{
\psfrag{B}{$B$}
\psfrag{R3}{Reidemeister III move}
\psfrag{phi}{$\phi$}
\psfrag{tDi+1}{$\tilde{D_{i+1}}$}
\psfrag{0handle}{$0$-handle}
\psfrag{1handle}{Sites of $1$-handle addition}
\psfrag{2handle}{$2$-handle}
\psfrag{R2}{Reidemeister II.1}
\psfrag{Rii.i}{Reidemeister II.1}
\includegraphics[height=3in,width=3.5in]{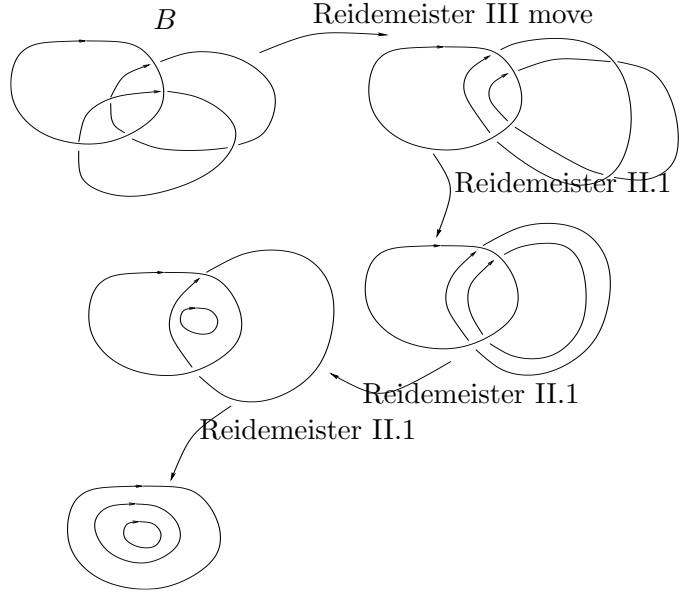}
}}
\caption{Uncrossing the unlink}
\label{FabstBlueRibbon}
\end{figure}

In Figure \ref{FabstBlueRibbon} we have shown the sequence of Reidemeister moves that we use to convert $B$ to the $0$-crossing $3$-component unlink.  Note that all the Reidemeister II moves are II.1 moves, that the Reidemeister III move involves three positive crossings, and that the number of circles in the oriented resolution never changes.

The cycle $\tilde{\alpha}_\xi$ in $CKh_n$ lies in the $0$th chain group summand $H(B^o)$.  Up to non-zero scalar multiplication:

\[ \tilde{\alpha}_\xi = x_1^{n-1}x_2^{n-1}x_3^{n-1} \in (\mathbb{C}[x_1, x_2, x_3]/x_1^nx_2^nx_3^n)\{ 3-3n \} = H(B^o) \]

Unpacking the proofs a little in \cite{KhovRoz1}, we see that (for the Reidemeister II.1 move which reduces the number of crossings and for the Reidemeister III move with all positive crossings) the maps induced from the chain group summands corresponding to the oriented resolutions are just the identity map of matrix factorizations.  (The fact that, in section 4.2.1 of this paper, the component of this section's chain map that runs from the oriented resolution of the two crossed strands to the uncrossed strands is the identity map, is a reflection of this fact).

This means that the element $\tilde{\alpha}_\xi$ always gets mapped to the top-degree element in the chain group summand that is the homology of the oriented resolution.  In the final diagram of the $0$-crossing $3$-component unlink, the homology of the oriented resolution is the \emph{only} chain group summand.  So the image of $\tilde{\alpha}_\xi$ represents a non-zero element in homology and, since the chain maps are all chain homotopy equivalences, $\tilde{\alpha}_\xi$ also represents a non-zero element in homology.

\subsection{Putting it all together}

We have shown that given a presentation of a connected cobordism $\Sigma$ between a knot $K$ and the unknot (starting with diagram $D$ for the knot, and ending with the $0$-crossing unknot diagram) we can convert this into a presentation of a surface with the same genus ending in a diagram of an $h$-component unlink $U^h$.  Furthermore, this presentation induces a map that is filtered of degree $(1-n)\chi$ (where $\chi$ is the Euler characteristic of the new multiply-punctured surface) from the homology $HKh_w(D)$ of $D$ to the homology of the $h$-component unlink.

If $\alpha_\xi$ is the element of $HKh_w(D)$ which is given by decorating $K$ with the root $\xi$ of $\partial w$, then we have seen that this element gets mapped to the top-degree part of $HKh_w(U^h)$, i.e. to a non-zero element of $HKh_w^{0,(n-1)h}$.  Since the map was filtered of degree $(1-n)\chi$, we see that $\alpha_\xi$ must have been in $HKh_w^{0,j}(D)$ where $j \geq (1-n)(2g(\Sigma) - 1)$.  This completes the proof of our main theorem \ref{maintheorem}.

\subsection{Computation for positive knots}

In this subsection we deduce Milnor's conjecture on the slice genus of torus knots as a corollary to our main theorem.

Suppose that a knot $K$ has a diagram $D$ in which all the crossings are positive (i.e. look like the crossing on the left of Figure \ref{writhe}).  Suppose that $D$ has $k$ crossings and $l$ circles in its oriented resolution.  We compute the top grading of $HKh_w(D)$.

Since all the crossings are positive, the homology of the oriented resolution is the leftmost chain group in the Khovanov cube (also the chain group of homological degree $0$).  Hence the homology $HKh_w(D)$ will just be the kernel of the differential coming from this chain group, and in particular the filtration grading of chain representatives of the homology will be the same taken in the chain group as taken in the homology since there is no boundary group by which to quotient.

Each basis element (got by decorating $K$ with a root $\xi$ of $\partial w$) of $HKh_w(D)$ lies in $HKh_w^{0,(1-n)(k-l)}(D)$.  In particular, $HKh_w^{0,(1-n)(k-l)}(D) \not= 0$.  Now, $(1-n)(k-l)$ is the highest filtration degree of the $0$th chain group, so our theorem says that

\[ (n-1)(2g^*(K)-1) \geq (n-1)(k-l) \]

There exists an explicit description of a Seifert surface for $K$ which consists of $l$ disks (filling the circles of the oriented resolution of $D$), connected by $k$ bands (where the crossings of $D$ are).  The Euler characteristic of this surface with boundary is $l-k$, so the surface is of genus $(k-l+1)/2$.  Pushing this Seifert surface slightly into the $4$-ball yields a slice surface for $K$ of the same genus.  Since we have seen that

\[ 2g^*(K) \geq k-1+1 \]

\noindent it follows that $g^*(K) = k-l+1$.  Performing this computation in the case of the standard diagram of a torus knot yields Milnor's conjecture on the slice genus of torus knots.

\section{Appendix}

Here we provide a proof of basic results used (sometimes implicitly) throughout this paper.

\subsection{Removal of Marks}
\label{removal of marks}

In this section we omit mention of the various quantum grading shifts that occur, although readers friendly with matrix factorizations can easily assure themselves that the filtration gradings of each map is as expected.

To first explain Figure \ref{markremoval}: the matrix factorization $N$ (resp. $N'$) is obtained from $Q$ (resp. $Q'$) by tensoring with the matrix factorization $P$.  Also, the matrix factorization $M$ (resp. $M'$) is the same as $Q$ (resp. $Q'$) with the formal replacement of the variable $y$ by $x$.

\begin{theorem} Removal of marks
\label{losingmarks}

We prove that the matrix factorizations $M$ (resp. $M'$) and $N$ (resp. $N'$) (in which we intend the circles to contain the same arbitrary trivalent graph with thick edges) in Figure \ref{markremoval} are equal in the homotopy category of matrix factorizations.

Furthermore, suppose $M'$, $N'$, $Q'$ are factorizations, maybe different from $M$, $N$, $Q$ but with the same boundary labels.  If we then have a map of matrix factorizations $\alpha_x : M \rightarrow M'$ (and, replacing $x$ by $y$, $\alpha_y : Q \rightarrow Q'$), inducing a map of matrix factorizations $A : N \rightarrow N'$, then

\[ \alpha_x = L \circ A \circ K \]

\noindent where $L: N \rightarrow M$ and $K: M \rightarrow N$ are the homotopy equivalences constructed.

\end{theorem}

\begin{proof}

We exhibit filtered degree-$0$ maps $M \rightarrow N$, $N \rightarrow M$ which we show are homotopy equivalences.

The matrix factorizations $M$ and $N$ are both defined over the same ground ring $R = \mathbb{C}[x_1,x_2,...,x_r,x]$.  Let's write $M$ as

\[ M^0 \stackrel{f_x}{\rightarrow} M^1 \stackrel{g_x}{\rightarrow} M^0 \]

\noindent $Q$ is the same factorization but with any occurence of $x$ relabelled as $y$, we shall consequently write $Q$ as

\[ Q^0 \stackrel{f_y}{\rightarrow} Q^1 \stackrel{g_y}{\rightarrow} Q^0 \]

\noindent where the subscripts are to remind us that the difference between the factorizations $M$ and $Q$ is just the interchanging of $x$ with $y$.

\begin{figure}
\centerline{
{
\psfrag{x}{$x$}
\psfrag{y}{$y$}
\psfrag{M}{$M$}
\psfrag{N}{$N$}
\psfrag{P}{$P$}
\psfrag{Q}{$Q$}
\psfrag{M'}{$M'$}
\psfrag{N'}{$N'$}
\psfrag{P'}{$P'$}
\psfrag{Q'}{$Q'$}
\psfrag{x1}{$x_1$}
\psfrag{xr}{$x_r$}
\includegraphics[height=3in,width=4in]{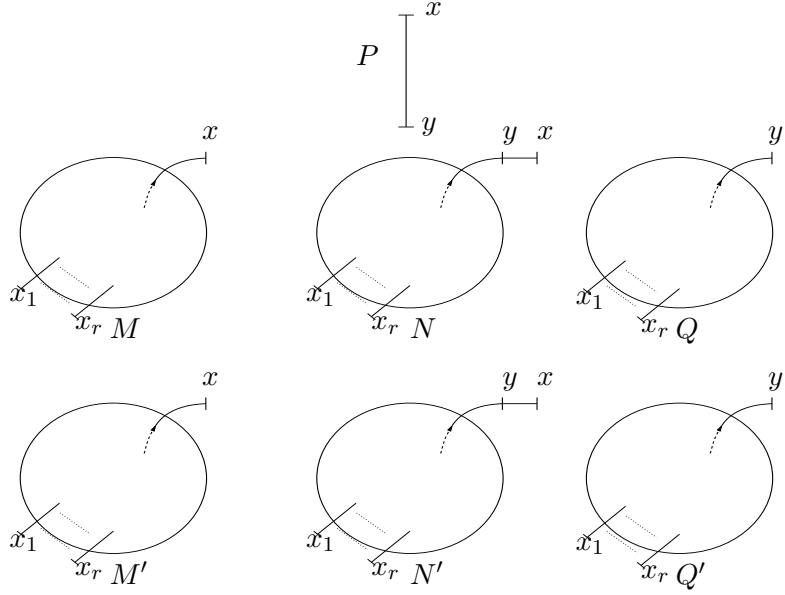}
}}
\caption{Removal of marks}
\label{markremoval}
\end{figure}

%\begin{figure}
%\centerline{
%{
%\psfrag{x}{$x$}
%\psfrag{y}{$y$}
%\psfrag{M}{$M,M'$}
%\psfrag{N}{$N,N'$}
%\psfrag{P}{$P,P'$}
%\psfrag{Q}{$Q,Q'$}
%\includegraphics[height=2in,width=3.2in]{markremoval2.eps}
%}}
%\label{markremoval2}
%\end{figure}

Then if $P$ is the factorization

\[ P^0 \stackrel{\pi_{xy}}{\rightarrow} P^1 \stackrel{x-y}{\rightarrow} P^0\]

\noindent where $P^0$ and $P^1$ are free rank-$1$ modules over $\mathbb{C}[x,y]$ we have that $N^0 = (Q^0 \otimes P^0) \oplus (Q^1 \otimes P^1)$, $N^1 = (Q^1 \otimes P^0) \oplus (Q^0 \otimes P^1)$ and $N$ can be written

\[ N^0 \stackrel{\left(\begin{array}{cc}
f_y & -(x-y) \\
\pi_{xy} & g_y \end{array} \right)}{\rightarrow} N^1 \stackrel{\left( \begin{array}{cc}
g_y & x-y \\
-\pi_{xy} & f_y \end{array} \right)}{\rightarrow} N^0 \]

\noindent in which we are implicitly thinking of $N$ as a factorization over the polynomial ring $\mathbb{C}[x_1,x_2,...,x_r,x,y]$ with degenerate potential.

Now we give the homotopy equivalence between $M$ and $N$.  We define maps of matrix factorizations $K: M \rightarrow N$ and $L: N \rightarrow M$ by the following $\mathbb{C}[x_1,x_2,...,x_r,x]$-module maps:

\[ K_0 = \left( \begin{array}{c} 1 \\ -\frac{f_x - f_y}{x-y} \end{array} \right), K_1 = \left( \begin{array}{c} 1 \\ \frac{g_x - g_y}{x-y} \end{array} \right) \]

\[ L_0 = ( \begin{array}{cc} P_{y \mapsto x} & 0 \end{array} ),  L_1 = ( \begin{array}{cc} P_{y \mapsto x} & 0 \end{array} ) \]

\noindent where by $P_{y \mapsto x}$ we mean the map which looks like the ring map

\[ \mathbb{C}[x_1,x_2,...,x_r,x,y] \rightarrow  \mathbb{C}[x_1,x_2,...,x_r,x,y]/(x-y) =  \mathbb{C}[x_1,x_2,...,x_r,x] \]

\noindent  on each module summand $\mathbb{C}[x_1,x_2,...,x_r,x,y]$.

 Note that $ L \circ K : M \rightarrow M $ is already the identity map.

We define an homotopy $H_0 : N^0 \rightarrow N^1, H_1 : N^1 \rightarrow N^0$ as follows:

\[ H_0 = \left( \begin{array}{cc} 0 & 0 \\ \frac{1 - P_{y \mapsto x}}{x-y} & 0 \end{array} \right), H_1 = \left( \begin{array}{cc} 0 & 0 \\ \frac{P_{y \mapsto x} - 1}{x-y} & 0 \end{array} \right) \]

It is now an exercise to see that

\[ id_{N^0} - K_0 \circ L_0 = H_1 \circ \left(\begin{array}{cc}
f_y & -(x-y) \\
\pi_{xy} & g_y \end{array} \right) + \left( \begin{array}{cc}
g_y & x-y \\
-\pi_{xy} & f_y \end{array} \right) \circ H_0 \]

\noindent and

\[ id_{N^1} - K_1 \circ L_1 = \left(\begin{array}{cc}
f_y & -(x-y) \\
\pi_{xy} & g_y \end{array} \right) \circ H_1 + H_0 \circ \left( \begin{array}{cc}
g_y & x-y \\
-\pi_{xy} & f_y \end{array} \right) \]

The second part of the theorem amounts to the observation that if $\alpha_{y,0}$ and $\alpha_{y,1}$  are the components of $\alpha_y : Q \rightarrow Q' $ then the induced map $A : N \rightarrow N'$ has components

\[ A_0 = \left( \begin{array}{cc} \alpha_{y,0} & 0 \\ 0 & \alpha_{y,1} \end{array} \right) , A_1 = \left( \begin{array}{cc} \alpha_{y,1} & 0 \\ 0 & \alpha_{y,0} \end{array} \right) \]

\noindent and then it is immediate that

\[ \alpha_{x,0} =  \left( \begin{array}{cc} P_{y \rightarrow x} & 0 \end{array} \right) \left( \begin{array}{cc} \alpha_{y,0} & 0 \\ 0 & \alpha_{y,1} \end{array} \right) \left( \begin{array}{c} 1 \\ -\frac{ f_x - f_y }{x - y} \end{array} \right) \]

\noindent and

\[ \alpha_{x,1} =  \left( \begin{array}{cc} P_{y \rightarrow x} & 0 \end{array} \right) \left( \begin{array}{cc} \alpha_{y,1} & 0 \\ 0 & \alpha_{y,0} \end{array} \right) \left( \begin{array}{c} 1 \\ \frac{g_x - g_y}{x-y} \end{array} \right) \]

\end{proof}

\subsection{The $\eta$ map}
\label{the eta map}

\begin{figure}
\centerline{
{
\psfrag{1}{$1$}
\psfrag{2}{$2$}
\psfrag{3}{$3$}
\psfrag{4}{$4$}
\psfrag{=}{$=$}
\psfrag{eta0}{$\eta_0$}
\psfrag{eta1}{$\eta_1$}
\psfrag{chi0}{$\chi_0$}
\psfrag{chi1}{$\chi_1$}
\includegraphics[height=2in,width=3.4in]{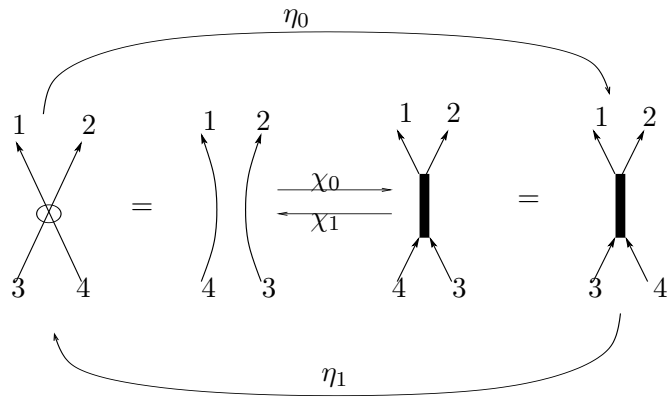}
}}
\caption{Definition of the $\eta$ maps}
\label{gorniketa}
\end{figure}

The $\eta_0$ and $\eta_1$ maps from Section 2 are defined as in Figure \ref{gorniketa}.

The equalities in Figure \ref{gorniketa} are equalities as bare matrix factorizations.  It is useful for us to know how this $\eta$ map would be different if it were defined as conjugation of the $\chi$ map by swapping $1$ and $2$ instead of $3$ and $4$.  For example, at the end of our proof that the generators of $HKh_w$ are preserved under the Reidemeister $2.1$ move, we implicitly use the fact that the $\eta$ map would merely change up to sign.  It is a justification of this that we give here.

\begin{figure}
\centerline{
{
\psfrag{1}{$1$}
\psfrag{2}{$2$}
\psfrag{3}{$3$}
\psfrag{4}{$4$}
\psfrag{=}{$=$}
\psfrag{chi0}{$\chi_0$}
\psfrag{chi1}{$\chi_1$}
\psfrag{-chi0}{$-\chi_0$}
\psfrag{-chi1}{$-\chi_1$}
\includegraphics[height=2in,width=3.4in]{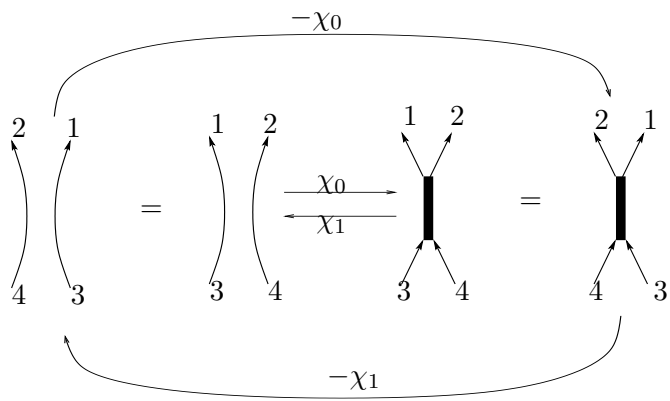}
}}
\caption{A commutative diagram}
\label{mickeymouse}
\end{figure}

We can restate this result in terms of Figure \ref{mickeymouse} which we assert to be commutative, as is easy to verify.

%%\subsection{The product of the $\chi$ maps}
%%\label{product of chi maps}

%%In fig \ref{chimaps} we show the $\chi_0$ and $\chi_1$ maps.  Here we show that $\chi_0 \chi_1$ and $\chi_1 \chi_0$ are homotopic to multiplication by $x_2 - x_3$.

\bibliography{slnras}
\bibliographystyle{hamsplain}
\nocite{*}

\end{document}